%% file: tisp.tex
\ifx\fltex\undefined
\documentclass[12pt,a4paper,english]{article}
\usepackage{amsmath,amssymb}
\usepackage{babel}
\usepackage{eucal}
\fi % \ifx\fltex\undefined

\input{ macros.tex}

\begin{document}

\title{{\Large{\bf Invariant Subspaces for Operators \\ in a General II$_1$-factor}}}
\author{Uffe~Haagerup\, and Hanne~Schultz \footnote{Supported by The Danish National Research
      Foundation.}}
  \date{}
\maketitle

\begin{abstract}
\noindent It is shown that to every operator $T$ in a general von Neumann
      factor $\CM$ of type II$_1$ and to every Borel set $B$ in the
      complex plane $\C$, one can associate a largest, closed, $T$-invariant
      subspace, $\CK = \CK_T(B)$, affiliated with $\CM$, such that the Brown
      measure of $T|_\CK$ is concentrated on $B$. Moreover, $\CK$ is
      $T$-hyperinvariant, and the Brown measure of
      $P_{\CK^\bot}T|_{\CK^\bot}$ is concentrated on $\C\setminus B$. In
      particular, if $T\in\CM$ has a Brown measure which is
      not concentrated on a singleton, then there exists a non-trivial, closed,
      $T$-hyperinvariant subspace. Furthermore, it is shown that for every $T\in\CM$ the limit $A=\lim_{n\rightarrow\infty}[(T^n)\cc T^n]^{\frac{1}{2n}}$ exists in the strong operator topology and $\CK_T(\overline{B(0,r)})=1_{[0,r]}(A)$, $r>0$.
\end{abstract}

\input{tispintroduction.tex}
\input{differentialeq.tex}
\input{spectralsubspaces.tex}
\input{tispfreeprob.tex}
\input{Lipschitz.tex}

\input{integration.tex}

\input{invariant.tex}
\input{stronglimit.tex}
\input{decomposable.tex}

\input{appendix.tex}
\input{bibliography.tex}

\noindent Uffe~Haagerup\\
Department of Mathematics and Computer Science\\
University of Southern Denmark\\
Campusvej~55, 5230~Odense~M\\
Denmark\\
{\tt haagerup@imada.sdu.dk}

\vspace{.5cm}

\noindent Hanne~Schultz\\
Department of Mathematics and Computer Science\\
University of Southern Denmark\\
Campusvej~55, 5230~Odense~M\\
Denmark\\
{\tt schultz@imada.sdu.dk}

\end{document}

%% file: macros.tex
\DeclareMathVersion{scriptstyle}

\usepackage{isolatin1,theorem}
\theoremstyle{change}
{\theorembodyfont{\slshape} %Alle saetningsformuleringer saettes i 
  \newtheorem{thm}{Theorem.}[section]
  \newtheorem{lemma}[thm]{Lemma.}
  \newtheorem{prop}[thm]{Proposition.}
  \newtheorem{cor}[thm]{Corollary.}
  \newtheorem{mainthm}[thm]{Main Theorem.}
}
{\theorembodyfont{\rmfamily} %Eksempler og bemaerkninger saettes i roman
  \newtheorem{definition}[thm]{Definition.}
  \newtheorem{remark}[thm]{Remark.}
  \newtheorem{example}[thm]{Example.}

}

\numberwithin{equation}{section} % Formelnummerering
\newcommand{\proof}[1][Proof. ]{{\it#1}}

\def\endproof{{\nobreak\qquad{\scriptstyle \blacksquare}}}

\def\C{{\mathbb C}}
\def\B{{\mathbb B}}
\def\D{{\mathbb D}}

\def\F{{\mathbb F}}

\def\N{{\mathbb N}}
\def\R{{\mathbb R}}

\def\Q{{\mathbb Q}}

\def\CA{{\mathcal A}}
\def\CB{{\mathcal B}}
\def\CC{{\mathcal C}}
\def\CD{{\mathcal D}}
\def\CE{{\mathcal E}}

\def\CH{{\mathcal H}}

\def\CJ{{\mathcal J}}
\def\CK{{\mathcal K}}

\def\CM{{\mathcal M}}
\def\CN{{\mathcal N}}

\def\CR{{\mathcal R}}

\def\CT{{\mathcal T}}
\def\CU{{\mathcal U}}

\def\unit{{\bf 1}}
\def\i{{\rm i}}
\def\e{{\rm e}}
\def\d{{\rm d}}
\def\eps{\varepsilon}
\def\<{{\langle}}
\def\>{{\rangle}}

\def\CMD{\CM^\Delta}

\def\im{{\text{\rm Im}}}
\def\re{{\text{\rm Re}}}
\def\supp{{\text{\rm supp}}}

\def\Prob{{\text{\rm Prob}}}

\def\difft{\frac{\rm d}{{\rm d}t}}
\def\Lp{L^p(\CM, \tau)}
\def\LpN{L^p(\CN, \tau)}
\def\cc{^*}

%% file: tispintroduction.tex
\section{Introduction}

Consider a von Neumann algebra $\CM$ acting on the Hilbert space $\CH$. A
closed subspace $\CH_0$ of $\CH$ is said to be {\it affiliated with $\CM$}
if the projection of $\CH$ onto $\CH_0$ belongs to $\CM$. $\CH_0$ is
said to be {\it non-trivial} if $\CH_0\neq 0$ and $\CH_0\neq \CH$.
  For $T\in\CM$, a subspace $\CH_0$ of $\CH$ is said to be {\it
  $T$-invariant}, if $T(\CH_0)\subseteq \CH_0$, i.e. if $T$ and the projection
  $P_{\CH_0}$ onto $\overline{\CH_0}$ satisfy
  \[
  P_{\CH_0}TP_{\CH_0}=TP_{\CH_0}.
  \]
$\CH_0$ is said to be {\it hyperinvariant} for $T$ (or $T$-{\it hyperinvariant}) if it is $S$-invariant
for every $S\in B(\CH)$ satisfying $ST=TS$. It is not hard to see that if
$\CH_0$ is hyperinvariant for $T$, then $P_{\CH_0}\in
W\cc(T)=\{T\}''$. However, the converse statement does not hold true. In fact, one can find $A\in M_3(\C)$ and an
$A$-invariant projection $P\in W\cc(A)$ which is not $A$-hyperinvariant
(cf. \cite{D}).

  The {\it invariant subspace problem relative to $\CM$} asks whether every
operator $T\in \CM\setminus \C\unit$ has a non-trivial invariant subspace
$\CH_0$ affiliated with $\CM$, and the {\it hyperinvariant subspace
  problem} asks whether one can always choose such an $\CH_0$ to be
hyperinvariant for $T$. Of course, if $\CM$ is not a factor, then both of
these questions may be answered in the affirmative. Also, if $\CM$ is a factor of finite dimension,
i.e. $\CM \cong M_n(\C)$ for some $n\in \N$, then every operator $T\in
\CM\setminus \C\unit$ has a non-trivial eigenspace, and therefore $\C^n$ has a
non-trivial $T$-invariant subspace. In this paper we shall focus on the invariant subspace problem relative to
a II$_1$-factor.

Recall
from \cite[Section~8]{KR2} that every II$_1$-factor $\CM$ has a unique
tracial state $\tau$, and $\tau$ is faithful and normal. The {\it
  Fuglede-Kadison determinant}, $\Delta:\CM\rightarrow [0,\infty[$, is given by
\begin{equation}
\Delta(T)= \exp\{\tau(\log|T|)\}, \qquad (T\in\CM),
\end{equation}
with $\exp\{-\infty\}:=0$ (cf. \cite{FuKa}). Also recall from \cite{Br} that for fixed $T\in\CM$, the function
\[
\lambda\mapsto \log\Delta(T-\lambda\unit)
\]
is subharmonic in $\C$, and its Laplacian 
\begin{equation}
\d\mu_T(\lambda_1+\i\lambda_2):= \frac{1}{2\pi}\nabla^2 \log\Delta[T-(\lambda_1+\i\lambda_2)\unit]\,\d\lambda_1\,\d\lambda_2
\end{equation}
(taken in the distribution sense) defines a probability measure $\mu_T$ on
$\C$, the {\it Brown measure of $T$}, with $\supp(\mu_T)\subseteq \sigma(T)$. Note that if $T\in\CM$ is normal, then $\mu_T=\tau\circ E_T$, where $E_T:\B(\C)\rightarrow {\rm Proj}(\CM)$ is the projection valued measure on $(\C, \B(\C))$ in the spectral resolution of $T$:
\[
T=\int_{\sigma(T)}\lambda\,\d E_T(\lambda).
\]
  If $\CM=M_n(\C)$ for some $n\in\N$, then the Fuglede-Kadison determinant
and the Brown measure are also defined for $T\in\CM$, and in this case we
have that
\[
\Delta(T)=|{\rm det}T|^{\frac 1n},
\]
and
\[
\mu_T= \frac 1n \sum_{i=1}^n \delta_{\lambda_i},
\]
where $\lambda_1, \ldots, \lambda_n$ are the eigenvalues of $T$, repeated
according to multiplicity. 

\vspace{1cm}

The main result of this paper is (cf. Theorem~\ref{Borelsets}):

  \begin{mainthm}\label{embedding} Let $\CM$ be a II$_1$-factor. Then for every $T\in\CM$ and every Borelset
  $B\subseteq \C$ there is a largest closed, $T$-invariant subspace, $\CK =
  \CK_T(B)$, affiliated with $\CM$, such that the Brown measure of $T|_\CK$,
  $\mu_{T|_\CK}$, is concentrated on $B$.\footnote{If $\CK = \{0\}$, then
  we define $\mu_{T|_\CK}:=0$. If $\CK \neq \{0\}$,  then $\mu_{T|_\CK}$ is
  computed 
  relative to the II$_1$-factor $P\CM P$, where $P\in\CM$ denotes the
  projection onto $\CK$.} Moreover, $\CK$ is hyperinvariant for $T$, and if $P=P_T(B)\in\CM$ denotes the
  projection onto $\CK$, then
  \begin{itemize}
    \item[(i)] $\tau(P)=\mu_T(B)$,
    \item[(ii)] the Brown measure of $P^\bot TP^\bot$, considered as an
    element of $P^\bot \CM P^\bot$, is concentrated on $\C\setminus B$.
  \end{itemize}
  \end{mainthm}

\vspace{.2cm}

It is then easily seen that if $T\in \CM$, and if $\mu_T$ is not a Dirac measure, then $T$ has a
non-trivial hyperinvariant subspace (cf. Corollary~\ref{partial solution}).

A II$_1$-factor $\CM$ on a separable Hilbert space is said to have {\it the
  embedding property}, if it embeds in the ultrapower $\CR^\omega$ of the
  hyperfinite II$_1$-factor $\CR$ for some free ultrafilter $\omega$ on
  $\N$.  In 1976 Connes (cf. \cite{Co}) raised the question whether every
  II$_1$-factor on a separable Hilbert space has the embedding
  property. This problem remains unsolved.

    In his unpublished lecture notes \cite{H} from MSRI~2001 (see also \cite{H2}), the first
  author proved our Main Theorem in the special case where $\CM$ has the embedding property. The contents of sections~2 and 3 are by and large taken from \cite{H}, but for the rest of the paper, we resort to a completely different line of proof in order to treat the general case. Our proof is based on free probability theory (cf. sections~4 and 5) and the Turpin--Waelbroek method of integration in quasinormed spaces (cf. section~6). 

The construction of the spectral subspaces $\CK_T(B)$ referred to in the Main~Theorem is carried out in several steps. In section~3 we introduce the closed $T$--hyperinvariant subspaces $E(T,r)$ and $F(T,r)$ in the following way: $E(T,r)$ is defined as the set of vectors $\xi\in\CH$, for which there is a sequence
$(\xi_n)_{n=1}^\infty$ in $\CH$ such that
\[
\lim_{n\rightarrow \infty}\|\xi_n-\xi\|=0 \;\; {\rm and}
\;\;\limsup_{n\rightarrow \infty}\|(T-\lambda\unit)^n \xi_n\|^{\frac
  1n}\leq r,
\]
and $F(T,r)$ is defined as the set of vectors $\eta\in \CH$, for which there is a sequence
$(\eta_n)_{n=1}^\infty$ in $\CH$ such that
\[
\lim_{n\rightarrow \infty}\|(T-\lambda\unit)^n\eta_n-\eta\|=0 \;\; {\rm and}
\;\;\limsup_{n\rightarrow \infty}\| \eta_n\|^{\frac
  1n}\leq \frac1r.
\]
In section~7 we combine the results of sections~2 through 6 and prove that with
\[
\CK_T(\overline{B(0,r)}):= E(T,r), \qquad r>0,
\]
and
\[
\CK_T(\C\setminus B(0,r)):= F(T,r), \qquad r>0,
\]
$\CK_T(\overline{B(0,r)})$ ($\CK_T(\C\setminus B(0,r))$, resp.) satisfies the conditions listed in our Main~Theorem in the case $B=\overline{B(0,r)}$ ($B=\C\setminus B(0,r)$, resp.). 

Then for a closed subset $F$ of $\C$, $F\neq \C$, we write $\C\setminus F$
as a countable union of open balls:
\[
\C\setminus F = \bigcup_{n=1}^\infty B(\lambda_n, r_n),
\]
and we prove that the subspace
\[
\CK_T(F) := \bigcap_{n\in\N} F(T-\lambda_n\unit, r_n)
\]
has the properties mentioned in Theorem~\ref{embedding}.

  Finally, for arbitrary $B\in \B(\C)$ we show that
\[
\CK_T(B):= \overline{\bigcup_{K\subseteq B, \; K \;{\rm
      compact}}\CK_T(K)}
\]
does the job for us. 

\vspace{1cm}

We now describe the contents of the rest of the paper. The main result of
section~2 is that for $T$ in the II$_1$-factor $\CM$, the push-forward
measure of $\mu_{(T\cc)^nT^n}$ under the map $t\mapsto t^{\frac 1n}$
converges weakly (as $n\rightarrow \infty$) to the measure $\nu\in
\Prob([0,\infty[)$ which is uniquely determined by
\[
\nu([0,t^2])=\mu_T(\overline{B(0,t)}), \qquad (t>0).
\]

In section~3 we define for $T\in \CM$ and $r>0$ the $T$-hyperinvariant
subspaces $E(T,r)$ and $F(T,r)$ mentioned above, and we explore some of
their properties. We shall see that these
subspaces seem
to be good candidates for the desired  $T$-invariant subspaces $\CK_T(\overline{B(0,r)})$ and $\CK_T(\C\setminus
{B(0,r)})$, respectively. However, in order to prove that $E(T,r)$ and
$F(T,r)$ fulfill the requirements listed in Theorem~\ref{embedding}, some
more work has to be done.

We begin by considering the case $r=1$ and assume that $\mu_T(\partial
B(0,1))=0$. The idea of proof is the following: If $\sigma(T)\cap \partial
B(0,1) = \emptyset$, then one can always define an
idempotent $e\in\CM$ by
\begin{equation}\label{integral}
e= \frac{1}{2\pi \i} \int_{\partial B(0,1)}
(\lambda\unit-T)^{-1}\,\d\lambda
\end{equation}
- as a Banach space valued integral in $\CM$. Then it is a fact that the range
  projection $p$ of $e$ is $T$-invariant with
\[
\sigma(T|_{p(\CH)})\subseteq {B(0,1)} \quad {\rm and}\quad \sigma((\unit-p)T|_{p(\CH)^\bot})\subseteq \C\setminus\overline{B(0,1)}.
\]
Hence, the Brown measure of $T|_{p(\CH)}$ and $(\unit-p)T|_{p(\CH)^\bot}$ are concentrated on $B(0,1)$ and $\C\setminus \overline{B(0,1)}$, respectively.

However, in general one can not make sense of the integral in \eqref{integral}, but this can be remedied by adding a small perturbation
to $T$. We consider $\CM$ as a subfactor of the II$_1$-factor $\CN =
\CM\ast L(\F_4)$ (with a tracial state which we also denote by $\tau$), and we note that $\CN$ contains a circular system
$\{x,y\}$ which is $\ast$-free from $\CM$. Moreover, we can define the unbounded operator $z=xy^{-1}$ which by \cite[Theorem~5.2]{HS} belongs to $L^p(\CN,\tau)$ for $0<p<1$. We will consider the perturbations of $T$ given by
\[
T_n=T+\frac1n z, \qquad n\in\N.
\]
Since $T_n\in L^p(\CN,\tau)$, $0<p<1$, it has a well--defined Fuglede--Kadison determinant $\Delta(T_n)$ and a well--defined Brown measure $\mu_{T_n}$ (cf. \cite[Appendix]{Br} or \cite[Section~2]{HS}). In section~4 we will prove that
\[
\Delta(T_n)=\Delta(T\cc T+\textstyle{\frac{1}{n^2}}\unit)^\frac12
\]
and that
\[
\mu_{T_n}\rightarrow \mu_T \quad {\rm as} \quad n\rightarrow\infty
\]
in the weak topology on ${\rm Prob}(\C)$. 

In 1968, Turpin and Waelbroek introduced an approach to vector valued integration in quasinormed spaces such as $L^p(\CN,\tau)$, $0<p<1$ (cf. \cite{TuWa}, \cite{Wa}, \cite{Ka}). In particular, one can define the integral
\[
\int_a^b f(x)\,\d x
\] 
for every function $f: [a,b]\rightarrow L^p(\CN,\tau)$ which satisfies the H\"older condition
\[
\|f(x)-f(y)\|_p \leq C|x-y|^\alpha
\]
with exponent $\alpha>\frac1p -1$ (cf. section~10 for a selfcontained proof). Based on results from sections~5 and 6, we prove in section~7 that $(\lambda\unit-T_n)^{-1}\in L^p(\CN,\tau)$ for $\lambda\in\C$ and $0<p<1$. Moreover, if $0<p<\frac23$, then
\[
\|(\lambda\unit-T_n)^{-1}-(\mu\unit-T_n)^{-1}\|_p\leq C_{p,n}|\lambda-\mu|
\]
for some constant $C_{p,n}>0$. It follows that he integral
\[
e_n=\frac{1}{2\pi\i}\int_{\partial B(0,1)}(\lambda\unit-T_n)^{-1}\,\d\lambda
\]
makes sense as a Turpin--Waelbroek integral in $L^p(\CN,\tau)$ for $\frac12 < p <\frac23$. Let $P_n$ denote the range projection of $e_n$. Then
\begin{itemize}
  \item[(i)]$P_nT_nP_n = T_nP_n$,
  \item[(ii)]  $\supp(\mu_{T_n|_{P_n(\CH)}})\subseteq \overline{B(0,1)}$
  \item[(iii)]  $\supp(\mu_{P_n^\bot T_n|_{P_n(\CH)^\bot}})\subseteq \C\setminus B(0,1)$
\end{itemize}
(cf. Theorem~\ref{P_n}). Then we take a free ultrafilter $\omega$ on $\N$ and define
$P\in\CN^\omega$ to be the image of $(P_n)_{n=1}^\infty$ under the quotient
mapping $\rho : \ell^\infty(\CN)\rightarrow \CN^\omega$. Using (i), (ii)
and (iii), one can prove that $P$ is $T$-invariant,
\[
\tau_\omega(P)=\mu_T(\overline{B(0,1)}),
\]
\[
\supp(\mu_{T|_{P(\CH)}})\subseteq \overline{B(0,1)},
\]
and
\[
\supp(\mu_{(\unit-P)T|_{P(\CH)^\bot}})\subseteq \C\setminus B(0,1).
\]
We then prove that $P=P_{E(T,1)}$ (cf. Lemma~\ref{Pr(H)=E(T,r)}), and thus the
$T$-hyperinvariant subspace $E(T,1)$ has the desired properties. The last
part of section~7 takes care of a general Borel set $B\subseteq \C$, as was
outlined above.

In section~8 we realize the $E(T,r)$ and $F(T,r)$ as spectral
projections of the positive operators $SO-\lim_{n\rightarrow
  \infty}((T\cc)^n T^n)^{\frac{1}{2n}}$ and $SO-\lim_{n\rightarrow
  \infty}(T^n(T\cc)^n)^{\frac{1}{2n}}$, respectively. In particular, we
prove that these two limits actually exist for every $T\in\CM$ when $\CM$ is a II$_1$--factor. There are examples of bounded operators on a Hilbert space which do {\it
  not} have this property (cf. Example~8.4). 
  
Finally, in section~9, we show that for operators $T$ which are decomposable in the sense of \cite[Section~1.2]{LN}, our spectral subspaces $\CK(T,F)$ for closed subsets $F$ of $\C$ coincide with the subspaces in the spectral capacity of $T$. Moreover, we show that every II$_1$--factor contains a non--decomposable operator. 

Throughout this paper we assume that $\CM$ is a II$_1$--factor. However, our main results can easily be generalized to the case where $\CM$ is a finite von Neumann algebra with a specified normal faithful trace $\tau$. Indeed, such a von Neumann algebra $\CM$ can be embedded into a II$_1$--factor $\CN$ in such a way that the restriction of $\tau_\CN$ to $\CM$ agrees with $\tau$ (cf. \cite[Proof of Theorem~2.6]{HW}).

%% file: differentialeq.tex
\section{Some results on the Brown measure of a bounded operator}

Consider a II$_1$-factor $\CM$ with faithful tracial state $\tau$. As was
mentioned in the introduction, one can associate to every $T\in\CM$ a
probability measure $\mu_T$ on $\C$, the Brown measure of $T$. 

\begin{remark}\label{remarks}
\begin{itemize}
\item[(i)] Note that if $T\in\CM$, and if $A\in \CM$ is invertible, then $\mu_T=\mu_{ATA^{-1}}$. 
\item[(ii)] According to \cite{Br}, we have the following generalization of Weil's Theorem: For $T\in\CM$ and $0<p<\infty$,
\begin{equation}\label{Weil}
\int_{\sigma(T)}|\lambda|^p\,\d\mu_T(\lambda)\leq \|T\|_p^p := \tau(|T|^p)
\end{equation}
\end{itemize}
\end{remark}

\vspace{.2cm}

\noindent Fack and Kosaki (cf. \cite{FK}) proved the generalized H\" older inequality 
\begin{equation}\label{Hoelder}
\|ST\|_r \leq \|S\|_p\|T\|_q, \qquad (S,T\in\CM),
\end{equation}
which holds for all $0<p,q,r\leq \infty$ with $\frac 1r = \frac 1p + \frac 1q$. As a consequence of \eqref{Weil} and \eqref{Hoelder} we have (cf. \cite{Br}):
\begin{equation}\label{Hoelder+Weil}
\|T^n\|_{\frac pn}\leq \|T\|_p^n, \qquad (T\in\CM, \; p>0).
\end{equation}

\vspace{.2cm}

The main result of this section is

\begin{thm}\label{MSRI2.5} Let $T\in\CM$, and for $n\in\N$, let
  $\mu_n\in{\rm Prob}([0,\infty[)$ denote the distribution of $(T^n)\cc
  T^n$ w.r.t. $\tau$, and let $\nu_n$ denote the push-forward measure of
  $\mu_n$ under the map $t\mapsto t^{\frac 1n}$. Moreover, let $\nu$ denote
  the push-forward measure of $\mu_T$ under the map $z\mapsto |z|^2$,
  i.e. $\nu$ is determined by
  \[
  \nu([0,t^2])= \mu_T(\overline{B(0,t)}), \qquad (t>0).
  \]
  Then $\nu_n\rightarrow \nu$ weakly in ${\rm Prob}([0,\infty[)$.
\end{thm}

\vspace{.2cm}

We will obtain Theorem~\ref{MSRI2.5} as a consequence of

\begin{thm}\label{existence of A_k} Let $T\in\CM$. Then there is a sequence $(A_k)_{k=1}^\infty$ in $\CM_{\rm inv}$, such that
\begin{itemize}
\item[(i)] $\|A_kTA_k^{-1}\|\leq \|T\|$ for all $k\in\N$,
\item[(ii)] $A_kTA_k^{-1}$ converges in $\ast$-distribution to a normal
  operator $N$ in an ultrapower of $\CM$, $\CM^\omega$, with $\mu_N=\mu_T$.
\end{itemize}
\end{thm}

\vspace{.2cm}

Before proving Theorem~\ref{existence of A_k} we
state and prove some of its additional consequences:

\begin{cor} For every $T\in\CM$ and every $p>0$,
\begin{equation}\label{infimum-inv}
\int_\C |\lambda|^p\,\d\mu_T(\lambda)=\inf_{A\in\CM_{\rm inv}}\|ATA^{-1}\|_p^p.
\end{equation}
\end{cor}

\noindent \proof According to Remark~\ref{remarks}~(i) and \eqref{Weil}, 
\[
\int_\C |\lambda|^p\,\d\mu_T(\lambda) \leq \inf_{A\in\CM_{\rm inv}}\|ATA^{-1}\|_p^p.
\]
To see that $\geq$ holds in \eqref{infimum-inv}, take $(A_k)_{k=1}^\infty$ as in Theorem~\ref{existence of A_k}. Then $\mu_{(A_kTA_k^{-1})\cc A_kTA_k^{-1}}\rightarrow \mu_{N\cc N}$ in moments, and since all of the measures $(\mu_{(A_kTA_k^{-1})\cc A_kTA_k^{-1}})_{k=1}^\infty$ are supported on $[0,\|T\|^2]$, this implies weak convergence. In particular,
\begin{eqnarray*}
\|A_kTA_k^{-1}\|_p^p &=& \int_0^\infty t^{\frac p2}\,\d\mu_{(A_kTA_k^{-1})\cc A_kTA_k^{-1}}(t)\\
&\rightarrow & \int_0^\infty t^{\frac p2}\,\d\mu_{N\cc N}(t)\\
&=& \tau(|N|^p)\\
&=& \int_\C |\lambda|^p \,\d\mu_N(\lambda)\\
& = & \int_\C |\lambda|^p\,\d\mu_T(\lambda),
\end{eqnarray*}
proving that $\geq$ holds in \eqref{infimum-inv}. $\endproof$

\vspace{.2cm}

\begin{thm}\label{MSRI2.3}
For every $T\in\CM$ and every $p>0$,
\begin{equation}\label{2.6}
\int_\C |\lambda|^p\,\d\mu_T(\lambda)=\lim_{n\rightarrow \infty}\|T^n\|_{\frac pn}^{\frac pn}.
\end{equation}
\end{thm}

\noindent \proof According to \cite{Br}, $\mu_{T^n}$ is the push-forward measure of $\mu_T$ under the map $z\mapsto z^n$. Hence by \eqref{Weil}, for all $n\in\N$,
\begin{equation}\label{jetlag1}
\int_\C |\lambda|^p \,\d\mu_T(\lambda) = \int_\C |\lambda|^{\frac pn}\,\d\mu_{T^n}(\lambda)\leq \|T^n\|_{\frac pn}^{\frac pn}.
\end{equation}
 
Let $\eps>0$. It follows from \eqref{infimum-inv} that we may choose $A\in\CM_{\rm inv}$, such that
\[
\|ATA^{-1}\|_p^p \leq \int_\C|\lambda|^p\,\d\mu_T(\lambda) + \eps.
\]
With $S=ATA^{-1}$ we then have (cf. \eqref{Hoelder+Weil}) that
\begin{eqnarray*}
\|S^n\|_{\frac pn}^{\frac pn} &\leq & \|S\|_p^p\\
&\leq & \int_\C|\lambda|^p\,\d\mu_T(\lambda) + \eps.
\end{eqnarray*}
Since $T^n=A^{-1}S^n A$, the generalized H\"older inequality \eqref{Hoelder} implies that
\[
\|T^n\|_{\frac pn}\leq \|A^{-1}\|\|A\|\|S^n\|_{\frac pn}.
\]
Thus,
\[
\|T^n\|_{\frac pn}^{\frac pn}\leq \|A^{-1}\|^{\frac pn}\|A\|^{\frac pn}\Big( \int_\C|\lambda|^p\,\d\mu_T(\lambda) + \eps\Big),
\]
and it follows that
\begin{equation}\label{jetlag2}
\limsup_{n\rightarrow\infty}\|T^n\|_{\frac pn}^{\frac pn}\leq \int_\C|\lambda|^p\,\d\mu_T(\lambda) + \eps.
\end{equation}
Combining \eqref{jetlag1} and \eqref{jetlag2} we find that
\[
\int_\C|\lambda|^p\,\d\mu_T(\lambda)\leq \liminf_{n\rightarrow\infty}\|T^n\|_{\frac pn}^{\frac pn}\leq  \limsup_{n\rightarrow\infty}\|T^n\|_{\frac pn}^{\frac pn}\leq \int_\C|\lambda|^p\,\d\mu_T(\lambda),
\]
and the theorem follows. $\endproof$

\vspace{.2cm}

\begin{cor} For $T\in\CM$ define $r'(T)$, {\emph the modified spectral radius
    of $T$}, by
\[
r'(T):=\max\{|\lambda|\,|\,\lambda\in\supp(\mu_T)\}.
\]
Then
\begin{equation}
r'(T)=\lim_{p\rightarrow\infty}\Big(\lim_{n\rightarrow\infty}\|T^n\|_{\frac pn}^{\frac 1n}\Big).
\end{equation}
\end{cor}

\noindent \proof $r'(T)$ is the essential supremum (w.r.t. $\mu_T$) of the map $\lambda\mapsto |\lambda|$. Hence,
\[
r'(T) = \lim_{p\rightarrow \infty}\Big(\int_\C |\lambda|^p\,\d\mu_T(\lambda)\Big)^{\frac 1p}.
\]
Now apply Theorem~\ref{MSRI2.3}. $\endproof$

\vspace{.2cm}

\noindent {\it Proof of Theorem~\ref{existence of A_k}.} If
$A:[0,\infty[\rightarrow \CM_{\rm inv}$ is a differentiable map, define $X:[0,\infty[\rightarrow \CM$ by
\[
X(t)=A(t)TA(t)^{-1}, \qquad (t>0).
\]
Note that 
\[
\frac{\d}{\d t} A(t)^{-1}= -A(t)^{-1}A'(t)A(t)^{-1}.
\]
Hence, $X$ is differentiable with
\begin{eqnarray*}
\frac{\d}{\d t}X(t)&=& A'(t)TA(t)^{-1} + A(t)T \frac{\d}{\d t}A(t)^{-1}\\
&=& [A'(t)A(t)^{-1}, X(t)].
\end{eqnarray*}
We will choose $A(t)$ to be the solution of the differential equation 
\begin{equation}
\left\{
\begin{array}{lll}
A'(t)= B(t)A(t)&,& t\geq 0\\
A(0)=\unit
\end{array}
\right .
\end{equation}
for a suitable function $B:[0,\infty[\rightarrow \CM_{{\rm sa}}$, chosen in such a way that the identity
\[
\frac{\d}{\d t}X(t)=[B(t),X(t)]
\]
implies that $t\mapsto \|X(t )\|_p^p$ is decreasing for all $p\in\N$. At first we consider the case $p=2$:
\begin{eqnarray*}
\frac{\d}{\d t}\|X(t)\|_2^2 &=& \frac{\d}{\d t} \tau(X(t)\cc X(t))\\
&=& \tau([B(t),X(t)]\cc X(t) + X(t)\cc [B(t),X(t)])\\
&=&2\tau\big(B(t)(X(t)X(t)\cc- X(t)\cc X(t))\big).
\end{eqnarray*}
Hence, if $B(t)=[X(t)\cc, X(t)]$, then 
\[
\frac{\d}{\d t} \|X(t)\|_2^2 = -2\|[X(t), X(t)\cc]\|_2^2\leq 0.
\]
Therefore, in the following we will try to solve the differential equation
\begin{equation}\label{DE1}
\left\{
\begin{array}{lll}
X'(t)= [[X(t)\cc, X(t)], X(t)]&,& t\geq 0\\
X(0)=T
\end{array}
\right .
\end{equation}
and then solve
\begin{equation}\label{DE2}
\left\{
\begin{array}{lll}
A'(t)= [X(t)\cc, X(t)]A(t)&,& t\geq 0\\
A(0)=\unit.
\end{array}
\right .
\end{equation}
We are going to apply \cite[Chapter~14, Theorem~3.1]{La} to \eqref{DE1}. Assuming that $T\neq 0$, set $\CU=B(0,2\|T\|)\subseteq \CM$, and let
\[
f(X)= [[X\cc, X],X], \qquad X\in\CM.
\]
For fixed $X\in\CM$, $f'(X)$ is a bounded operator on $\CM$ given by
\[
f'(X): H\mapsto [[H\cc, X],X] + [[X\cc, H],X] + [[X\cc,X],H].
\]
For $X\in \CU$,
\[
\|f'(X)\|\leq 12 (2\|T\|)^2=48\|T\|^2.
\]
Hence, $K=48\|T\|^2$ is a Lipschitz constant for $f|_\CU$. $\|f(X)\|$ is bounded on $\CU$ by
\[
L=4(2\|T\|)^3=32\|T\|^3.
\]
With $a=\frac13\|T\|$, $\overline{B_{2a}(T)}\subseteq\CU$. Set
\[
b_0= \min\Big\{\frac1K, \frac aL\Big\} = \min\Big\{\frac{1}{48\|T\|^2}, \frac{1}{96\|T\|^2}\Big\}.
\]
Then by \cite[Chapter~14, Theorem~3.1]{La}, for all $b\in (0,b_0)$ there is a unique solution to 
\begin{equation*}
\left\{
\begin{array}{lll}
X'(t)= f(X(t))\\
X(0)=T
\end{array}
\right .
\end{equation*}
defined on the interval $(-b,b)$. In particular, such a solution exists and is unique on $\Big[0,\frac{c}{\|T\|^2}\Big]$, where $c=\frac{1}{100}$. Next, with $B=[X\cc, X]$ we get that for every $p\in \N$,
\begin{eqnarray*}
\frac{\d}{\d t}\|X(t)\|_{2p}^{2p}&=& \frac{\d}{\d t}\tau((X(t)\cc X(t))^p)\\
&=& p\, \tau\Big(\frac{\d}{\d t}\big(X(t)\cc X(t)\big)(X(t)\cc X(t))^{p-1}\Big)\\
&=& p\,\tau([2X(t)\cc B(t)X(t)-B(t)X(t)\cc X(t)-X(t)\cc X(t)B(t)](X(t)\cc X(t))^{p-1})\\
&=& 2p\,\tau(B(t)[(X(t)X(t)\cc)^p - (X(t)\cc X(t))^p])\\
&=& 2p\, \tau([X(t)\cc X(t)-X(t)X(t)\cc][(X(t)X(t)\cc)^p - (X(t)\cc X(t))^p]).
\end{eqnarray*}
For fixed $t\in \Big[0, \frac{c}{\|T\|^2}\Big]$ there is a unique compactly supported probability measure $\mu_t$ on $[0,\infty[\times [0,\infty[$, such that for all $f, g\in C([0,\infty[)$,
\[
\tau(f(X(t)\cc X(t))g(X(t)X(t)\cc))=\int_{[0,\infty[\times [0,\infty[} f(u)g(v)\,\d\mu_t(u,v)
\]
(cf. \cite[Proposition~1.1]{Co}). In particular, 
\begin{equation}\label{circstar1}
\frac{\d}{\d t}\|X(t)\|_{2p}^{2p} = 2p \int_{[0,\infty[\times [0,\infty[}(u-v)(v^p-u^p)\,\d\mu_t(u,v) \leq 0,
\end{equation}
and it follows that $t\mapsto \|X(t)\|_p$ is decreasing for all $p\in\N$. Thus, also $t\mapsto \|X(t)\|_\infty= \lim_{p\rightarrow\infty}\|X(t)\|_p$ is decreasing. We can therefore extend the solution to \eqref{DE1} to the interval $\Big[\frac{c}{\|T\|^2},\frac{2c}{\|T\|^2}\Big]$, and repeating the argument, we find that \eqref{DE1} has a solution defined on all of $[0,\infty[$ and satisfying $\|X(t)\|\leq \|T\|$, $t\geq 0$. Then with $B(t)=[X(t)\cc, X(t)]$ we have
\begin{equation}\label{25/11a}
\|B(t)\|\leq 2\|T\|^2, \qquad t\geq 0.
\end{equation}
In order to solve \eqref{DE2}, we apply the method of proof of \cite[Chap.~14, Theorem~3.1]{La} and define $(A_n(t))_{n=0}^\infty$ recursively by
\begin{eqnarray*}
A_0(t) &=& \unit, \qquad t\geq 0,\\
A_n(t) & = & \unit + \int_0^t B(t')A_{n-1}(t')\,\d t', \qquad t\geq 0.
\end{eqnarray*}
By induction on $n$, we find that
\[
\|A_n(t)-A_{n-1}(t)\|\leq \frac{(2\|T\|^2 t)^n}{n!}, \qquad n\in\N.
\]
Hence,
\[
1+ \sum_{n=1}^\infty \|A_n(t)-A_{n-1}(t)\|\leq e^{2\|T\|^2t} <\infty,
\]
and it follows that the limit
\[
A(t)=\lim_{n\rightarrow\infty}A_n(t)
\]
exists for all $t\geq 0$. Moreover, the convergence is uniform on compact subsets of $[0,\infty[$. Therefore,
\[
A(t)= \unit + \int_0^t B(t')A(t')\,\d t', \qquad t\geq 0,
\]
showing that $A(t)$ is a $C^1$ solution to \eqref{DE2}. By a similar argument, the problem
\begin{equation}
\left\{
\begin{array}{lll}
C'(t)= -C(t)B(t)&,& t\geq 0,\\
C(0)=\unit
\end{array}
\right .
\end{equation}
has a $C^1$ solution defined on $[0,\infty)$. Moreover,
\[
\frac{\d}{\d t}\big(C(t)A(t)\big) = C'(t)A(t) + C(t)A'(t) = 0,
\]
which implies that $C(t)A(t)=\unit$, $t\geq 0$. Hence, $A(t)$ is invertible with inverse $C(t)$ for all $t\geq 0$. Now, observe that both $X(t)$ and $A(t)TA(t)^{-1}$ are $C^1$ solutions to the problem
\begin{equation}
\left\{
\begin{array}{lll}
Y'(t)= [B(t),Y(t)]&,& t\geq 0,\\
Y(0)=T.
\end{array}
\right .
\end{equation}
Then by \cite[Chap.~14, Theorem~3.3]{La}, $X(t)=A(t)TA(t)^{-1}$ for all $t\geq 0$. That is, $A(t)$, $B(t)$ and $X(t)$ have all the properties we were asking for in the beginning of this proof. 

Note that the map $t^p\mapsto p\cdot t^p$, $t\in [0,\infty[$, extends to a
map from $C^1([0,\infty[)$ into $C([0,\infty[)$ given by
\[
\phi(t)\mapsto t\cdot \phi'(t).
\]
Let $\phi\in C^1([0,\infty[)$. Choosing a sequence of polynomials
$(p_n)_{n=1}^\infty$, such that $p_n\rightarrow \phi$ and $p_n'\rightarrow
\phi'$ uniformly on $[0,\|T\|^2]$, we get from \eqref{circstar1} that
\begin{equation}\label{circstar2}
  \difft \tau(\phi(X(t)\cc X(t)))= -2\, \int_{[0,\infty[\times [0,\infty[}
  (u-v)(u\phi'(u) -v\phi'(v))\,\d\mu_t(u,v).
\end{equation}
We shall need this identity later on.

Next we prove that as $t\rightarrow \infty$, $X(t)$ converges in
$\ast$-moments to a normal operator $N$. Recall that
\begin{equation}\label{(a)}
  \difft \|X(t)\|_2^2=
  -2\,\|[X(t)\cc,X(t)]\|_2^2.
\end{equation}
Moreover,
\begin{equation}\label{(b)}
  X'(t)= [[X(t)\cc, X(t)],X(t)]= X(t)\cc X(t)^2 - 2X(t)X(t)\cc X(t)+
  X(t)^2X(t)\cc.
\end{equation}
Since $\tau(Y[Z,W])=\tau(Z[W,Y])=\tau(W[Y,Z])$, we have that
\begin{eqnarray}
  \difft \|[X(t)\cc, X(t)]\|_2^2 &=& 2\,\tau\Big([X(t)\cc, X(t)]\cdot
  \difft [X(t)\cc, X(t)]\Big)\nonumber\\
  &=& 2\,\tau\Big([X(t)\cc, X(t)]\Big([X'(t)\cc, X(t)]+ [X(t)\cc,
  X'(t)]\Big)\Big)\nonumber\\
  &=& 2\,\tau(X'(t)\cc[X(t),[X(t)\cc, X(t)]])+2\,\tau(X'(t)[[X(t)\cc,
  X(t)], X(t)\cc])\nonumber\\
  &=& -2\,\tau(X'(t)\cc[[X(t)\cc, X(t)], X(t)])-2\,\tau(X'(t)[X(t)\cc,[X(t)\cc,
  X(t)]])\nonumber\\
  &=& -4 \|X'(t)\|_2^2 \label{(c)}.
\end{eqnarray}
Hence, $t\mapsto \|[X(t)\cc, X(t)]\|_2^2$ is decreasing, and it follows
that for all $t>0$,
\[
\|[X(t)\cc, X(t)]\|_2^2\leq \frac 1t \int_0^t \|[X(u)\cc, X(u)]\|_2^2\,\d
u,
\]
so by \eqref{(a)},
\[
\|[X(t)\cc, X(t)]\|_2^2 \leq \frac{1}{2t} \Big(\|X(0)\|_2^2-
\|X(t)\|_2^2\Big)\leq \frac{1}{2t} \|T\|_2^2.
\]
This shows that
\begin{equation}\label{(d)}
  \lim_{t\rightarrow\infty}\|[X(t)\cc, X(t)]\|_2^2=0.
\end{equation}
According to \eqref{circstar1}, $t\mapsto \|X(t)\|_{2p}^{2p}$ is decreasing
for every $p\in\N$. Therefore,
\[
\lim_{t\rightarrow\infty}\tau((X(t)\cc X(t))^p)
\]
exists for every $p\in\N$. Combining this with \eqref{(d)} and
the fact that $\|X(t)\|_\infty \leq \|T\|$ for every $t\geq 0$, we get that the
trace of any monomial in $X(t)$ and $X(t)\cc$ with the same number of
$X(t)$'s and $X(t)\cc$'s also converges as $t\rightarrow\infty$. To prove
that this holds for any monomial, it is, because of \eqref{(d)}, sufficient
to show that for all $m, n\in \N_0$, the limit
\begin{equation}\label{(e)}
\lim_{t\rightarrow\infty}\tau((X(t)\cc)^m X(t)^n)
\end{equation}
exists. To see this, note that for every $\lambda\in\C$,
$X(t)-\lambda\unit$ satisfies the same differential equation as $X(t)$,
because
\begin{eqnarray*}
  \difft (X(t)-\lambda\unit)&=& [[X(t)\cc, X(t)],X(t)]\\
  &=&  [[X(t)\cc-\overline\lambda \unit,
  X(t)-\lambda\unit],X(t)-\lambda\unit].
\end{eqnarray*}
Hence, according to the above, the limit
\begin{equation}\label{(f)}
  \lim_{t\rightarrow\infty}\tau((X(t)\cc-\overline\lambda
  \unit)^n(X(t)-\lambda\unit)^n)
\end{equation}
exists for all $n\in\N$. Expanding $\tau((X(t)\cc-\overline\lambda
  \unit)^n(X(t)-\lambda\unit)^n)$ in powers of $\lambda$ and
  $\overline\lambda$, it is readily seen that the limit \eqref{(e)}
  exists. Hence, $(X(t))_{t\geq 0}$ converges in $\ast$-moments to an
  operator $N$ ($N$ may be realized in the ultrapower $\CM^\omega$ of
  $\CM$). According to \eqref{(d)}, $N$ is normal.

In order to show that $\mu_N=\mu_T$, it suffices to show that
\begin{equation}
  \Delta(T-\lambda\unit)=\Delta(N-\lambda\unit), \qquad (\lambda\in\C).
\end{equation}
At first we consider the case $\lambda=0$. Since $X(t)=A(t)TA(t)^{-1}$,
$\Delta(T)=\Delta(X(t))$ for all $t\geq 0$. Hence,
\begin{eqnarray}
  \log\Delta(T)&=& \log\Delta(X(t))\nonumber\\
  &=& \inf_{\eps>0}\Big\{\frac 12 \tau(\log(X(t)\cc X(t)+\eps \unit))\Big\}.
  \label{(g)}
\end{eqnarray}
Applying now \eqref{circstar2} to $\phi(u)=\log(u+\eps)$, $u\geq 0$, we get
that
\begin{equation}\label{(h)}
  \difft \tau(\log(X(t)\cc X(t)+\eps\unit))=-2\, \int_0^\infty \int_0^\infty
  (u-v)\Big(\frac{u}{u+\eps}-\frac{v}{v+\eps} \Big)\,\d\mu_t(u)\d\mu_t(v)
  \leq 0.
\end{equation}
Now, according to \eqref{(g)},
\begin{eqnarray*}
  \log\Delta(T)&=& \inf_{t>0}\inf_{\eps>0}\Big\{\frac 12 \tau(\log(X(t)\cc
  X(t)+\eps \unit))\Big\}\\
  &=&  \inf_{\eps>0}\inf_{t>0}\Big\{\frac 12 \tau(\log(X(t)\cc
  X(t)+\eps \unit))\Big\}\\
  &=& \inf_{\eps>0}\lim_{t\rightarrow\infty}\Big\{\frac 12 \tau(\log(X(t)\cc
  X(t)+\eps \unit))\Big\},
\end{eqnarray*}
where the last identity follows from \eqref{(h)}. But $X(t)\cc X(t)$
converges in moments to $N\cc N$, and since $\sup_{t>0}\|X(t)\cc
X(t)\|<\infty$, it follows from the Weierstrass approximation theorem that
\begin{equation}\label{(j)}
  \forall\; \phi\in C([0,\infty[): \quad
  \lim_{t\rightarrow\infty}\tau(\phi(X(t)\cc X(t)))= \tau_\omega(\phi(N\cc
  N)),
\end{equation}
where $\tau_\omega$ denotes the tracial state on $\CM^\omega$. Hence,
\[
\log \Delta(T)= \inf_{\eps>0}\Big\{\frac 12 \tau(\log(N\cc
  N+\eps \unit))\Big\}= \log\Delta(N).
\]
The same arguments apply to $X(t)-\lambda\unit$, and we obtain that
\[
\log\Delta(T-\lambda\unit)=\log\Delta(N-\lambda\unit), \qquad
(\lambda\in\C).
\]
Hence, $\mu_T=\mu_N$. Theorem~\ref{existence of A_k} now follows by taking
$A_k= A(k)$, $k\in\N$. $\endproof$
%Applying \eqref{(j)} to $\phi(u)=u^{\frac p2}$, it
%also follows that $\|X(t)\|_p\rightarrow \|N\|_p$ as
%$t\rightarrow\infty$. Hence,
%\[
%\lim_{t\rightarrow\infty}\|X(t)\|_p^p = \int_\C
%|\lambda|^p\,\d\mu_N(\lambda)= \int_\C
%|\lambda|^p\,\d\mu_T(\lambda).
%\]
\vspace{.2cm}

{\it Proof of Theorem~\ref{MSRI2.5}.} Let $p\in [0,\infty[$. Then, according to Theorem~\ref{MSRI2.3},
\begin{eqnarray*}
  \int_0^\infty t^{\frac p2}\,\d\nu(t)&=& \int_\C
  |\lambda|^p\,\d\mu_T(\lambda)\\
   &=& \lim_{n\rightarrow\infty} \|T^n\|_{\frac pn}^{\frac pn}\\
   &=& \lim_{n\rightarrow\infty} \tau(((T^n)\cc T^n)^{\frac{p}{2n}})\\
   &=& \lim_{n\rightarrow\infty} \int_0^\infty
  t^{\frac{p}{2n}}\,\d\mu_n(t)\\
  &=& \lim_{n\rightarrow\infty} \int_0^\infty
  t^{\frac{p}{2}}\,\d\nu_n(t).
\end{eqnarray*}
In particular, $\nu_n\rightarrow \nu$ in moments. Since $\|(T^n)\cc
T^n\|\leq \|T\|^{2n}$ and $\supp(\mu_T)\subseteq \overline{B(0, \|T\|)}$,
 $\nu$ and the $\nu_n$ are all supported on $[0,\|T\|^2]$. Hence,
 $\nu_n\rightarrow \nu$ weakly. $\endproof$

\vspace{.2cm}

\begin{cor} Let $T\in\CM$. Then $\mu_T=\delta_0$ if and only if $((T^n)\cc
  T^n)^{\frac 1n}$ tends to $0$ in the strong operator topology as
  $n\rightarrow \infty$.
\end{cor}

\proof If  $\mu_T=\delta_0$, then with the same notation as in
Theorem~\ref{MSRI2.5}, $\nu_n\rightarrow \delta_0$ weakly as $n\rightarrow
\infty$. In particular,
\[
\lim_{n\rightarrow \infty}\int_0^{\|T\|^2}t^4\,\d\nu_n(t)= 0,
\]
i.e.
\[
\lim_{n\rightarrow\infty} \tau([(T^n)\cc T^n]^{\frac 4n})=0.
\]
It follows that for every vector $\xi$ in the dense subspace $\CM$ of
$L^2(\CM,\tau)$,
\begin{eqnarray*}
  \|[(T^n)\cc T^n]^{\frac 1n}\xi\|_2^2 &=& \tau([(T^n)\cc T^n]^{\frac
  2n}\xi\xi\cc)\\
&\leq & \|[(T^n)\cc T^n]^{\frac 2n}\|_2\|\xi\xi\cc\|_2\\
&\rightarrow & 0,
\end{eqnarray*}
and hence, $((T^n)\cc T^n)^{\frac 1n}$ tends to $0$ in the strong operator
topology as $n\rightarrow \infty$. On the other hand, if  $((T^n)\cc
T^n)^{\frac 1n}\rightarrow 0$ strongly, then $((T^n)\cc
T^n)^{\frac pn}\rightarrow 0$ strongly for every $p\in\N$. Hence, with $\xi = \unit\in\CM\subseteq L^2(\CM,\tau)$,
\[
  \tau([(T^n)\cc T^n]^{\frac{2p}{n}}) =  \|[(T^n)\cc
  T^n]^{\frac{p}{n}}\xi\|_2^2 \rightarrow 0
\]
as $n\rightarrow \infty$ for every $p\in\N$, i.e.
\[
\lim_{n\rightarrow\infty}\int_0^{\|T\|^2}t^p\,\d\nu_n(t)= 0, \qquad (p\in\N).
\]
It follows then from the Weierstrass approximation theorem that
$\nu_n\rightarrow \delta_0$ weakly as $n\rightarrow \infty$. That is, $\nu= \delta_0$, and therefore
$\mu_T=\delta_0$ as well. $\endproof$

%% file: spectralsubspaces.tex
\section{Constructing certain hyperinvariant subspaces}\label{sec6}

Consider a II$_1$-factor $\CM$ acting on the Hilbert space $\CH$ and with
faithful, tracial state $\tau$.

\begin{definition}\label{E(T,r)}
  \begin{itemize}
    \item[(i)] For $T\in\CM$ and $r>0$, let $E(T,r)$ denote the set of
    $\xi\in\CH$, for which there exists a sequence $(\xi_n)_{n=1}^\infty$
    in $\CH$, such that
    \begin{equation}\label{UH4}
    \lim_{n\rightarrow\infty}\|\xi_n-\xi\|=0\quad {\rm and}\quad
    \limsup_{n\rightarrow\infty} \|T^n\xi_n\|^{\frac 1n}\leq r.
    \end{equation}
    \item[(ii)] For $T\in\CM$ and $r>0$, let $F(T,r)$ denote the set of
    $\eta\in\CH$, for which there exists a sequence $(\eta_n)_{n=1}^\infty$
    in $\CH$, such that
    \begin{equation}\label{UH5}
    \lim_{n\rightarrow\infty}\|T^n\eta_n-\eta\|=0\quad {\rm and}\quad
    \limsup_{n\rightarrow\infty} \|\eta_n\|^{\frac 1n}\leq \frac1r.
    \end{equation}
  \end{itemize}
\end{definition}

\vspace{.2cm}

\begin{lemma}\label{Properties} For $T\in\CM$ and $r>0$ one has:
  \begin{itemize}
    \item[(a)] $E(T,r)$ and $F(T,r)$ are closed subspaces of $\CH$.
    \item[(b)] $r\mapsto E(T,r)$ is increasing, and
    $E(T,r)=\bigcap_{s>r}E(T,s)$.
    
    $r\mapsto F(T,r)$ is decreasing, and
    $F(T,r)=\bigcap_{0<s<r}F(T,s)$.
    \item[(c)] $E(T,r)$ and $F(T,r)$ are hyperinvariant for $T$, i.e. for
    every $S\in \{T\}'$, $S(E(T,r))\subseteq E(T,r)$, and $S(F(T,r))\subseteq F(T,r)$.
    \item[(d)] The projections $P(T,r)$ and $Q(T,r)$ onto $E(T,r)$ and
    $F(T,r)$, respectively, belong to $W\cc(T)$, and they are independent
    of the particular representation of $\CM$ on a Hilbert space.
  \end{itemize}
\end{lemma}

\proof (a) It is easily seen that $E(T,r)$ and $F(T,r)$ are subspaces of
$\CH$. To see that $E(T,r)$ is closed, suppose $(\xi^{(k)})_{k=1}^\infty$ is a
sequence in $E(T,r)$ converging to some $\xi\in\CH$. We prove that $\xi$
belongs to $E(T,r)$ as well. For each $k\in\N$ there is a sequence
$(\xi_n^{(k)})_{n=1}^\infty$ in $\CH$, such that
 \[
    \lim_{n\rightarrow\infty}\|\xi_n^{(k)}-\xi^{(k)}\|=0\quad {\rm and}\quad
    \limsup_{n\rightarrow\infty} \|T^n\xi_n^{(k)}\|^{\frac 1n}\leq r.
\]

Choose a (strictly) increasing sequence of positive integers
$(n_k)_{k=1}^\infty$ such that for all $k\in\N$, 
\begin{equation}\label{UH1}
\|\xi_n^{(k)}-\xi^{(k)}\|<\frac1k, \qquad (n\geq n_k),
\end{equation}
and
\begin{equation}\label{UH2}
\|T^n\xi_n^{(k)}\|^{\frac 1n}<r+\frac1k, (n\geq n_k).
\end{equation}
Then define a sequence $(\xi_n)_{n=1}^\infty$ in $\CH$ by
\begin{equation}\label{UH3}
\xi_n=\left\{\begin{array}{lll} 0&, & n<n_k\;{\rm for}\;{\rm all} \; k\\
    \xi_n^{(k)}&, &n_k\leq n<n_{k+1}
    \end{array}\right.
\end{equation}

Then
\[
\|\xi_n-\xi^{k}\|\leq \frac1k, \quad (n_k\leq n<n_{k+1}),
\]
and
\[
\|T^n\xi_n\|^{\frac1n}\leq r+\frac1k, \quad (n_k\leq n<n_{k+1}).
\]
It follows that
\[
\xi = \lim_{k\rightarrow\infty}\xi^{(k)}=\lim_{n\rightarrow\infty}\xi_n,
\]
and
\[
\limsup_{n\rightarrow\infty}\|T^n\xi_n\|^{\frac1n}\leq r.
\]
Thus, $\xi\in E(T,r)$.  Modifying the above arguments a bit, one easily sees
that $F(T,r)$ is closed as well.

(b) Clearly, $r\mapsto E(T,r)$ is increasing, and therefore
\[
E(T,r)\subseteq \bigcap_{s>r}E(T,s).
\]
To see that the reverse inclusion holds, let $\xi\in
\bigcap_{s>r}E(T,s)$. For each $k\in\N$ take a sequence
$(\xi_n^{(k)})_{n=1}^\infty$ in $\CH$, such that
\[
    \lim_{n\rightarrow\infty}\|\xi_n^{(k)}-\xi\|=0\quad {\rm and}\quad
    \limsup_{n\rightarrow\infty} \|T^n\xi_n^{(k)}\|^{\frac 1n}<r+\frac1k.
\]

As in the above, choose a (strictly) increasing sequence of positive integers
$(n_k)_{k=1}^\infty$ such that for all $k\in\N$, \eqref{UH1} and \eqref{UH2} hold. Let $\xi_n$ be given by \eqref{UH3}, and
as in the proof of (a), note that $(\xi_n)_{n=1}^\infty$ satisfies
\eqref{UH4}, i.e. $\xi\in E(T,r)$. 

The second statement in (b) is proved in a similar way.

(c) This follows immediately from the definition of $E(T,r)$ and $F(T,r)$
and the fact that when $ST=TS$, then $ST^n=T^nS$.

(d) Since $E(T,r)$ and  $F(T,r)$ are hyperinvariant for $T$, $P(T,r)$ and
$Q(T,r)$ belong to $W\cc(T)=\{T,T\cc\}''$. The independence of the
representation of $\CM$ follows exactly as in \cite[Lemma~2.3 and
2.4]{DH}. $\endproof$

\vspace{.2cm}

\begin{remark} If $T\in\CM$ is invertible, then
  $E(T,r)=F(T^{-1},\frac1r)$. Indeed ``$\subseteq$'' follows by choosing
  $(\xi_n)_{n=1}^\infty$ as in Definition~\ref{E(T,r)} (i) and setting
  $\eta_n=T^n\xi_n$. Then \eqref{UH5} holds for $\eta =\xi$ and with $T$
  replaced by $T^{-1}$ and $r$ replaced by $\frac 1r$. The reverse
  inclusion follows by choosing $(\eta_n)_{n=1}^\infty$ as in
  Definition~\ref{E(T,r)} (ii) (with $T$ replaced by $T^{-1}$ and $r$
  replaced by $\frac 1r$) and setting $\xi_n= T^{-n}\eta_n$. 
\end{remark}

\vspace{.2cm}

\begin{lemma}\label{Lemma3.4} Let $T\in\CM$, and let $r>0$. Then
  \begin{itemize}
    \item[(i)] $\supp(\mu_T)\subseteq \overline{B(0,r)}$ iff $E(T,r)=\CH$.
    \item[(ii)] $\supp(\mu_T)\subseteq \C\setminus B(0,r)$ iff
    $F(T,r)=\CH$.
  \end{itemize}
\end{lemma}

\proof According to Lemma~\ref{Properties}~(d) we may assume that $\CM$
acts on the Hilbert space $\CH =L^2(\CM,\tau)$. For $a\in\CM$ we let
$[a]$ denote the corresponding element in $\CH$.

Throughout the proof we let $\mu_n= \mu_{(T^n)\cc T^n}=\mu_{T^n(T^n)\cc}$, and we let $\nu_n$ denote the push-forward
measure of $\mu_n$ under the map $t\mapsto t^{\frac 1n}$. Then, according
to Theorem~\ref{MSRI2.5}, $\nu_n\rightarrow \nu$ weakly in
$\Prob([0,\infty[)$, where $\nu$ is determined by
\[
  \nu([0,t^2])=\mu_T(\overline{B(0,t)}), \qquad (t>0).
\]

(i) Suppose that $\supp(\mu_T)\subseteq \overline{B(0,r)}$. 
Consider a fixed $s>r$. Since $\supp(\mu_T)\subseteq \overline{B(0,r)}$,
$\nu$ is supported on $[0,r^2]$, 
\[
\mu_n(]s^{2n},\infty[)= \mu_n(]s^2,\infty[)\rightarrow 0\;{\rm as}\;
n\rightarrow \infty.
\]
Define a sequence of projections in $\CM$, $(E_n)_{n=1}^\infty$, by
\[
E_n= 1_{[0,s^{2n}]}((T^n)\cc T^n).
\]
Note that for all $a\in\CM$,
\begin{eqnarray*}
\|(\unit-E_n)[a]\|_2^2&=& \tau(a\cc (\unit-E_n)a)\\
&\leq& \|a\|^2 \tau(\unit-E_n)\\
&=& \|a\|^2 \mu_n(]s^{2n},\infty[)\\
&\rightarrow& 0 \;{\rm as} \; n\rightarrow \infty.
\end{eqnarray*}
Since $\{[a]\,|\,a\in\CM\}$ is dense in $\CH$, this shows that
$E_n\rightarrow \unit$ in the strong operator topology, i.e. for all
$\xi\in\CH$,
\[
\lim_{n\rightarrow \infty}\|E_n\xi - \xi\|=0.
\]
Moreover, since $E_n (T^n)\cc T^n E_n \leq s^{2n}\unit$,
\[
(\|T^nE_n\xi\|)^{\frac 1n}\leq (\|T^n E_n\|\|\xi\|)^{\frac 1n} \leq s
\|\xi\|^{\frac 1n}
\]
for all $\xi\in\CH$, and hence $E(T,s)=\CH$. Since $s>r$ was arbitrary,
Lemma~\ref{Properties}~(b) now implies that $E(T,r)=\CH$. 

Next assume that $E(T,r)=\CH$. Let $\xi=[\unit]\in\CH$, and take a sequence
$(\xi_n)_{n=1}^\infty$ in $\CH$ such that
\[
    \lim_{n\rightarrow\infty}\|\xi_n-\xi\|=0\quad {\rm and}\quad
    \limsup_{n\rightarrow\infty} \|T^n\xi_n\|^{\frac 1n}\leq r.
\]
Since $\|\xi_n\|\rightarrow \|\xi\|=1$ as $n\rightarrow\infty$, we may as
well assume that $\|\xi_n\|=1$ for all $n\in\N$. Let $\mu_n'\in
\Prob([0,\infty[)$ denote the distribution of $(T^n)\cc T^n$ w.r.t. the
vector state $x\mapsto (x\xi_n, \xi_n)$. Let $s>r$, and take $s_1\in
(r,s)$. Then $\|T^n\xi_n\|\leq s_1^n$ eventually as $n\rightarrow \infty$. Hence for all large $n$,
\[
\int_0^\infty t\,\d\mu_n'(t) = ((T^n)\cc T^n \xi_n, \xi_n)\leq s_1^{2n},
\]
and therefore
\[
\mu_n'(]s^{2n}, \infty[)\leq s^{-2n}\int_{s^{2n}}^\infty t\,\d\mu_n'(t)\leq
\Big(\frac{s_1}{s}\Big)^{2n}.
\]
It follows that with $E_n= 1_{[0,s^{2n}]}((T^n)\cc T^n)$ as before,
\[
\|E_n\xi_n - \xi_n\|^2\leq \Big(\frac{s_1}{s}\Big)^{2n}
\]
for all large $n\in\N$, and hence
\[
\lim_{n\rightarrow\infty}\|E_n\xi_n - \xi_n\|=0.
\]
Since $\xi_n\rightarrow \xi$ as $n\rightarrow \infty$, it follows that
\[
\lim_{n\rightarrow\infty}\|E_n\xi_n - \xi\|=0
\]
as well. That is,
\[
\tau(\unit-E_n)= \|(\unit-E_n)\xi\|^2 \rightarrow 0 \;{\rm as}\;
n\rightarrow \infty,
\]
which is equivalent to saying that
\[
\mu_n(]s^{2n},\infty[)\rightarrow 0 \;{\rm as}\;
n\rightarrow \infty.
\]

We now have that for every $\phi\in
C_c^+(\R)$ with $\supp(\phi)\subseteq ]s^2,\infty[$,
\begin{eqnarray*}
  \int_0^\infty \phi(t)\,\d\nu(t) &= & \lim_{n\rightarrow \infty}
  \int_0^\infty \phi(t)\,\d\nu_n(t)\\
  &=& \leq \|\phi\|_\infty
  \limsup_{n\rightarrow\infty}\nu_n(]s^2,\infty[)\\
  & = & 0.
\end{eqnarray*}
Hence, $\supp(\nu)\subseteq [0,s^2]$, i.e. $\supp(\mu_T)\subseteq
\overline{B(0,s)}$ for all $s>r$, and we conclude that $\supp(\mu_T)\subseteq
\overline{B(0,r)}$.

(ii) At first assume that $\supp(\mu_T)\subseteq \C\setminus B(0,r)$. Let
$0<s<r$, and define $F_n\in \CM$ by
\[
F_n = 1_{[s^{2n},\infty[}(T^n(T^n)\cc).
\]
By an argument similar to the one given in the proof of (i), we find that
for all $\xi\in\CH$,
\[
\|F_n\xi -\xi\|\rightarrow 0\;{\rm as}\;n\rightarrow \infty.
\]
Now, consider a fixed $\xi\in\CH$, and with
\[
g_n(t)=\left\{\begin{array}{lll}
    0 &,& 0\leq t < s^{2n}\\
    \frac 1t &, & t\geq s^{2n}
    \end{array}\right.
\]
let
\[
\eta_n= (T\cc)^n g_n(T^n(T^n)\cc)\xi.
\]
Then, since $t\,g_n(t)= 1_{[s^{2n},\infty[}(t)$,
\[
T^n \eta_n = F_n\xi.
\]
Thus,
\[
\|T^n\eta_n - \xi\|\rightarrow 0 \;{\rm as} n\rightarrow \infty.
\]
Moreover, since $t\,g_n(t)^2=g_n(t)\leq \frac{1}{s^{2n}}$,
\[
\|\eta_n\|^2 =(g_n(T^n(T^n)\cc \xi, \xi)\leq \frac{1}{s^{2n}} \|\xi\|.
\]
Thus,
\[
\limsup_{n\rightarrow\infty}\|\eta_n\|^{\frac1n}\leq \frac1s,
\]
and this shows that $F(T,s)=\CH$. It now follows from
Lemma~\ref{Properties}~(b) that $F(T,r)=\CH$.

On the other hand, suppose that $F(T,r)=\CH$, and let $\xi =
[\unit]$. Choose a sequence $(\eta_n)_{n=1}^\infty$ in $\CH$ such that
\[
\lim_{n\rightarrow\infty}\|T^n\eta_n-\xi\|=0\quad {\rm and}\quad
    \limsup_{n\rightarrow\infty} \|\eta_n\|^{\frac 1n}\leq \frac1r.
\]
Since $\|\xi\|=1$, we may assume that $\|T^n\eta_n\|=1$ for all $n\in\N$.

Put $\xi_n=T^n\eta_n$, and let $\mu_n'\in\Prob([0,\infty[)$ denote the
distribution of $T^n(T^n)\cc$ w.r.t. the vector state $x\mapsto
(x\xi_n,\xi_n)$. Since for every $\eps >0$,
\[
\|(T^n(T^n)\cc + \eps\unit)^{-\frac12}T^n\|\leq 1,
\]
we have that
\begin{eqnarray*}
  \int_0^\infty \frac1t\,\d\mu_n'(t)&=& \sup_{\eps>0}\int_0^\infty
  \frac{1}{t+\eps}\,\d\mu_n'(t)\\
  &=& \sup_{\eps>0}\|(T^n(T^n)\cc + \eps\unit)^{-\frac12}\xi_n\|^2\\
  &=&  \sup_{\eps>0}\|(T^n(T^n)\cc + \eps\unit)^{-\frac12}T^n\eta_n\|^2\\
  &\leq & \|\eta_n\|^2.
\end{eqnarray*}

Now, let $s\in(0,r)$, and choose $s_1\in (s,r)$. Then
$\|\eta_n\|^{\frac1n}\leq \frac{1}{s_1}$ eventually as $n\rightarrow
\infty$, and hence
\[
\int_0^\infty \frac1t\,\d\mu_n'(t)\leq \frac{1}{s_1},
\]
eventually as $n\rightarrow \infty$. With $F_n=
1_{[s^{2n},\infty[}(T^n(T^n)\cc)$ as above we then have that
\begin{eqnarray*}
  \|F_n\xi_n - \xi_n\|^2 & = & \mu_n'([0,s^{2n}[)\\
  &\leq & \int_0^{2^n}\frac{s^{2n}}{t}\d\mu_n'(t)\\
  &\leq & \Big(\frac{s}{s_1}\Big)^{2n}\\
  &\rightarrow & 0 \;{\rm as}\;n\rightarrow\infty.
\end{eqnarray*}
Since $\xi_n\rightarrow \xi$, it follows that
\[
\|\unit-F_n\|_2^2 = \|\xi-F_n\xi\|^2\rightarrow 0\;{\rm
  as}\;n\rightarrow\infty,
\]
i.e.
\[
\lim_{n\rightarrow\infty}\nu_n([0,s[)=\lim_{n\rightarrow\infty}\mu_n([0,s^{2n}[)=0.
\]
Since $\nu_n\rightarrow \nu$ weakly, it follows that $\supp(\nu)\subseteq
[s^2,\infty[$, and hence $\supp(\mu_T)\subseteq \C\setminus B(0,s)$ for all
$s\in (0,r)$. We thus conclude that $\supp(\mu_T)\subseteq \C\setminus
B(0,r).$ $\endproof$

\vspace{.2cm}

%\begin{remark}\label{subspaces} Note that if $E_1$ and $E_2$ are $T$-invariant subspaces of
%  $\CH$ affiliated with $\CM$ and $E_1\subseteq E_2$, and if $P_i\in\CM$
%  denotes the projection onto $E_i$, then, according to Brown,
%  \[
%  \mu_{T|_{E_2}}= \tau_{P_2\CM P_2}(P_1)\mu_{T|_{E_1}}+ \tau_{P_2\CM
%  P_2}(P_2-P_1)\mu_{(P_2-P_1) T(P_2-P_1)}.
%  \]
%  Hence, $\supp(\mu_{T|_{E_1}})\subseteq \supp(\mu_{T|_{E_2}})$.
%\end{remark}

%\vspace{.2cm}

\begin{cor}\label{Cor3.5} Let $T\in\CM$, let $\lambda\in\C$, and let $r>0$. Then
  \begin{itemize}
    \item[(a)] If $E$ is any closed, $T$-invariant subspace affiliated with $\CM$, such that 
    $\supp(\mu_{T|_{E}})\subseteq \overline{B(\lambda,r)}$, then
    $E\subseteq E(T-\lambda\unit,r)$.
    \item[(b)] If $F$ is any closed, $T$-invariant subspace affiliated with $\CM$, such that 
    $\supp(\mu_{T|_{F}})\subseteq \C\setminus {B(\lambda,r)}$, then
    $F\subseteq F(T-\lambda\unit,r)$.
  \end{itemize}
\end{cor}

\proof (a) Suppose that  $E$ is a closed, $T$-invariant subspace affiliated
with $\CM$, such that $\supp(\mu_{T|_{E}})\subseteq
\overline{B(\lambda,r)}$. Then $E$ is invariant for $T-\lambda\unit$ as
well, and
\[
\supp(\mu_{T-\lambda\unit|_{E}})\subseteq
\overline{B(0,r)}.
\]
According to Lemma~\ref{Lemma3.4}, this implies that
\[
E(T-\lambda\unit|_{E},r)=E.
\]
It is easily seen that $E(T-\lambda\unit|_{E},r)\subseteq E(T-\lambda\unit
,r)$, and hence,
\[
E\subseteq E(T-\lambda\unit ,r).
\]

(b) follows in a similar way, so we leave out the proof of it. $\endproof$

\vspace{.2cm}

\begin{lemma}\label{perpendicular} For $T\in \CM$ and $s>r>0$ one has that
  $E(T,r)\,\bot\, F(T\cc,s)$.

\end{lemma}

\proof Let $\xi\in E(T,r)$, and $\eta\in F(T\cc,s)$. Take sequences
$(\xi_n)_{n=1}^\infty$ and $(\eta_n)_{n=1}^\infty$ in $\CH$, such that
$\xi_n\rightarrow \xi$, $(T\cc)^n \eta_n\rightarrow \eta$,
\[
\limsup_{n\rightarrow \infty}\|T^n\xi_n\|\leq r,
\]
and
\[
\limsup_{n\rightarrow \infty}\|\eta_n\|\leq \frac 1s.
\]

Choose positive numbers $r'$ and $s'$ such that $s>s'>r'>r$. Eventually as
$n\rightarrow \infty$, $\|T^n\xi_n\|\leq (r')^n$ and $\|\eta_n\|\leq
\frac{1}{(s')^n}$, whence
\begin{eqnarray*}
  |(\xi, \eta)| & = & \lim_{n\rightarrow\infty}|(\xi_n, (T\cc)^n\eta_n)|\\
  & = & \lim_{n\rightarrow\infty}|(T^n\xi_n, \eta_n)|\\
  & \leq & \limsup_{n\rightarrow\infty}\|T^n\xi_n\|\|\eta_n\|\\
  & \leq & \limsup_{n\rightarrow\infty}\Big(\frac{r'}{s'}\Big)^n\\
  & = & 0.
\end{eqnarray*}
That is, $\xi\bot \eta$. $\endproof$

%% file: tispfreeprob.tex
\section{Some results from free probability theory}

In order to proceed, we will need a number of results which are based on free probability theory. Recall from \cite{VDN} that a $C\cc$--probability space is a pair $(\CA,\varphi)$ where $\CA$ is a unital $C\cc$--algebra and $\varphi$ is a state on $\CA$. We say that $(\CA,\varphi)$ is a $W\cc$--probability space if $\CA$ is a von Neumann algebra and if $\varphi$ is a normal state on $\CA$. 

Let $\CH$ be a Hilbert space. The full Foch--space over $\CH$, $\CT(\CH)$, is the Hilbert space
\[
\CT(\CH)=\C\Omega\oplus\Big(\bigoplus_{n=1}^\infty \CH^{\otimes n}\Big),
\]
where $\Omega$, the vacuum vector, is a unit vector in the one--dimensional Hilbert space $\C\Omega$. $\Omega$ induces the vector state $\omega$ on $B(\CT(\CH))$. 

To a given orthonormal set $(e_i)_{i\in I}$ in $\CH$ we associate the creation operators $(\ell_i)_{i\in I}$ in $B(\CT(\CH))$ given by
\[
\ell_i\Omega = e_i \qquad {\rm and} \qquad \ell_i\xi=e_i\otimes\xi, \quad \xi\in\CH^{\otimes n}.
\]  
Then $\ell_i\cc\ell_j=\delta_{ij}\unit$ and
\begin{equation}\label{4.1}
\omega(\ell_i^m(\ell_i\cc)^n)=
\left\{\begin{array}{ll}
1, & (m,n)=(0,0),\\
0, & (m,n)\in\N_0^2\setminus\{(0,0)\}.
\end{array}\right.
\end{equation}
Moreover, the $\ast$--algebras $\CA_i={\rm alg}(\ell_i,\ell_i\cc)$, $i\in I$, are free in $B(\CT(\CH))$ w.r.t. $\omega$ (cf. \cite{VDN}). 

A family of creation operators $(\ell_i)_{i\in I}$ obtained from an orthonormal set $(e_i)_{i\in I}$ as above is called a \emph{free family of creation operators}. 

\begin{lemma}\label{lemma4.1} Let $(\CA,\varphi)$ be a $C\cc$--probability space, and let $(s_i)_{i\in I}$ be operators in $\CA$ which satisfy
\begin{itemize}
\item[(i)] $s_i\cc s_j = \delta_{ij}\unit$, $i,j\in I$,
\item[(ii)] $\varphi(s_is_i\cc)=0$, $i\in I$.
\end{itemize}
Then $(s_i)_{i\in I}$ has the same $\ast$--distribution as a free family of creation operators $(\ell_i)_{i\in I}$.
\end{lemma}

\proof Let $\CB_i={\rm alg}(s_i,s_i\cc)\subseteq \CA$, $i\in I$. Since $s_i\cc s_i=\unit$, 
\begin{equation}\label{4.2}
\CB_i={\rm span}\{s_i^m(s_i\cc)^n\,|\,(m,n)\in \N_0^2\}.
\end{equation}
By (ii) and the Schwarz inequality
\[
|\varphi(ab)|\leq \varphi(aa\cc)^\frac12\varphi(b\cc b)^\frac12, \qquad a,b\in\CA,
\]
we have that
\begin{equation}\label{4.3}
\varphi(s_i^m(s_i\cc)^n) =\left\{\begin{array}{ll}
1, & (m,n)=(0,0),\\
0, & (m,n)\in\N_0^2\setminus\{(0,0)\}.
\end{array}\right.
\end{equation}
Let $(\ell_i)_{i\in I}$ be the set of creation operators coming from an orthonormal set $(e_i)_{i\in I}$ as described above. Then by \eqref{4.1}, $\ell_i$ and $s_i$ have the same $\ast$--distribution for all $i\in I$. Hence, in order to show that the families $(\ell_i)_{i\in I}$ and $(s_i)_{i\in I}$ have the same $\ast$--distributions, it suffices to show that $(\CB_i)_{i\in I}$ are free subalgebras in $(\CA,\varphi)$.

Put 
\[
\CB_i^0=\{b\in\CB_i\,|\,\varphi(b)=0\}.
\]
It follows from \eqref{4.2} and \eqref{4.3} that
\[
\CB_i^0={\rm span}\{s_i(s_i\cc)^n\,|\,(m,n)\in \N_0^2\setminus\{(0,0)\}\}.
\]
Thus, in order to prove that $(\CB_i)_{i\in I}$ are free, it suffices to show that
\begin{equation}\label{4.4}
\varphi(s_{i_1}^{m_1}(s_{i_1}\cc)^{n_1}s_{i_2}^{m_2}(s_{i_2}\cc)^{n_2}\cdots s_{i_k}^{m_k}(s_{i_k}\cc)^{n_k})=0
\end{equation}
for all $k\in\N$, all $(m_1,n_1),\ldots (m_k,n_k)\in \N_0^2\setminus\{(0,0)\}$ and all $i_1,\ldots, i_k\in I$ such that $i_1\neq i_2\neq \cdots \neq i_{k-1}\neq i_k$. 

Given $k\in\N$, $(m_1,n_1),\ldots (m_k,n_k)\in \N_0^2\setminus\{(0,0)\}$ and $i_1,\ldots, i_k\in I$ such that $i_1\neq i_2\neq \cdots \neq i_{k-1}\neq i_k$, assume that 
\[
\varphi(s_{i_1}^{m_1}(s_{i_1}\cc)^{n_1}s_{i_2}^{m_2}(s_{i_2}\cc)^{n_2}\cdots s_{i_k}^{m_k}(s_{i_k}\cc)^{n_k})\neq 0.
\]
By (ii) and the Schwarz inequality, 
\begin{equation}\label{4.5}
\varphi(as_i\cc)=0, \qquad a\in\CA, \,\,i\in I.
\end{equation}
Hence $n_k=0$, which implies that $m_k\neq 0$. By (i), $s_{i_{k-1}}\cc s_{i_k}=0$, and therefore $n_{k-1}$ must be zero. Then $m_{k-1}\neq 0$, and continuing this way we find that $n_1=n_2=\cdots =n_k=0$, $m_1, m_2,\ldots,m_k\geq 1$ and
\[
\varphi(s_{i_1}^{m_1}s_{i_2}^{m_2}\cdots s_{i_k}^{m_k})\neq 0.
\]
However, by \eqref{4.5}, $\varphi(s_ia)=\overline{\varphi(a\cc s_i\cc)}=0$ for $i\in I$ and $a\in \CA$. Thus, since $m_1\geq 1$, $\varphi(s_{i_1}^{m_1}s_{i_2}^{m_2}\cdots s_{i_k}^{m_k})=0$ and we have reached a contradiction. This proves \eqref{4.4} and completes the proof of the lemma. $\endproof$

\vspace{.2cm}

Next we will apply the following result due to Shlyakhtenko:

\begin{lemma}\cite[Corollary~2.5]{Sh}\label{lemma4.2} Let $(\CA,\varphi)$ be a $C\cc$--probability space, $\CB\subseteq \CA$ a $C\cc$--subalgebra of $\CA$ with $\unit\in\CB$, and let $\ell_1,\ldots,\ell_n\in\CA$. Suppose the following three conditions hold: 
\begin{itemize}
\item[(1)]$\ell_1,\ldots,\ell_n$ have the same $\ast$--distribution w.r.t. $\varphi$ as a free family of $n$ creation operators.
\item[(2)] For all $p,q\geq 0$ with $p+q>0$, for all $b_1,\ldots,b_p, b_1',\ldots, b_{q+1}'\in\CB$ and for all $1\leq i_1,\ldots, i_p,j_1,\ldots, j_q\leq n$,
\[
\varphi(b_1\ell_{i_1}b_2\cdots\ell_{i_p}b_1'\ell_{j_1}\cc b_2'\cdots \ell_{j_q}\cc b_{q+1}')=0.
\]
\item[(3)] The GNS representation of $C\cc(\CB,\ell_1,\ldots,\ell_n)$ associated with $\varphi$ is faithful.
\end{itemize}
Then the following two conditions are equivalent:
\begin{itemize}
\item[(a)] $\CB$ and $C\cc(\ell_1,\ldots,\ell_n)$ are free in $(\CA,\varphi)$
\item[(b)] $\ell_i\cc b\ell_j=\delta_{ij}\varphi(b)\unit$ for $b\in\CB$, $1\leq i,j\leq n$
\end{itemize}
\end{lemma}

\vspace{.2cm}

Combining Lemma~\ref{lemma4.1} and Lemma~\ref{lemma4.2}, we get:

\begin{lemma}\label{lemma4.3} Let $(\CA,\varphi)$ be a $C\cc$--probability space, $\CB\subseteq \CA$ a $C\cc$--subalgebra of $\CA$ with $\unit\in\CB$, and let $\ell_1,\ldots,\ell_n\in\CA$. Suppose the GNS representation of $C\cc(\CB,\ell_1,\ldots,\ell_n)$ associated to $\varphi$ is faithful. 
\begin{itemize}
\item[(i)] If the following two conditions hold:
\begin{itemize}
\item[(a1)]$\ell_1,\ldots,\ell_n$ have the same $\ast$--distribution w.r.t. $\varphi$ as a free family of $n$ creation operators
\item[(a2)] $\CB$ and $C\cc(\ell_1,\ldots,\ell_n)$ are free in $(\CA,\varphi)$
\end{itemize}
Then the following two conditions also hold:
\begin{itemize}
\item[(b1)] $\ell_i\cc b\ell_j=\delta_{ij}\varphi(b)\unit$ for $b\in\CB$, $1\leq i,j\leq n$
\item[(b2)] $\varphi(b\cc\ell_i\ell_i\cc b)=0$ for $b\in\CB$, $1\leq i\leq n$
\end{itemize}
\item[(ii)] Conversely, (b1) together with (b2) implies (a1) and (a2).
\end{itemize}
\end{lemma}

\proof (i) Suppose (a1) and (a2) hold. Then (1) and (a) in Lemma~\ref{lemma4.2} are fulfilled. Moreover, by assumption, (3) in Lemma~\ref{lemma4.2} also holds. We prove that (b2) holds: Let $\CC=C\cc(\ell_1,\ldots, \ell_n)$ and let
\[
\CB^0=\{b\in\CB\,|\,\varphi(b)=0\} \quad {\rm and} \quad \CC^0=\{c\in\CC\,|\,\varphi(c)=0\}.
\]
For $b\in\CB$, let $b^0=b-\varphi(b)\unit$. Then for $1\leq i\leq n$,
\begin{eqnarray*}
\varphi(b\cc\ell_i\ell_i\cc b)&=& \varphi(((b\cc)^0+\varphi(b\cc)\unit)\ell_i\ell_i\cc(b^0+\varphi(b)\unit))\\
&=& \varphi((b\cc)^0\ell_i\ell_i\cc b^0)+\overline{\varphi(b)}\varphi(\ell_i\ell_i\cc b^0)+ \varphi(b)\varphi((b\cc)^0\ell_i\ell_i\cc)+ |\varphi(b)|^2\varphi(\ell_i\ell_i\cc).
\end{eqnarray*}
According to (a1), $\varphi(\ell_i\ell_i\cc)=0$, that is $\ell_i\ell_i\cc\in\CC^0$. Moreover, $b^0,(b\cc)^0\in\CB^0$ so by the freeness assumption (a2),
\[
\varphi((b\cc)^0\ell_i\ell_i\cc b^0)=\varphi(\ell_i\ell_i\cc b^0)= \varphi((b\cc)^0\ell_i\ell_i\cc)= \varphi(\ell_i\ell_i\cc)=0.
\]
This proves (b2). We will now show that (2) of Lemma~\ref{lemma4.2} holds: By the Schwarz inequality and by (b2),
\begin{equation}\label{4.6}
\varphi(b\ell_i a)=0, \qquad b\in\CB,\, a\in\CA,
\end{equation}
because
\[
|\varphi(b\ell_i a)|\leq \varphi(b\ell_i(b\ell_i)\cc)^\frac12 \varphi(a\cc a)^\frac12=0.
\]
Moreover, by complex conjugation of \eqref{4.6}, it follows that
\begin{equation}\label{4.7}
\varphi(a\ell_i\cc b)=0, \qquad b\in\CB,\,a\in\CA.
\end{equation}
This proves that (2) (in addition to (1), (3) and (a)) of Lemma~\ref{lemma4.2} holds, and thus (b1)=(b) holds. This proves (i).

(ii) Assume now that (b1) and (b2) hold. Then (b) of Lemma~\ref{lemma4.2} holds. Moreover, by assumption, condition (3) of Lemma~\ref{lemma4.2} holds. From the proof of (i) we know that (b2) implies \eqref{4.6} and \eqref{4.7}, and the latter two imply (2) of Lemma~\ref{lemma4.2}. 

By (b1) and (b2), $\ell_1,\ldots,\ell_n$ satisfy the conditions in Lemma~\ref{lemma4.1}. Hence, $(\ell_1,\ldots,\ell_n)$ has the same $\ast$--distribution as a set of $n$ free creation operators. This proves (a1) as well as (1) in Lemma~\ref{lemma4.2}. Thus, (1), (2), (3) and (b) of Lemma~\ref{lemma4.2} hold, and then by Lemma~\ref{lemma4.2}, (a2)=(a) holds. This completes the proof of (ii). $\endproof$

\vspace{.2cm}

Recall from \cite{VDN} that when $(\CA_1,\varphi_1)$ and $(\CA_2,\varphi_2)$ are two $W\cc$--probability spaces, one can define a \emph{reduced} free product,
\[
(\CA,\varphi)= (\CA_1,\varphi_1)\ast (\CA_2,\varphi_2),
\]
which is a von Neumann algebra $\CA$ realized on the free product of Hilbert spaces with distinguished unit vector,
\[
(\CH,\xi)=(\CH_{\varphi_1},\xi_{\varphi_1})\ast (\CH_{\varphi_2},\xi_{\varphi_2}),
\]
where $\CH_{\varphi_i}$ with distinguished unit vector $\xi_{\varphi_i}$ is the GNS Hilbert space in the representation $\pi_{\varphi_i}$ of $\CA_i$. $\varphi$ is the vector state on $\CA$ induced by the unit vector $\xi\in\CH$ which is cyclic for $\CA$. In particular, the GNS representation $\pi_\varphi$ of $\CA$ is one--to--one and can be identified with the inclusion map $\CA\hookrightarrow B(\CH)$. If $\pi_{\varphi_1}$ and $\pi_{\varphi_2}$ are one--to--one, then $\CA_1$ and $\CA_2$ are naturally embedded in $\CA$ as a free pair of subalgebras w.r.t. $\varphi$, and $\varphi_i=\varphi|_{\CA_i}$, $i=1,2$. 

\vspace{.2cm}

Let $(\CA,\varphi)$ be a $C\cc$--probability space. Following the notation of \cite{V1} (or \cite{VDN}), we say that $s_1,\ldots,s_n$ in $(\CA,\varphi)$ form a \emph{semicircular family} in $(\CA,\varphi)$ if $s_1,\ldots,s_n$ are free self--adjoint elements with
\[
\varphi(s_i^p)=\frac{1}{2\pi}\int_{-2}^2 t^p\sqrt{4-t^2}\,\d t, \qquad p\in\N.
\]
Elements $x_1,\ldots,x_n$ in $(\CA,\varphi)$ are said to form a \emph{circular family} if there is a semicircular family $(s_1,\ldots, s_{2n})$ in $(\CA,\varphi)$ such that
\[
x_j=\frac{s_j+\i s_{n+j}\cc}{\sqrt2}, \qquad 1\leq j\leq n.
\]

\begin{prop}\label{circular} Let $(\CN,\tau)$ be a finite von Neumann algebra with a faithful, normal, tracial state $\tau$, and let $\CM\subseteq \CN$ be a von Neumann subalgebra of $\CN$. Let $n\in\N$, and let $u= (u_{ij})_{i,j=1}^n \in
  \CU(M_n(\CM))$. Moreover, suppose that $(x_1, \ldots, x_n)$ is a
  circular family in $\CN$ which is $\ast$-free from $\CM$. Then with
  $y_1, \ldots, y_n\in \CN$ given by
  \begin{equation}
    y_i = \sum_{j=1}^n u_{ij}x_j, \qquad (1\leq i\leq n),
  \end{equation}
  $(y_1, \ldots, y_n)$ is a circular family which is $\ast$-free from $\CM$.
\end{prop}

\proof We may without loss of generality assume that 
\[
\CN=\CM\ast W\cc(x_1,\ldots,x_n).
\]
By \cite{V1} (or \cite{VDN}), $W\cc(x_1,\ldots,x_n)\cong L(\F_{2n})$, the von Neumann algebra of the free group on $2n$ generators. Moreover, we can realize the circular family as follows:

Take
an orthonormal set $e_1, e_1', \ldots, e_n, e_n'$ in a Hilbert space $\CH$, and let $\ell_i= \ell(e_i)$ and $\ell_i'= \ell(e_i')$
be the corresponding creation operators on $\CT(\CH)$. Moreover, let $\omega$ be the state on $B(\CT(\CH))$ given by the vacuum
vector $\Omega$. Then
\[
x_i = \ell_i + (\ell_i')\cc, \qquad 1\leq i\leq n,
\]
is a circular family in $(B(\CT(\CH)),\omega)$ and the restriction $\omega'$ of $\omega$ to $W\cc(x_1,\ldots,x_n)$ is a faithful, normal, tracial state. Furthermore, $\Omega$ is cyclic for $W\cc(x_1,\ldots,x_n)$ and we can therefore identify the GNS representation $\pi_{\omega'}$ with the inclusion map $W\cc(x_1,\ldots,x_n)\hookrightarrow B(\CT(\CH))$. The GNS vector $\xi_{\omega'}$ then corresponds to $\Omega$. Hence, the reduced free products
\[
(\CN,\tau)=(\CM,\tau|_\CM)\ast(W\cc(x_1,\ldots,x_n),\omega')
\]
and
\[
(\tilde\CN,\varphi)=(\CM,\tau|_\CM)\ast(B(\CT(\CH)),\omega)
\]
are both realized on the free product of Hilbert spaces,
\[
(\CH,\xi)=(L^2(\CM,\tau|_\CM),\unit)\ast(\CT(\CH),\Omega).
\]
We may therefore regard $\CN$ as a subalgebra of $\tilde\CN$ such that $\varphi|_{\CN}=\tau$. We will now apply Lemma~\ref{lemma4.3}~(i) to $\CA=\tilde\CN$, $\CB=\CM$ and the family of creation operators $(\ell_1,\ldots,\ell_n,\ell_1',\ldots,\ell_n')$. Since $\xi$ is cyclic for $\tilde\CN = W\cc(\CM,\ell_1,\ldots,\ell_n,\ell_1',\ldots,\ell_n')$, the GNS representation of $\tilde\CN$ associated to $\varphi$ is one--to--one, which implies that the GNS representation of the restriction of $\varphi$ to $C\cc(\CM,\ell_1,\ldots,\ell_n,\ell_1',\ldots,\ell_n')$ is also one--to--one. Hence, the assumptions of Lemma~\ref{lemma4.3} and (a1) and (a2) in Lemma~\ref{lemma4.3} are fulfilled. It then follows that (b1) and (b2) hold, that is
\begin{itemize}
  \item[(i)] $\forall x\in \CM \; \forall i, j\in\{1,\ldots, n\}$:
    \begin{eqnarray*}
      \ell_i\cc x \ell_j &=& \delta_{ij}\cdot\tau(x)\unit,\\
       (\ell_i')\cc x \ell_j' &=& \delta_{ij}\cdot\tau(x)\unit,\\
      \ell_i\cc x \ell_j' &=&(\ell_i')\cc x \ell_j = 0.
    \end{eqnarray*}
  \item[(ii)] $\forall x\in \CM\; \forall i\in\{1,\ldots, n\}$:
    \[
    \varphi(x\cc \ell_i \ell_i\cc x)= \varphi(x\cc \ell_i' (\ell_i')\cc x) =0.
    \]
\end{itemize}

Put
\begin{equation}\label{4.8}
s_i = \sum_{j=1}^n u_{ij}\ell_j, \qquad (1\leq i\leq n)
\end{equation}
and
\begin{equation}\label{4.9}
s_i' = \sum_{j=1}^n \ell_j' u_{ij}\cc , \qquad (1\leq i\leq n).
\end{equation}
Then
\begin{equation}\label{4.10}
y_i = s_i + (s_i')\cc, \qquad (1\leq i\leq n).
\end{equation}
We will now apply Lemma~\ref{lemma4.3}~(ii) to $\CA=\tilde\CN$, $\CB=\CM$ and the family of operators $(s_1,\ldots,s_n,s_1',\ldots, s_n')$. Since $u=(u_{ij})_{i,j=1}^n$ is unitary, \eqref{4.8} and \eqref{4.9} imply that
\begin{eqnarray*}
\ell_j&=&\sum_{i=1}^n u_{ij}\cc s_i, \qquad 1\leq j\leq n,\\
\ell_j'&=&\sum_{i=1}^n s_i'u_{ij}\cc, \qquad 1\leq j\leq n.
\end{eqnarray*}
Hence, 
\[
C\cc(\CM,s_1,\ldots,s_n,s_1',\ldots,s_n')=C\cc(\CM,\ell_1,\ldots,\ell_n,\ell_1',\ldots,\ell_n'),
\]
and the conditions of Lemma~\ref{lemma4.3} are fulfilled by $s_1,\ldots,s_n,s_1',\ldots,s_n'$ as well. We will now make sure that (b1) and (b2) in Lemma~\ref{lemma4.3} are fulfilled, i.e. that
\begin{itemize}
  \item[(i')] $\forall x\in \CM \; \forall i, j\in\{1,\ldots, n\}$:
    \begin{eqnarray*}
      s_i\cc x s_j &=& \delta_{ij}\cdot\tau(x)\unit,\\
       (s_i')\cc x s_j' &=& \delta_{ij}\cdot\tau(x)\unit,\\
      s_i\cc x s_j' &=&(s_i')\cc x s_j = 0.
    \end{eqnarray*}
  \item[(ii')] $\forall x\in \CM\; \forall i\in\{1,\ldots, n\}$:
    \[
    \varphi(x\cc s_i s_i\cc x)= \varphi(x\cc s_i' (s_i')\cc x) =0.
    \]
\end{itemize}
If (i') and (ii') hold, then according to Lemma~\ref{lemma4.3}~(ii), $(s_1,\ldots,s_n,s_1',\ldots,s_n')$ has the same $\ast$--distribution w.r.t. $\varphi$ as free creation operators, and $\CM$ and $C\cc(s_1,\ldots,s_n,s_1',\ldots,s_n')$ are free in $(\tilde\CN,\varphi)$. Hence, $y_i=s_i+(s_i')\cc$, $1\leq i\leq n$, is a circular family which is $\ast$--free from $\CM$. 

To prove (i') and (ii'), let $x\in\CM$, and let $1\leq i,j\leq n$. Then (i)
implies that
\begin{eqnarray*}
s_i\cc x s_j &=& \sum_{p,q=1}^n \ell_p\cc u_{ip}\cc x u_{jq}\ell_q\\
& = & \sum_{p,q=1}^n \delta_{pq}\varphi(u_{ip}\cc x u_{jq})\unit\\
& = & \sum_{p=1}^n \varphi(x u_{jp}(u\cc)_{pi})\unit\\
& = & \varphi(x (uu\cc)_{ji})\unit \\
& = & \delta_{ij}\cdot\tau(x)\unit,
\end{eqnarray*}
and similarly,
\[
(s_i')\cc x s_j' = \delta_{ij}\cdot\tau(x)\unit.
\]
Moreover, according to (i),
\[
s_i\cc x s_j' = \sum_{p,q=1}^n \ell_p\cc u_{ip}\cc x \ell_q' u_{jq}\cc =
\sum_{p,q=1}^n 0\cdot u_{jq}\cc =0,
\]
and hence
\[
(s_i')\cc x s_j = (s_j\cc x\cc s_i')\cc = 0.
\]
This proves (i'). Next we prove that (ii') holds, making use of (i) and (ii):

Let $x\in\CM$. Then
\begin{equation}\label{4.11}
\varphi(x\cc s_is_i\cc x)= \sum_{p,q=1}^n\varphi(x\cc u_{ip}\ell_p\ell_q\cc u_{iq}\cc x),
\end{equation}
and
\begin{equation}\label{4.12}
\varphi(x\cc s_i'(s_i')\cc x)= \sum_{p,q=1}^n\varphi(x\cc \ell_p' u_{ip}\cc (\ell_q')\cc x).
\end{equation}
Since (ii) holds, we get, using the Schwarz inequality, that
\begin{equation}\label{4.13}
\varphi(y\cc \ell_i a)=\varphi(y\cc\ell_i'a)=0
\end{equation}
for $y\in\CM$, $a\in\tilde\CN$ and $1\leq i\leq n$. Applying this to $y=x\cc u_{ip}$ and $a=\ell_p\cc u_{iq}\cc x$, we get, using \eqref{4.11}, that
\[
\varphi(x\cc s_is_i\cc x)=0.
\]
Similarly, applying \eqref{4.13} to $y=x$ and $a=u_{ip}\cc u_{iq}(\ell_q')\cc x$, we get from \eqref{4.12} that
\[ 
\varphi(x\cc s_i'(s_i')\cc x)=0.
\]
This proves (ii'). $\endproof$

\vspace{.2cm}

Let $\CM$ be a II$_1$--factor, and let $T\in\CM$. We regard $\CM$ as a subfactor of $\CN=\CM\ast L(\F_4)$ with tracial state $\tau$, and we choose a circular system $\{x,y\}$ that generates $L(\F_4)$ and which therefore is free from $\CM$. Then by \cite[Theorem~5.2]{HS}, the unbounded operator $z=xy^{-1}$ is in $L^p(\CN,\tau)$, $0<p<1$. Thus, for all $a\in(0,\infty)$, $T+az\in L^p(\CN,\tau)$, $0<p<1$, and therefore $T+az$ has a well--defined Fuglede--Kadison determinant and a well--defined Brown measure (cf. \cite[Section~2]{HS}). 

\begin{prop}\label{T_a}
  Let $T\in\CM$, and let $z=xy^{-1}$ as above. Then for each $a>0$,
  \begin{equation}
    \Delta(T+az)= \Delta(T\cc T+a^2\unit)^\frac12, \qquad (a>0).
  \end{equation}
\end{prop}

\proof Take a unitary $u\in\CM$ such that $T=u|T|$, and for fixed $a>0$
put
\begin{equation*}
  w=(T\cc T+a^2\unit)^{-\frac12}(|T|y+au\cc x).
\end{equation*}
Then
\begin{equation*}
T+az = u (T\cc T+a^2\unit)^{\frac12} w y^{-1},
\end{equation*}
and hence
\begin{eqnarray*}
\Delta(T+az)& = &\Delta(u) \Delta((T\cc
T+a^2\unit)^{\frac12})\Delta(w)\Delta( y^{-1})\nonumber\\
&=& \Delta(T\cc T+a^2\unit)^{\frac12}\Delta(w)\Delta( y)^{-1}.
\end{eqnarray*}

It follows that if $w$ is circular, then $\Delta(T+az)= \Delta(T\cc
T+a^2\unit)^\frac12$. To see that $w$ is circular, define $v\in
\CU(M_2(\CM))$ by
\begin{equation*}
v = \begin{pmatrix}
    (T\cc T+a^2\unit)^{-\frac12}|T| & a(T\cc T+a^2\unit)^{-\frac12}u\cc \\
    a(T\cc T+a^2\unit)^{-\frac12} & -(T\cc T+a^2\unit)^{-\frac12}|T|u\cc
    \end{pmatrix},
\end{equation*}
and define $x', y'\in\CM$ by
\begin{equation*}
\begin{pmatrix}
  y' \\ x'
\end{pmatrix} = v \begin{pmatrix}
  y \\ x
\end{pmatrix}.
\end{equation*}
Then, according to Proposition~\ref{circular}, $\{x', y'\}$ is a circular
system (which is $\ast$-free from $\CM$). Since $w=y'$, $w$ is then circular, and hence $\Delta(w)=\Delta(y)$. $\endproof$

\vspace{.2cm}

\begin{cor}\label{cor4.6}The Brown measure of $T+az$ is given by
\begin{equation}\label{mu_a}
\d\mu_{T+az}(\lambda)=\frac{a^2}{\pi}\tau\big((T(\lambda)\cc T(\lambda)+a^2\unit)^{-1}(T(\lambda)T(\lambda)\cc + a^2\unit)^{-1}\big)\,\d\lambda_1\,\d\lambda_2,
\end{equation}
where $\lambda_1=\re\lambda$, $\lambda_2=\im\lambda$, and $T(\lambda)=T-\lambda\unit$.
\end{cor}

\proof By \cite[Definition~2.13]{HS} and Proposition~\ref{T_a}, 
\begin{eqnarray*}
\d\mu_{T+az}(\lambda)&=&\frac{1}{2\pi}\nabla^2\big(\log\Delta(T+az-\lambda\unit)\big)\,\d\lambda_1\,\d\lambda_2\\
&=& \frac{1}{2\pi}\nabla^2\big(\textstyle{\frac12}\log\Delta\big((T-\lambda\unit)\cc(T-\lambda\unit)+a^2\unit\big)\big)\,\d\lambda_1\,\d\lambda_2,
\end{eqnarray*}
where the Laplacian $\nabla^2=\frac{\partial^2}{\partial\lambda_1^2}+\frac{\partial^2}{\partial\lambda_2^2}$ is taken in the distribution sense. By the proof of \cite[Lemma~2.8]{HS}, 
\[
L_a(\lambda):=\textstyle{\frac12}\log\Delta\big((T-\lambda\unit)\cc(T-\lambda\unit)+a^2\unit\big)
\]
is a $C^2$--function in $(\lambda_1,\lambda_2)$ and its Laplacian is given by
\[
(\nabla^2L_a)(\lambda)=2a^2\tau\big((T(\lambda)\cc T(\lambda)+a^2)^{-1}(T(\lambda)T(\lambda)\cc+a^2)^{-1}\big)
\]
(replace $(a,b,\eps)$ by $(T,\unit,a^2)$ in the function $g_\eps$ in the proof of \cite[Lemma~2.8]{HS}). This proves \eqref{mu_a}. $\endproof$

\vspace{.2cm}

\begin{remark} In \cite[Theorem~5.2]{HS} it was shown that the Brown measure of $z$ is given by
\[
\d\mu_z(\lambda)=\frac{1}{\pi(1+|\lambda|^2)}\,\d\lambda_1\,\d\lambda_2.
\]
Note that this can also be obtained as a special case of \eqref{mu_a} with $T=0$ and $a=1$.
\end{remark}

\vspace{.2cm}

\begin{cor}\label{cor4.8} Let $T$ and $T+az$ be as in Proposition~\ref{T_a}. Then
\[
\mu_{T+az}\overset{w\cc}{\rightarrow}\mu_T \quad {\rm as} \quad a\rightarrow 0+.
\]
\end{cor}

\proof It suffices to show that for all $\varphi\in C_c^\infty(\C)$,
\begin{equation}\label{4.18}
\int_\C\varphi\,\d\mu_{T+az}\rightarrow \int_\C\varphi\,\d\mu_T \quad {\rm as} \quad a\rightarrow 0+.
\end{equation} 
As in the proof of Corollary~\ref{cor4.6}, put
\begin{equation}
\begin{split}
L_a(\lambda)& = \log\Delta(T+az-\lambda\unit)\\
& = \frac12 \log\Delta\big((T-\lambda\unit)\cc(T-\lambda\unit)+a^2\unit\big),\label{4.19}
\end{split}
\end{equation}
and let
\begin{equation}\label{4.20}
L(\lambda)=\log\Delta(T-\lambda\unit)=\frac12\log\Delta((T-\lambda\unit)\cc(T-\lambda\unit)).
\end{equation}
Then 
\[
\d\mu_{T+az}(\lambda)=\frac{1}{2\pi}\nabla^2L_a(\lambda)\,\d\lambda_1\,\d\lambda_2,
\]
and
\[
\d\mu_T(\lambda)=\frac{1}{2\pi}\nabla^2L(\lambda)\,\d\lambda_1\,\d\lambda_2.
\]
Hence, for all $\varphi\in C_c^\infty(\C)$,
\begin{equation}\label{4.21}
\int_\C\varphi\,\d\mu_{T+az}=\frac{1}{2\pi}\int_\C(\nabla^2\varphi)(\lambda)L_a(\lambda)\,\d\lambda_1\,\d\lambda_2,
\end{equation}
and
\begin{equation}\label{4.22}
\int_\C\varphi\,\d\mu_{T}=\frac{1}{2\pi}\int_\C(\nabla^2\varphi)(\lambda)L(\lambda)\,\d\lambda_1\,\d\lambda_2.
\end{equation}
By \eqref{4.19} and \eqref{4.20}, $L_a(\lambda) \searrow L(\lambda)$ as $a\searrow 0$. Since $L$ and $L_a$ are subharmonic and hence locally integrable, \eqref{4.18} then follows from \eqref{4.21}, \eqref{4.22} and Lebesgue's dominated convergence theorem. $\endproof$

%% file: Lipschitz.tex
\section{A certain Lipschitz condition}

Consider a von Neumann algebra $\CN$ equipped with a
faithful normal tracial state $\tau$, and assume that $\CN$ contains the von
Neumann algebra $\CM$ as a sub-algebra and a circular system $\{x,y\}$
which is $\ast$-free from $\CM$. We let $z=xy^{-1}\in L^p(\CN, \tau)$,
$(0<p<1)$. By \cite[Theorem~5.2 and Theorem~5.4]{HS}, 
\begin{equation}\label{6.-1}
z,z^{-1}\in L^p(\CN,\tau) \quad {\rm for \,\,\, all} \quad p\in (0,1).
\end{equation}
Moreover, if $p\in (0,\frac23)$ and $\lambda\in\C$, then $z^2,z^{-2},(z^2-\lambda\unit)^{-1}\in L^p(\CN,\tau)$, and
\begin{equation}\label{6.0}
\|(z^2-\lambda\unit)^{-1}\|_p\leq \|z^{-2}\|_p=\|z^2\|_p<\infty.
\end{equation}
As an application of \eqref{6.-1}, \eqref{6.0}, and Proposition~\ref{circular}, in this section we prove:

\begin{thm}\label{Lipschitz} For every $T\in\CM$ and every $p\in(0,1)$,
  $T+z$ has an inverse $(T+z)^{-1}\in L^p(\CN, \tau)$. Moreover, for each $p\in(0, \frac 23)$
  there is a constant $C_p^{(1)} >0$ such that
  \begin{equation}\label{eqLip}
    \forall S, T\in\CM:\quad \|(S+z)^{-1}-(T+z)^{-1}\|_p \leq C_p^{(1)} \|S-T\|.
  \end{equation}
\end{thm}

\vspace{.2cm}

Theorem~\ref{Lipschitz} is a consequence of the following three results,
which we prove later on.

\begin{lemma}\label{ABC} Let $S,T\in \CM$. Then
  \begin{eqnarray}
    \|(\unit +S\cc S)^{-\frac 12}S\cc - (\unit +T\cc T)^{-\frac 12}T\cc \|&\leq&
    \textstyle{\frac 54}\, \|S-T\|,\label{eq4-20}
  \end{eqnarray}
  and
   \begin{eqnarray}
     \|(\unit +S\cc S)^{-\frac 12} - (\unit+ T\cc T)^{-\frac 12}\|&\leq&
    \textstyle{\frac 2\pi}\, \|S-T\|.\label{eq4-21}
  \end{eqnarray}
\end{lemma}

\vspace{.2cm}

\begin{prop}\label{u(S)-u(T)}For $T\in\CM$ define $u(T)\in\CU(M_2(\CM))$ by
\begin{equation}\label{eq4-23}
  u(T)=\begin{pmatrix}
    (\unit + T\cc T)^{-\frac 12}T\cc & -(\unit + T\cc T)^{-\frac 12}\\
    (\unit + T T\cc)^{-\frac 12} &  (\unit + T T\cc)^{-\frac 12}T
  \end{pmatrix}.
\end{equation}
Then for all $S,T\in\CM$,
\begin{equation}
  \|u(S)-u(T)\|\leq 2\|S-T\|.
\end{equation}
\end{prop}

\vspace{.2cm}

According to Proposition~\ref{circular}, if $u=
(u_{ij})_{i,j=1}^2\in \CU(M_2(\CM))$, then $u_{11}x + u_{12}y$ and
$u_{21}x+u_{22}y$ are $\ast$-free circular elements in $\CN$, and they are
$\ast$-free from $\CM$. In particular, we may define $g_u(z)\in L^p(\CN,
\tau)$, $(0<p<1)$, by
\begin{equation}
  g_u(z)= ( u_{11}x+ u_{12}y)( u_{21}x+ u_{22}y)^{-1}= ( u_{11}z+ u_{12})( u_{21}z+ u_{22})^{-1}.
\end{equation}

\begin{prop}\label{g_U} 
\begin{itemize}
  \item[(i)] Let $p\in(0,1)$. Then for every $u\in \CU(M_2(\CM))$,
    \begin{equation}
      \|g_u(z)\|_p = \|z\|_p<\infty.
    \end{equation}
  \item[(ii)] Let $p\in(0,\frac23)$. Then there is a constant $C_p^{(2)}>0$ such that for all
$u, v\in \CU(M_2(\CM))$,
\begin{equation}\label{eq4-13}
\|g_u(z)-g_v(z)\|_p \leq C_p^{(2)} \|u-v\|.
\end{equation}
\end{itemize}
\end{prop}

\vspace{.2cm}

\vspace{.2cm}

{\it Proof of Theorem~\ref{Lipschitz}.} Let $T\in\CM$, and let $u(T)\in
\CU(M_2(\CM))$ be given by \eqref{eq4-23}. Moreover define (a circular
system) $\{v,w\}\subseteq \CN$ by
\begin{equation}\label{eq4-24}
  \begin{pmatrix} v \\ w  \end{pmatrix} =
  u(T)\begin{pmatrix} x \\ y  \end{pmatrix} =
  \begin{pmatrix}
    (\unit + T\cc T)^{-\frac 12}(T\cc x - y) \\
    (\unit + T T\cc)^{-\frac 12}(x+Ty)
  \end{pmatrix}.
\end{equation}

Then
\begin{eqnarray*}
  T+z & = & (Ty+x)y^{-1} \\
  & = & (\unit + TT\cc)^{\frac12} w y^{-1},
\end{eqnarray*}
so for each $p\in(0,1)$, $T+z$ has an inverse $(T+z)^{-1}\in L^p(\CN,
\tau)$ given by
\begin{equation}\label{eq4-30}
   (T+z)^{-1} = y  w^{-1}  (\unit + TT\cc)^{-\frac12}.
\end{equation}

Moreover,
\begin{eqnarray*}
  \begin{pmatrix} x \\ y  \end{pmatrix} &=&
  u(T)\cc\begin{pmatrix} v \\ w  \end{pmatrix} \\
  &=&
  \begin{pmatrix}
   T (\unit + T\cc T)^{-\frac 12}& (\unit + T T\cc)^{-\frac 12}\\
   -(\unit + T\cc T)^{-\frac 12} & T\cc  (\unit + T T\cc)^{-\frac 12}
  \end{pmatrix}
  \begin{pmatrix} v \\ w  \end{pmatrix}\\
  & = &
   \begin{pmatrix}
   (\unit + T T\cc)^{-\frac 12}T & (\unit + T T\cc)^{-\frac 12}\\
   -(\unit + T\cc T)^{-\frac 12} & (\unit + T\cc T)^{-\frac 12}T\cc
  \end{pmatrix}
  \begin{pmatrix} v \\ w  \end{pmatrix}.
\end{eqnarray*}

In particular,
\begin{equation}
y= (\unit + T\cc T)^{-\frac 12}(T\cc w -v),
\end{equation}
and hence by \eqref{eq4-30}, 
\begin{eqnarray*}
(T + z)^{-1} &=&  y w^{-1} (\unit + TT\cc)^{-\frac 12}\\
& = & (\unit + T\cc T)^{-\frac 12}(T\cc w -v) w^{-1} (\unit +
TT\cc)^{-\frac 12}\\
& = & T\cc  (\unit +TT\cc)^{-1} - (\unit +
T\cc T)^{-\frac 12} v w^{-1} (\unit +TT\cc)^{-\frac 12}\\
& = & T\cc  (\unit +TT\cc)^{-1} - (\unit +
T\cc T)^{-\frac 12} g_{u(T)}(z) (\unit +TT\cc)^{-\frac 12} .
\end{eqnarray*}

It follows that for $S,T\in \CM$ and $0<p<\frac 23$,
\begin{equation}\label{eq4-22}
  \|(S+z)^{-1}-(T+z)^{-1}\|_p^p \leq \|A\|_p^p + \|B\|_p^p,
\end{equation}
where
\begin{equation}
A = S\cc(\unit+ SS\cc)^{-1} - T\cc(\unit+ TT\cc)^{-1},
\end{equation}
and
\begin{equation}
B = (\unit + S\cc S)^{-\frac 12} g_{u(S)}(z) (\unit +SS\cc)^{-\frac 12} - (\unit +
T\cc T)^{-\frac 12} g_{u(T)}(z) (\unit +TT\cc)^{-\frac 12}.
\end{equation}

According to \eqref{eq4-20},
\[
\|A\|_p\leq \|A\| = \| (\unit +SS\cc)^{-\frac12}S - (\unit+
TT\cc)^{-\frac12}T\| \leq \textstyle{\frac 54}\, \|S\cc - T\cc\| =  \textstyle{\frac 54}\, \|S - T\|.
\]

To get an estimate of $\|B\|_p$, we write $B$ as a sum, $B=B_1+B_2 + B_3$, where $B_1, B_2$ and $B_3$ are given by 
\begin{eqnarray*}
B_1 &=& [(\unit + S\cc S)^{-\frac 12}-(\unit +
T\cc T)^{-\frac 12}] g_{u(S)}(z) (\unit +SS\cc)^{-\frac 12},\\
B_2 &=& (\unit + T\cc T)^{-\frac 12}( g_{u(S)}(z)-  g_{u(T)}(z)) (\unit
+SS\cc)^{-\frac 12},\\
B_3 &=& (\unit + T\cc T)^{-\frac 12} g_{u(T)}(z)[(\unit +SS\cc)^{-\frac
  12}-(\unit +TT\cc)^{-\frac 12} ].
\end{eqnarray*}
 According to \eqref{eq4-21},
\[
\|(\unit + S\cc S)^{-\frac 12}-(\unit +
T\cc T)^{-\frac 12}\|\leq \textstyle{\frac 2\pi}\,\|S-T\|\leq \|S-T\|,
\]
and
\[
\|(\unit + SS\cc)^{-\frac 12}-(\unit +
TT\cc)^{-\frac 12}\|\leq \textstyle{\frac 2\pi}\,\|S\cc-T\cc\|\leq \|S-T\|.
\]

Moreover, by Proposition~\ref{g_U} and Proposition~\ref{u(S)-u(T)},
\[
\| g_{u(S)}(z)-  g_{u(T)}(z)\|_p \leq C_p^{(2)}\|u(S)-u(T)\| \leq 2
C_p^{(2)}\|S-T\|,
\]
and
\[
\| g_{u(S)}(z)\|_p = \|z\|_p.
\]
It follows that
\begin{equation}
\|B_1\|_p \leq \|z\|_p\|S-T\|,
\end{equation}
\begin{equation}
\|B_2\|_p \leq 2 C_p^{(2)}\|S-T\|
\end{equation}
and
\begin{equation}
\|B_3\|_p \leq  \|z\|_p\|S-T\|.
\end{equation}
Altogether we have shown that
\begin{eqnarray*}
\|A\|_p^p+\|B\|_p^p &\leq& \|A\|_p^p + \|B_1\|_p^p +  \|B_2\|_p^p +  \|B_3\|_p^p\\
& \leq & \Big( (\textstyle{\frac 54} )^p + 2  \|z\|_p^p + (2 C_p^{(2)})^p\Big)\|S-T\|^p,
\end{eqnarray*}
so by \eqref{eq4-22}, \eqref{eqLip} holds with $C_p^{(1)} = \Big( (\textstyle{\frac 54} )^p + 2  \|z\|_p^p + 2^p
{C_p^{(2)}}^p\Big)^\frac1p.$ $\endproof$

\vspace{.2cm}

{\it Proof of Lemma~\ref{ABC}.} For all $x>0$ we have:
\[
x^{-\frac 12} = \frac 1\pi \int_0^\infty \frac{a^{-\frac 12}}{x+a}\,\d a.
\]
Therefore every $R\in\CM$ satisfies
\[
(\unit+R\cc R)^{-\frac 12} = \frac 1\pi \int_0^\infty (R\cc R
+(a+1)\unit)^{-1}a^{-\frac 12} \,\d a,
\]
(as a Banach space valued integral).

For $a\geq 1$ define
\begin{eqnarray*}
   A(a)&=& (a\unit +S\cc S)^{-1}S\cc - (a\unit +T\cc T)^{-1}T\cc,\\
   B(a) & = &  (a\unit +S\cc S)^{-1} - (a\unit+ T\cc T)^{-1},\\
   C(a) & = &  (a\unit +SS\cc)^{-1} - (a\unit +TT\cc)^{-1}.
  \end{eqnarray*}

Then
\begin{eqnarray*}
  A(a)& = & (a\unit +S\cc S)^{-1}S\cc - T\cc(a\unit +TT\cc)^{-1}\\
  & = & (a\unit +S\cc S)^{-1}[S\cc (a\unit +TT\cc) -(a\unit +S\cc S)T\cc
  ](a\unit +TT\cc)^{-1}\\
  & = & a (a\unit +S\cc S)^{-1}(S\cc - T\cc)(a\unit + TT\cc)^{-1} + \\
  & & \qquad \qquad \qquad 
   (a\unit +S\cc S)^{-1}S\cc (T-S)T\cc (a\unit +TT\cc)^{-1}.
\end{eqnarray*}

Now, for every $R\in\CM$,
\[
\|(a\unit +R\cc R)^{-1}\|\leq \textstyle{\frac 1a},
\]
and
\begin{eqnarray}
\|(a\unit +R\cc R)^{-1}R\cc\| & = & \|(a\unit +R\cc R)^{-1}R\cc R(a\unit
+R\cc R)^{-1}\|^{\frac12}\nonumber\\
&\leq & \Big(\sup_{t\geq 0}\frac{t}{(a+t)^2}\Big)^{\frac 12}\label{Lip1}\\
& = & \frac{1}{2\sqrt a}.
\end{eqnarray}
Hence
\begin{equation}
\|A(a)\| \leq  \textstyle{\frac 1a} \,\|S\cc - T\cc\| +
\textstyle{\frac{1}{4a}}\, \|T-S\| = \textstyle{\frac{5a}{4}}\,\|S-T\|.
\end{equation}
It follows that
\begin{eqnarray}
 \|(\unit +S\cc S)^{-\frac 12}S\cc - (\unit +T\cc T)^{-\frac 12}T\cc \|
 & \leq & \frac 1\pi \int_0^\infty \|A(a+1)\|a^{-\frac 12}\,\d a\nonumber \\
& \leq & \textstyle{\frac{5}{4\pi}}\, \|S-T\|\int_0^\infty (a+1)a^{-\frac 12}\,\d a
 \nonumber \\
& = & \textstyle{\frac{5}{4}}\, \|S-T\|, \label{eq4-2}
\end{eqnarray}
which is \eqref{eq4-20}.

Also, since
\begin{eqnarray*}
  B(a)& = & (a\unit +S\cc S)^{-1}(T\cc T - S\cc S)(a\unit +T\cc T)^{-1}\\
  & = & (a\unit +S\cc S)^{-1}(T\cc -S\cc)T(a\unit +T\cc T)^{-1}\\
  & & \qquad \qquad \qquad + (a\unit +S\cc S)^{-1}S\cc(T -S)(a\unit +T\cc T)^{-1},
\end{eqnarray*}
where
\begin{eqnarray*}
\|(a\unit +S\cc S)^{-1}(T\cc -S\cc)T(a\unit +T\cc T)^{-1}\|&\leq
&\textstyle{\frac 1a} \|T-S\|\|T(a\unit +T\cc T)^{-1}\|\\
& = & \textstyle{\frac 1a}\|T-S\|\|(a\unit +T\cc T)^{-1}T\cc \|\\
&\leq & \textstyle{\frac 1a}\textstyle{\frac{1}{2\sqrt a}} \|T-S\|,
\end{eqnarray*}
and
\begin{eqnarray*}
\|(a\unit +S\cc S)^{-1}S\cc(T -S)(a\unit +T\cc T)^{-1}\| &\leq &
\textstyle{\frac 1a}\|(a\unit +S\cc S)^{-1}S\cc\|\|T-S\|\\
& \leq & \textstyle{\frac 1a}\textstyle{\frac{1}{2\sqrt a}} \|T-S\|,
\end{eqnarray*}
we find that
\begin{equation}
  \|B(a)\| \leq  \textstyle{\frac{1}{a\sqrt a}}\, \|S-T\|,
\end{equation}
whence
\begin{eqnarray}
   \|(\unit +S\cc S)^{-\frac 12} - (\unit+ T\cc T)^{-\frac 12}\|     & \leq
  & \frac1\pi \int_0^\infty \|B(a+1)\|a^{-\frac 12}\,\d a\nonumber \\
  & \leq &  \|S-T\|\, \frac1\pi \int_0^\infty (a+1)^{-\frac32}a^{-\frac 12}\,\d a
 \nonumber \\
  & = & \textstyle{\frac 2\pi}\, \|S-T\|,
\end{eqnarray}
and this is \eqref{eq4-21}. $\endproof$

\vspace{.2cm}

{\it Proof of Proposition~\ref{u(S)-u(T)}.} With
\begin{eqnarray*}
    A&=&(\unit +S\cc S)^{-\frac 12}S\cc - (\unit +T\cc T)^{-\frac
    12}T\cc,\\
    B&=& (\unit +S\cc S)^{-\frac 12} - (\unit+ T\cc T)^{-\frac 12},\\
    C & = & (\unit +SS\cc )^{-\frac 12} - (\unit+ TT\cc)^{-\frac 12}
  \end{eqnarray*}
one has that
\[
u(S)-u(T)= \begin{pmatrix} A & B\\ C & -A\cc \end{pmatrix},
\]
and hence by Lemma~\ref{ABC},
\begin{eqnarray*}
\|u(S)-u(T)\| &\leq &\Big\| \begin{pmatrix} A & 0\\ 0 & -A\cc \end{pmatrix}
\Big\| + \Big\| \begin{pmatrix} 0 & B\\ C & 0 \end{pmatrix}
\Big\|\\
& = & \|A\| + \max\{\|B\|, \|C\|\}\\
& \leq & \textstyle{\frac54}\, \|S-T\| + \textstyle{\frac 2\pi}\,\|S-T\|\\
& \leq & 2 \|S-T\|. \qquad\qquad\qquad \endproof
\end{eqnarray*}

\vspace{.2cm}

Next we want to prove Proposition~\ref{g_U}. To this end we need a series
of lemmas which we state and prove in the following. 

\vspace{.2cm}

\begin{lemma}\label{w_1w_2l}
  Let $u\in \CU(M_2(\CM))$ with $\|u-\unit\|<1$. Then there exist $w_1,
  w_2\in \CU(\CM)$ and $\ell\in\CM$ such that
  \begin{equation}\label{eq4-7}
    u=wv,
  \end{equation}
  where
  \begin{equation}\label{eq4-8}
    w=\begin{pmatrix} w_1 & 0 \\ 0 & w_2 \end{pmatrix},
  \end{equation}
  and
  \begin{equation}\label{eq4-9}
    v = \begin{pmatrix} (\unit-\ell\ell\cc)^\frac12 & l \\ -\ell\cc &  (\unit-\ell\cc\ell)^\frac12 \end{pmatrix}.
  \end{equation}
  $w_1$, $w_2$ and $\ell$ are uniquely determined by
  \eqref{eq4-7},  \eqref{eq4-8} and  \eqref{eq4-9}. Moreover, $\|\ell\|<1$,
  $\|v-\unit\|\leq 2\|u-\unit\|$ and $\|w_i-\unit\|\leq 3 \|u-\unit\|$,
  $i=1,2$.
\end{lemma}

\proof $u=(u_{ij})_{i,j=1}^2$ with $u_{ij}\in\CM$, where
$\|u_{ii}-\unit\|<1$, $i=1,2$. Thus, $u_{11}$ and $u_{22}$ are invertible
and have unique polar decompositions $u_{ii}=w_i|u_{ii}|$, where
$w_i\in\CU(\CM)$, $i=1,2$. Put $\ell = w_1\cc u_{12}$, $m=w_2\cc u_{21}$ and
\[
v= \begin{pmatrix} w_1\cc & 0 \\ 0 & w_2\cc \end{pmatrix}u =
\begin{pmatrix} |u_{11}| & \ell \\ m & |u_{22}| \end{pmatrix}.
\]
Since $v$ is unitary, $ |u_{11}|^2 + \ell\ell\cc = \unit$, $\ell\ell\cc
+  |u_{22}|^2 = \unit$ and $m |u_{11}| +  |u_{22}|\ell\cc = 0$. Hence
\begin{equation}
 |u_{11}| = (\unit - \ell\ell\cc)^\frac12,
\end{equation}
\begin{equation}
|u_{22}| = (\unit -\ell\cc\ell)^\frac12
\end{equation}
and
\begin{equation}
m = -|u_{22}|\ell\cc|u_{11}|^{-1} = -\ell\cc.
\end{equation}
Thus, $v$ is given by \eqref{eq4-9}. Moreover, due to the uniqueness of the
polar decompositions of $u_{11}$ and $u_{22}$, $w_1$, $w_2$ and $\ell$ are uniquely determined
by \eqref{eq4-7},  \eqref{eq4-8} and  \eqref{eq4-9}.

Since $\|u-\unit\|<1$, we find that
\begin{equation}
\|\ell\| = \|u_{12}\| \leq \|u-\unit\|<1.
\end{equation}
Also, when $t\in [0,1]$,  then $0\leq 1-\sqrt t \leq 1-t$, and hence
\[
\|\unit - (\unit-\ell\ell\cc)^\frac 12\|\leq \|\unit -
(\unit-\ell\ell\cc)\| = \|\ell\|^2,
\]
and
\[
\|\unit - (\unit-\ell\cc\ell)^\frac 12\|\leq \|\unit -
(\unit-\ell\cc\ell)\| = \|\ell\|^2.
\]
Since
\[
\unit-v =  \begin{pmatrix} \unit - (\unit-\ell\ell\cc)^\frac12 & 0 \\
  0 &  \unit - (\unit-\ell\cc\ell)^\frac12 \end{pmatrix} +
 \begin{pmatrix} 0 & -l \\ \ell\cc &  0 \end{pmatrix},
\]
we may use the estimates obtained above to see that
\begin{eqnarray}
  \|v-\unit\| & \leq & \max\{\| \unit - (\unit-\ell\ell\cc)^\frac12\|, \,
  \| \unit - (\unit-\ell\cc\ell)^\frac12\|\} + \|\ell\| \nonumber\\
  &\leq &\|\ell\|^2 + \|\ell\|\nonumber\\
  & \leq & 2\|\ell\|\nonumber \\
  & \leq & 2\|u-\unit\|.
\end{eqnarray}

Finally, $w$ given by \eqref{eq4-8} satisfies
\begin{eqnarray*}
  \|w-\unit\| & = & \|uv\cc -\unit\|\\
  & = & \|(u-v)v\cc\|\\
  & = & \|u-v\|\\
  & \leq & \|u-\unit\| + \|v-\unit\|\\
  & \leq & 3\|u-\unit\|,
\end{eqnarray*}
and hence $\|w_i-\unit\|\leq  3\|u-\unit\|$, $i=1,2$. $\endproof$

\vspace{.2cm}

\begin{lemma}\label{w_3,m}
  Let $\ell\in\CM$ with $\|\ell\|<1$. Then there exists $w_3\in \CU(\CM)$  such that
  \begin{equation}\label{eq4-10}
     \begin{pmatrix} (\unit-\ell\ell\cc)^\frac12 & l \\ -\ell\cc &
  (\unit-\ell\cc\ell)^\frac12 \end{pmatrix}=
     \begin{pmatrix} w_3 & 0 \\ 0 & \unit \end{pmatrix}
     \begin{pmatrix} (\unit-|\ell|^2)^\frac12 & |\ell| \\ -|\ell| &
  (\unit-|\ell|^2)^\frac12 \end{pmatrix}
     \begin{pmatrix} w_3\cc & 0 \\ 0 & \unit \end{pmatrix}.
  \end{equation}
\end{lemma}

\proof Since $\CM$ is finite, we may take a unitary $w_3\in \CU(\CM)$ such
that $\ell = w_3|\ell|$. Let  $m_1 = w_3\cc  (\unit-\ell\ell\cc)^\frac12 w_3$
and  $m_2 =  (\unit-\ell\cc\ell)^\frac12$. Then 
\[
\begin{pmatrix} w_3\cc & 0 \\ 0 & \unit \end{pmatrix}
\begin{pmatrix} (\unit-\ell\ell\cc)^\frac12 & l \\ -\ell\cc &
  (\unit-\ell\cc\ell)^\frac12 \end{pmatrix} 
\begin{pmatrix} w_3 & 0 \\ 0 & \unit \end{pmatrix} =
 \begin{pmatrix} m_1 & |\ell| \\ -|\ell| & m_2 \end{pmatrix}.
\]
In particular, $ \begin{pmatrix} m_1 & |\ell| \\ -|\ell| & m_2 \end{pmatrix}$ is
unitary, and hence $m_1^2 + |\ell|^2 = \unit$, so $m_1 = (\unit-
|\ell|^2)^\frac12$. This proves \eqref{eq4-10}. $\endproof$

\vspace{.2cm}

\begin{lemma} \label{w_4,w_5,theta} Let $m\in \CM_+$ with $\|m\|<1$. Then
  there exist unitaries $w_4, w_5\in \CM$ and $\theta\in\R$ such that
  \begin{eqnarray}
    v'&:=&\begin{pmatrix} (\unit -m^2)^\frac12 & m \\ -m &
  (\unit-m^2)^\frac12 \end{pmatrix}\nonumber \\
    & = & 
     \begin{pmatrix} w_4\cc & 0 \\ 0 & \unit \end{pmatrix}
     \begin{pmatrix} \cos\theta \unit & \sin\theta w_5 \\ -\sin\theta
       w_5\cc  & \cos\theta \unit \end{pmatrix}
     \begin{pmatrix} \cos\theta \unit & \sin\theta w_5\cc \\ -\sin\theta
       w_5  & \cos\theta \unit \end{pmatrix}
     \begin{pmatrix} w_4 & 0 \\ 0 & \unit \end{pmatrix},\label{eq4-11}
   \end{eqnarray}
   $0\leq\theta \leq \|v'-\unit\|$ and $\|w_4-\unit\|\leq 2\|v'-\unit\|$.
\end{lemma}

\proof If $m=0$, then $v' =\unit$, and \eqref{eq4-11} holds with
$\theta=0$, $w_4=w_5$. Now, assume that $m\neq 0$. Then we may choose
$\theta\in (0,\frac\pi 4)$ such that $\sin(2\theta)=\|m\|$. Let
\begin{equation}
m'=\frac{1}{\sin(2\theta)}m = \frac{1}{2\cos\theta\sin\theta}m.
\end{equation}
Then $m'=(m')\cc$, $\|m'\|=1$, and with
\begin{equation}
w_5 = m' + \i\sqrt{\unit-(m')^2}\in \CU(\CM)
\end{equation}
one has that
\[
m' = \textstyle{\frac12} (w_5+w_5\cc).
\]
Let
\[
v''= \begin{pmatrix} \cos\theta \unit & \sin\theta w_5 \\ -\sin\theta
       w_5\cc  & \cos\theta \unit \end{pmatrix}
     \begin{pmatrix} \cos\theta \unit & \sin\theta w_5\cc \\ -\sin\theta
       w_5  & \cos\theta \unit \end{pmatrix}
     = \begin{pmatrix} k & m \\ 
       -m & k\cc \end{pmatrix},
\]
where
\begin{equation}\label{eq4-31}
k = \cos^2\theta\unit -\sin^2\theta w_5^2.
\end{equation}
Since $v''$ is unitary, $kk\cc =k\cc k = \unit-m^2$. Then let
\begin{equation}
k = w_4 (\unit-m^2)^\frac12
\end{equation}
be the polar decomposition of $k$.
$\sin^2\theta < \cos^2\theta$, when $0<\theta<\frac\pi 4$, so $k$ is
invertible, and $w_4$ must be unitary. Moreover, $w_4$ belongs to the abelian von Neumann algebra
$W\cc(k)$, and therefore it commutes with $m=(\unit-k\cc k)^\frac12$. It follows that
\begin{equation}
\begin{pmatrix} w_4\cc & 0 \\ 0 & \unit \end{pmatrix}
     \begin{pmatrix} k  & m \\ -m & k\cc \end{pmatrix}
     \begin{pmatrix} w_4 & 0 \\ 0 & \unit \end{pmatrix}=
     \begin{pmatrix}  (\unit-m^2)^\frac12 & m \\ -m & (\unit-m^2)^\frac12
     \end{pmatrix},
\end{equation}
that is, \eqref{eq4-11} holds. Also, since $2\theta\in (0,\frac\pi 2)$,
$\|m\|=\sin(2\theta)\geq \textstyle{\frac 2\pi} \, 2\theta$, and we get that
\begin{equation}
\theta\leq \textstyle{\frac\pi 4}\, \|m\|\leq \|m\|\leq \|v'-\unit\|.
\end{equation}

$w_5^2$ is unitary, and hence by \eqref{eq4-31} and the spectral mapping theorem, the spectrum of $k$ must be contained in the boundary of the
closed ball $\overline{B(\cos^2\theta, \sin^2\theta)}\subseteq \C$. As
$\cos^2\theta > \sin^2\theta$, this ball is contained in $\{z\in\C
\,|\, \re z >0\}$, so there is a continuous determination of $\arg z$ in
$\sigma(k)$ taking values in $(-\frac\pi 2, \frac\pi 2)$. In
particular,
\begin{equation}
w_4 = \e^{\i \arg(k)}
\end{equation}
and we get that
\begin{eqnarray}
  \|w_4-\unit\| &\leq &\|\arg(k)\|\nonumber\\
  & \leq & \max\{|\arg (z)|\,|\, z\in\partial B(\cos^2\theta,
  \sin^2\theta)\}\nonumber\\
  & = & \arcsin(\tan^2\theta).
\end{eqnarray}
Making use of the estimates
\[
0\leq \tan t \leq \frac 4\pi t, \qquad (0\leq t\leq \frac\pi 4),
\]
and
\[
0\leq \arcsin t\leq \frac\pi 2 t, \qquad (0\leq t\leq 1),
\]
we finally find that
\begin{eqnarray*}
\|w_4-\unit\|&\leq  &\arcsin(\tan^2\theta)\\
&\leq & \arcsin(\textstyle{\frac{16}{\pi^2}}\theta^2)\\
&\leq &\textstyle{\frac\pi 2}\textstyle{\frac{16}{\pi^2}}\theta^2 \\
&\leq & 2\theta \\
&\leq & 2\|v'-\unit\|. \qquad \endproof
\end{eqnarray*}

\vspace{.2cm}

\begin{definition}
  For $p\in (0,1)$ define $d_p: \CU(M_2(\CM))\times \CU(M_2(\CM))\rightarrow [0,\infty[$ by
  \begin{equation}
    d_p(u, v)= \|g_u(z)-g_v(z)\|_p^p, \qquad (u,v\in  \CU(M_2(\CM)).
  \end{equation}
\end{definition}

\vspace{.2cm}

Clearly, for all $u,v,w\in  \CU(M_2(\CM))$, we have:
\begin{equation}\label{eq4-25}
  d_p(u,w) \leq d_p(u,v) + d_p(v,w).
\end{equation}

\vspace{.2cm}

\begin{prop}\label{i,ii,iii} Let $0<p<\frac 23$. Then
  \begin{itemize}
    \item[(i)] $d_p$ is right-invariant on $ \CU(M_2(\CM))$.
    \item[(ii)] $d_p$ is left-invariant w.r.t. multiplication by diagonal
    unitaries.
    \item[(iii)] For all $w_1, w_2\in \CU(\CM)$ one has that $w =
    \begin{pmatrix} w_1 & 0 \\ 0 & w_2 \end{pmatrix}$ satisfies
    \begin{equation}
      d_p(w,\unit)\leq \big(\|w_1-\unit\|^p +
      \|w_2-\unit\|^p\big)\|z\|_p^p.
    \end{equation}
  \end{itemize}
\end{prop}

\proof (i) Let $u,v,w\in \CU(M_2(\CM))$. Then with
\[
\begin{pmatrix} x' \\ y' \end{pmatrix} = w \begin{pmatrix} x \\ y
\end{pmatrix}
\]
and $z'= x'(y')^{-1}$ we have that
\[
g_{uw}(z)-g_{vw}(z) = g_u(z')-g_v(z') .
\]
According to  Proposition~\ref{circular}, $\{x', y'\}$ is a circular
system which is $\ast$-free from $\CM$, so
\begin{eqnarray*}
d_p(uw,vw) & = & \|g_{uw}(z)-g_{vw}(z)\|_p\\
& = & \| g_u(z')-g_v(z')\|_p \\
& = & \| g_u(z)-g_v(z)\|_p\\
& = & d_p(u,v).
\end{eqnarray*}

(ii) Let $u= (u_{ij})_{i,j=1}^2,v=(v_{ij})_{i,j=1}^2\in \CU(M_2(\CM))$, and let $w= \begin{pmatrix} w_1 & 0 \\ 0 &
  w_2\end{pmatrix}\in \CU(M_2(\CM))$. Then $w_1, w_2\in\CU(\CM)$, and hence
\begin{eqnarray*}
  d_p(wu, wv) & = &  \|g_{wu}(z)-g_{wv}(z)\|_p^p \\
  & = & \|w_1 ( u_{11}x+ u_{12}y)( u_{21}x+
  u_{22}y)^{-1} w_2\cc -w_1 ( v_{11}x+ v_{12}y)( v_{21}x+
  v_{22}y)^{-1} w_2\cc \|_p^p\\
  & = & \|( u_{11}x+ u_{12}y)( u_{21}x+
  u_{22}y)^{-1} -( v_{11}x+ v_{12}y)( v_{21}x+
  v_{22}y)^{-1}\|_p^p\\
  & = &  \|g_{u}(z)-g_{v}(z)\|_p^p\\
  & = & d_p(u,v).
\end{eqnarray*}

(iii) Let $w_1, w_2\in \CU(\CM)$, and let $w=\begin{pmatrix} w_1 & 0 \\ 0 &
         w_2 \end{pmatrix}$. Then
       \begin{eqnarray*}
         d_p(w, \unit) &=& \|w_1x (w_2y)^{-1}- xy^{-1}\|_p^p\\
         & = & \|w_1xy^{-1} - xy^{-1}w_2\|_p^p\\
         & \leq & \|(w_1-\unit)xy^{-1}\|_p^p + \|xy^{-1}(\unit-w_2)\|_p^p\\
         & \leq & (\|(w_1-\unit)\|^p + \|(w_2-\unit)\|^p)\|z\|_p^p. \qquad \endproof
       \end{eqnarray*}

\vspace{.2cm}

\begin{cor}\label{iv,v} Let $n\in\N$, let $u, u_1, \ldots, u_n \in\CU(M_2(\CM))$,
      and let $w\in\CU(M_2(\CM))$ be diagonal. Then
      \begin{itemize}
      \item[(iv)] $d_p(u_1\cdots u_n, \unit)\leq \sum_{i=1}^n
      d_p(u_i,\unit)$,
       \item[(v)] $d_p(wuw\cc, \unit) = d_p(u,\unit)$.
       \end{itemize}
\end{cor}

\proof (iv) follows from \eqref{eq4-25} and from Proposition~\ref{i,ii,iii}~(i). Indeed,
if $u_1, u_2\in \CU(M_2(\CM))$, then
\begin{eqnarray*}
  d_p(u_1 u_2, \unit ) & = & d_p(u_1, u_2\cc)\\
  & \leq & d_p(u_1, \unit) + d_p(\unit, u_2\cc)\\
  & = & d_p(u_1, \unit) + d_p(u_2, \unit).
\end{eqnarray*}
Then (iv) follows for general $n\in\N$ by induction.

(v) is a consequence of Proposition~\ref{i,ii,iii}~(i) and (ii). $\endproof$

\vspace{.2cm}

{\it Proof of Proposition~\ref{g_U}}.(i) Let $u= (u_{ij})_{i,j=1}^2 \in \CU(M_2(\CM))$. Then, as
already mentioned, $u_{11}x+u_{12}y$ and $u_{21}x+u_{22}y$ are $\ast$-free
circular elements. In particular, $g_u(z)=
(u_{11}x+u_{12}y)(u_{21}x+u_{22}y)^{-1}$ has the same $\ast$-distribution
as $z$ (in the sense of \cite[Definition~3.2]{HS}), whence $\|g_u(z)\|_p=\|z\|_p<\infty$ for $p\in(0,1)$.

(ii) Now assume that $p\in(0,\frac23)$. At first we consider
$u\in\CU(M_2(\CM))$ with  $\|u-\unit\|<1$. Take $w= \begin{pmatrix}w_1
  &0\\ 0 & w_2\end{pmatrix}$ and $v$ as in
Lemma~\ref{w_1w_2l}. Then, according to Corollary~\ref{iv,v}~(iv),
\begin{equation}
d_p(u,\unit)\leq d_p(w,\unit)+d_p(v,\unit),
\end{equation}
where
\begin{eqnarray}
d_p(w,\unit) &\leq& \big(\|w_1-\unit\|^p + \|w_2-\unit\|^p\big)\|z\|_p^p\nonumber\\
&\leq & 2\cdot 3^p \|u-\unit\|^p \|z\|_p^p
\end{eqnarray}
(cf. Proposition~\ref{i,ii,iii} and Lemma~\ref{w_1w_2l}).

As in Lemma~\ref{w_3,m}, take
$m=|\ell|\in\CM_+$ an $w_3\in\CU(\CM)$ such that
  \[
     v =
     \begin{pmatrix} w_3 & 0 \\ 0 & \unit \end{pmatrix}
     \begin{pmatrix} (\unit-m^2)^\frac12 & m \\ -m &
  (\unit-m^2)^\frac12 \end{pmatrix}
     \begin{pmatrix} w_3\cc & 0 \\ 0 & \unit \end{pmatrix},
\]
and put
\begin{equation}
v' =  \begin{pmatrix} (\unit-m^2)^\frac12 & m \\ -m &
  (\unit-m^2)^\frac12 \end{pmatrix}.
\end{equation}
Then
\begin{equation}
\|v-\unit\|=\|v'-\unit\| \leq 2\|u-\unit\|,
\end{equation}
and because of Corollary~\ref{iv,v}~(v), $d_p(v',\unit)=d_p(v,\unit)$, so
\begin{equation}
d_p(u,\unit)\leq 2\cdot 3^p \|u-\unit\|^p \|z\|_p^p + d_p(v',\unit).
\end{equation}

As in Lemma~\ref{w_4,w_5,theta} we may take $w_4, w_5\in \CU(\CM)$ and
$\theta\in \R$ such that $0\leq \theta \leq \|v'-\unit\|$,
$\|w_4-\unit\|\leq 2 \|v'-\unit\|$ and
\begin{equation}
v' = \begin{pmatrix} w_4\cc & 0 \\ 0 & \unit \end{pmatrix}
     \begin{pmatrix} \cos\theta \unit & \sin\theta w_5 \\ -\sin\theta
       w_5\cc  & \cos\theta \unit \end{pmatrix}
     \begin{pmatrix} \cos\theta \unit & \sin\theta w_5\cc \\ -\sin\theta
       w_5  & \cos\theta \unit \end{pmatrix}
     \begin{pmatrix} w_4 & 0 \\ 0 & \unit \end{pmatrix}.
\end{equation}
Then, according to Proposition~\ref{i,ii,iii}~(iii) and
Corollary~\ref{iv,v}~(iv), with
\begin{equation}
v'' = \begin{pmatrix} \cos\theta \unit & \sin\theta w_5 \\ -\sin\theta
       w_5\cc  & \cos\theta \unit \end{pmatrix}
     \begin{pmatrix} \cos\theta \unit & \sin\theta w_5\cc \\ -\sin\theta
       w_5  & \cos\theta \unit \end{pmatrix}
\end{equation}
we have that
\begin{eqnarray}
  d_p(v',\unit) & \leq & 2\|w_4-\unit\|_p^p \|z\|_p^p + d_p(v'',\unit)\nonumber\\
  & \leq & 2\cdot 2^p \|v'-\unit\|^p  \|z\|_p^p + d_p(v'',\unit)\nonumber\\
  & \leq & 2\cdot 4^p \|u-\unit\|^p  \|z\|_p^p + d_p(v'',\unit).
\end{eqnarray}

Altogether we have shown that
\begin{equation}
  d_p(u,\unit)  \leq  2\cdot (3^p+4^p) \|u-\unit\|^p  \|z\|_p^p + d_p(v'',\unit).
\end{equation}

Since
\[
\begin{pmatrix} \cos\theta \unit & \sin\theta w_5 \\ -\sin\theta
       w_5\cc  & \cos\theta \unit \end{pmatrix} =  \begin{pmatrix} w_5 & 0 \\ 0  & \unit \end{pmatrix}
     \begin{pmatrix} \cos\theta  & \sin\theta \\ -\sin\theta  & \cos\theta
     \end{pmatrix}
     \begin{pmatrix} w_5\cc & 0 \\ 0  & \unit \end{pmatrix},
\]
and
\[
 \begin{pmatrix} \cos\theta \unit & \sin\theta w_5\cc \\ -\sin\theta
       w_5  & \cos\theta \unit \end{pmatrix} =
  \begin{pmatrix} w_5\cc & 0 \\ 0  & \unit \end{pmatrix}
     \begin{pmatrix} \cos\theta  & \sin\theta \\ -\sin\theta  & \cos\theta
     \end{pmatrix}
     \begin{pmatrix} w_5 & 0 \\ 0  & \unit \end{pmatrix},    
\]
Corollary~\ref{iv,v} implies that with
\begin{equation}
u(\theta)= 
\begin{pmatrix} \cos\theta  & \sin\theta \\ -\sin\theta  & \cos\theta
     \end{pmatrix},
\end{equation}
one has that
\begin{equation}
d_p(v'',\unit)\leq 2 d_p(u(\theta),\unit).
\end{equation}
Thus, when $\|u-\unit\|<1$,
\begin{equation}\label{eq4-12}
  d_p(u,\unit)\leq 4^{p+1}\|u-\unit\|^p \|z\|_p^p + 2d_p(u(\theta),\unit),
\end{equation}
where $0\leq \theta \leq \max\{2\|u-\unit\|, \, \frac\pi 4\}$.

If $\|u-\unit\|\geq 1$, then
\begin{eqnarray*}
  d_p(u,\unit) & = & \|g_u(z)-z\|_p^p\\
  & \leq & \|g_u(z)\|_p^p + \|z\|_p^p\\
  & = & 2 \|z\|_p^p\\
  & \leq & 4^{p+1}\|u-\unit\|^p \|z\|_p^p,
\end{eqnarray*}
so \eqref{eq4-12} still holds in this case.

To get an estimate of $d_p(u(\theta),\unit)$, where  $0\leq \theta \leq
\max\{2\|u-\unit\|, \, \frac\pi 4\}$, let $\phi=\frac\theta 2$. Then
$u(\theta)= u(\phi)^2$, so by the right-invariance of $d_p$ and the inequality
\[
1+\tan^2\phi \leq 2, \qquad 0\leq\phi\leq \frac{\pi}{8},
\]
we get that
\begin{eqnarray*}
  d_p(u(\theta),\unit)&=& d_p(u(\phi), u(\phi)^{-1})\\
  & = & d_p\Big(u(\phi)\begin{pmatrix}0 & 1\\ 1& 0\end{pmatrix},
  u(\phi)^{-1}\begin{pmatrix}0 & 1\\ 1& 0\end{pmatrix}\Big)\\
  & = & \Big\|\frac{\sin\phi \,z+ \cos\phi}{\cos\phi \,z -\sin\phi } -
  \frac{-\sin\phi \,z+ \cos\phi}{\cos\phi \,z +\sin\phi }\Big\|_p^p\\
  & = & \Big\| \frac{2\tan\phi \,(z^2+\unit) }{z^2-\tan^2\phi\,  } \Big\|_p^p\\
  & = &  \Big\| 2\tan\phi \,  + \frac{2\tan\phi
  \, (1+\tan^2\phi)}{z^2-\tan\phi\, }\Big\|_p^p\\
  & \leq & (2\tan\phi)^p + (4\tan\phi)^p \Big\| \frac{1}{z^2-\tan\phi}
  \Big\|_p^p.
\end{eqnarray*}
Hence, by \eqref{6.0} and the fact that 
\[
\tan\phi \leq 2\phi =\theta, \qquad 0\leq\phi\leq \frac{\pi}{8},
\]
\begin{eqnarray*}
d_p(u(\theta),\unit)&\leq & (2\tan\phi)^p + (4\tan\pi)^p\|z^2\|_p^p\\
&\leq & (2\theta)^p + (4\theta)^p\|z^2\|_p^p\\
&\leq & 2\theta^p \big(1+2\|z^2\|_p^p\big)\\
&\leq & 4\cdot 2^p\|u-1\|^p\big(1+2\|z^2\|_p^p\big).
\end{eqnarray*}
Finally, by insertion of the above estimate into
\eqref{eq4-12}, we find that with
\begin{equation}
C_p^{(2)} = \big(4^{p+1}\|z\|_p^p + 2^{p+2} (1+2\|z^2\|_p^p)\big)^\frac1p,
\end{equation}
one has that
\begin{equation}
d_p(u, \unit)^\frac 1p \leq C_p \|u-\unit\|.
\end{equation}
Then for all $u,v\in \CU(M_2(\CM))$,
\begin{eqnarray*}
\|g_u(z)-g_v(z)\|_p &= &d_p(u,v)^\frac 1p \\
& = & d_p(uv^{-1},\unit)^\frac 1p \\
& \leq &  C_p^{(2)} \|uv^{-1}-\unit\|\\
& = &  C_p^{(2)}\|u-v\|,
\end{eqnarray*}
and this is \eqref{eq4-13}.  $\endproof$

%% file: integration.tex
\section{Integration in $\Lp$ for $0<p<1$}

In this section we consider a II$_1$-factor $\CM$ with faithful tracial
state $\tau$ and a fixed $p\in (0,1)$. In this case
\[
\|A\|_p=\big(\tau(|A|^p)\big)^{\frac1p}, \qquad A\in\CM,
\]
does not define a norm on $\CM$ but only a quasi--norm satisfying
\begin{equation}\label{6.1'}
\|A+B\|_p^p \leq \|A\|_p^p+ \|B\|_p^p,\qquad A,B\in\CM
\end{equation}
(cf. \cite{FK}). $L^p(\CM,\tau)$ is complete w.r.t. the metric
\[
d_p(A,B)=\|A-B\|_p^p, \qquad A,B\in\CM.
\]

Although $L^p(\CM,\tau)$ is not a locally convex vector space, one can define a Riemann type integral of a vector valued function
\[
f: [a,b]\rightarrow\Lp,
\]
provided $f$ fulfills the condition that for some $\alpha >
\frac{1-p}{p}$ and for some positive constant $C$, 
\begin{equation}\label{eq3-15}
  \forall \, x,y\in [a,b]: \quad \|f(x)-f(y)\|_p\leq C|x-y|^\alpha.
\end{equation}
This is a special case of an integral introduced by Turpin and Waelbroeck in 1968--71 (cf. \cite{TuWa}, \cite{Wa}) and further developed by Kalton in 1985 (cf. \cite[Section~3]{Ka}). We shall only need the properties of the Turpin--Waelbroeck integral stated in Definition~\ref{define int} and Theorem~\ref{partition1} below. For the convenience of the reader we have enclosed a self--contained proof of Theorem~\ref{partition1} in the appendix (section~10).

Suppose \eqref{eq3-15} is fulfilled. Then for each $n\in\N_0$ define $S_n\in L^p(\CM,\tau)$ by
\begin{equation}
S_n = \frac{b-a}{2^n}\sum_{k=1}^{2^n} f\Big(a+k \frac{b-a}{2^n}\Big), \qquad (n\in\N_0).
\end{equation}
Then for every $n\geq 2$,
\[
\begin{split}
S_n-S_{n-1}&= \frac{b-a}{2^n}\sum_{k=1}^{2^n} (-1)^{k-1}f\Big(a+k \frac{b-a}{2^n}\Big)\\
&= \frac{b-a}{2^n}\sum_{j=1}^{2^{n-1}}T_{n,j},
\end{split}
\]
where
\[
T_{n,j}= f\Big(a+(2j-1) \frac{b-a}{2^n}\Big)-f\Big(a+2j \frac{b-a}{2^n}\Big).
\]
Hence, by \eqref{eq3-15}
\[
\|T_{n,j}\|_p\leq C\Bigg(\frac{b-a}{2^n}\Bigg)^\alpha, \qquad 1\leq j\leq 2^{n-1},
\]
and then by \eqref{6.1'},
\begin{equation}\label{eq3-4}
\|S_n-S_{n-1}\|_p^p \leq \frac{2^{n-1}}{2^{np}}C^p
\Big(\frac{b-a}{2^n}\Big)^{p\alpha + p} = \frac{C^p(b-a)^p}{2\cdot
  2^{n(p+p\alpha +1)}}.
\end{equation}
Since $\alpha >\frac{1-p}{p}$, $p+p\alpha -1>0$, and hence
\[
\sum_{n=2}^\infty \|S_n-S_{n-1}\|_p^p <\infty.
\]
It follows that $(S_n)_{n=1}^\infty$ is a Cauchy sequence in $\Lp$ w.r.t. the
metric $d_p$, and therefore $\lim_{n\rightarrow\infty}S_n$ exists in $L^p(\CM,\tau)$.

\vspace{.2cm}

\begin{definition}\label{define int} For $f:\R\rightarrow \Lp$ as above we define $ \int_a^b f(x)\,\d x \in\Lp$ by
  \begin{equation}
    \int_a^b f(x)\,\d x := \lim_{n\rightarrow\infty} S_n.
  \end{equation}
\end{definition}

\vspace{.2cm}

\begin{thm}\label{partition1} Suppose $f:[a,b]\rightarrow \Lp$ satisfies \eqref{eq3-15} for some $\alpha>\frac1p -1$. Let $I=(a=x_0<x_1<\cdots<x_{n-1}<x_{n}=b)$ be an arbitrary
partition of $[a,b]$, and for arbitrary $t_i\in [x_{i-1},x_i]$ define
\begin{equation}\label{eq3-6}
  M=\sum_{i=1}^n f(t_i)(x_i-x_{i-1}).
\end{equation}

Then with $\delta(I) = \max_{1\leq i\leq n}(x_i-x_{i-1})$,
\begin{equation}\label{partition}
  \Big\|M- \int_a^b f(x)\,\d x \Big\|_p^p \leq
  \frac{C^p(b-a)\delta(I)^{p+\alpha p-1}}{p+p\alpha-1}.
\end{equation}
It follows that if $(I_n)_{n=1}^\infty$ is a
sequence of partitions of $[a,b]$, such that the fineness
$\delta(I_n)$ of the $n$'th partition $I_n$ tends to
zero as $n$ tends to infinity, and if for each $n\in\N$ we associate to
$I_n$ a finite sum $M_n\in\Lp$ as in \eqref{eq3-6}, then $\|M_n-\int_a^b
f(x)\,\d x\|_p\rightarrow 0$ as $n\rightarrow \infty$. 
\end{thm}

\proof See section~10. 

\vspace{.2cm}

By \cite[Section~2]{HS}, the definition of the Fuglede--Kadison determinant $\Delta(T)$ and the Brown measure $\mu_T$ can be extended to all unbounded operators in $\CMD$, where
\[
\CMD=\{T\in\tilde\CM\,|\,\int_0^\infty \log^+ t\,\d\mu_{|T|}(t)<\infty\}.
\]
Note that $L^p(\CM,\tau)\subseteq \CMD$ for all $p\in (0,\infty)$.

\begin{definition}
  For $T\in\CMD$ define $r'(T)$, the {\it modified spectral radius} of $T$, by
  \begin{equation}
    r'(T)=\sup\{|z|\,|\,z\in \supp(\mu_T)\}.
  \end{equation}
\end{definition}

\vspace{.2cm}

\begin{prop}\label{modspecrad}
  For each $T\in\CMD$ and each $p\in (0,\infty)$ we have that
  \begin{equation}
    r'(T)\leq \limsup_{m\rightarrow\infty}\|T^m\|_p^{\frac 1m}.
  \end{equation}
\end{prop}

\proof $r'(T)$ is the essential supremum of $|\lambda|$ w.r.t. $\mu_T$. Hence by
application of \cite[Proposition~2.15]{HS} and
\cite[Theorem~2.19]{HS}, 
\begin{eqnarray*}
  r'(T) & = & \lim_{q\rightarrow\infty}\Big(\int_\C |z|^q
  \,\d\mu_{T}(z)\Big)^{\frac 1q}\\
  & = &  \lim_{m\rightarrow\infty}\Big(\int_\C |z|^{pm}
  \,\d\mu_{T}(z)\Big)^{\frac{1}{pm}}\\
  & = & \lim_{m\rightarrow\infty}\Big(\int_\C |z^m|^p
  \,\d\mu_{T}(z)\Big)^{\frac{1}{pm}}\\
  & = & \lim_{m\rightarrow\infty}\Big(\int_\C |z|^p
  \,\d\mu_{T^m}(z)\Big)^{\frac{1}{pm}}\\
  &\leq &  \limsup_{m\rightarrow\infty}\Big(\|T^m\|_p^p)^{\frac{1}{pm}}\\
  & = & \limsup_{m\rightarrow\infty}\|T^m\|_p^{\frac{1}{m}}. \qquad \endproof
\end{eqnarray*}

\vspace{.2cm}

\begin{prop}\label{prop6.9} Let $T\in\CMD$ and let $P\in\CM$ be a non--trivial projection such that $PTP=TP$. Then
\begin{equation}\label{prop6.9a}
\Delta(T)=\Delta_{P\CM P}(PTP)^{\tau(P)}\,\Delta_{P^\bot\CM P^\bot}(P^\bot TP^\bot)^{1-\tau(P)},
\end{equation}
and
\begin{equation}\label{prop6.9b}
\mu_T=\tau(P)\mu_{PTP}+ (1-\tau(P))\mu_{P^\bot T P^\bot},
\end{equation}
where the Brown measures $\mu_{PTP}$ and $\mu_{P^\bot T P^\bot}$ are computed relative to the II$_1$--factors $P\CM P$ and 
$P^\bot\CM P^\bot$, respectively.
\end{prop}

\proof \eqref{prop6.9a} was proven in \cite[Proposition~2.24]{HS}. Since $PTP=TP$, we have that for all $\lambda\in\C$, $P(T-\lambda\unit)P=(T-\lambda\unit)P$, and thus by \eqref{prop6.9a}, 
\begin{equation}\label{prop6.9c}
\log\Delta(T-\lambda\unit)=\tau(P)\log \Delta_{P\CM P}(P(T-\lambda\unit)P)+(1-\tau(P))\log\Delta_{P^\bot\CM P^\bot}(P^\bot (T-\lambda\unit)P^\bot), \qquad \lambda\in\C.
\end{equation}
\eqref{prop6.9b} now follows by taking the Laplacian (in the Schwartz distribution sense) on both sides of \eqref{prop6.9c} (cf. \cite[Definition~2.13]{HS}). $\endproof$

\vspace{.2cm}

\begin{thm}\label{projection} Suppose $T\in\CMD$ and that $T$ has empty
  point spectrum.\footnote{The {\it point spectrum} of $T$ is the set of
  eigenvalues of $T$.} Moreover, assume that for some $p\in(\frac12, 1)$ there
  exist $\alpha\in(\frac1p -1, 1]$ and a positive constant $C$ such that
  \begin{itemize}
    \item[(i)] $(T-\lambda\unit)^{-1}\in\Lp$ for all $\lambda\in\C$,
    \item[(ii)] $\|(T-\lambda\unit)^{-1}-(T-\mu\unit)^{-1}\|_p \leq
    C|\lambda-\mu|^\alpha$ for all $\lambda, \mu\in\C$.
  \end{itemize}

  Then there is a $T$-invariant subspace $\CK$ affiliated with
  $\CM$ such that when $P\in\CM$ denotes the projection onto  $\CK$, and
  when we write
  \[
  T=\begin{pmatrix}
    T_{11} & T_{12}\\
    0 & T_{22}
    \end{pmatrix},
  \]
  according to the decomposition $\CH=\CK \oplus \CK^\bot$, then
  $T_{11}\in L^p(P\CM P)$, $T_{22}^{-1}\in L^p(P^\bot \CM P^\bot)$,
  $\supp(\mu_{T_{11}})\subseteq \overline{B(0,1)}$ and
  $\supp(\mu_{T_{22}})\subseteq \{z\in\C\,|\,|z|\geq 1\}$.
\end{thm}

\proof Define $f:\R\rightarrow \Lp$ by
\[
f(t)=(\e^{\i t}\unit - T)^{-1}\e^{\i t}, \qquad (t\in\R).
\]

Then $f$ is H\"older continuous with exponent $\alpha$. Indeed,
\begin{eqnarray}\label{eq3-12}
  \|f(t)-f(s)\|_p^p & = & \|(\e^{\i t}\unit - T)^{-1}\e^{\i t}- (\e^{\i
  s}\unit - T)^{-1}\e^{\i s}\|_p^p \nonumber \\
  &\leq & \|((\e^{\i t}\unit - T)^{-1}-(\e^{\i s}\unit - T)^{-1})\e^{\i
  t}\|_p^p+  \|(\e^{\i s}\unit - T)^{-1}(\e^{\i s}-\e^{\i t})\|_p^p \nonumber\\
& \leq & C^p|\e^{\i t}-\e^{\i s}|^{\alpha p} + (C')^p|s-t|^p\nonumber  \\
& \leq & C^p|s-t|^{\alpha p} + (C')^p|s-t|^p,
\end{eqnarray}

where
\[
C'=\max\{\|(\lambda\unit-T)^{-1}\|_p\,|\, |\lambda|=1\}<\infty.
\]

It follows that $f$ is H\"older continuous with exponent
$\min\{\alpha,1\}=\alpha$ and that we may define $E\in\Lp$ by
\begin{equation}
  E=\frac{1}{2\pi\i}\int_{\partial B(0,1)}(\lambda\unit-T)^{-1}\,\d\lambda
  :=\frac{1}{2\pi}\int_0^{2\pi}f(t)\,\d t.
\end{equation}

We are going to prove that the range projection $P$ of $E$ has the
properties stated in Theorem~\ref{projection}.

At first we prove that $E^2=E$. To see this, note that with
\begin{eqnarray*}
s_k^{(n)} &=& \frac{2\pi k}{n}, \qquad (1\leq k\leq n),\\
t_k^{(n)} &=& \frac{2\pi (k-\frac12)}{n}, \qquad (1\leq k\leq n)
\end{eqnarray*}

\begin{equation}\label{eq3-6a}
  E_n = \frac 1n \sum_{k=1}^n(\e^{\i s_k^{(n)}}\unit -T)^{-1}\e^{\i
  s_k^{(n)}}
\end{equation}
and
\begin{equation}
  F_n = \frac 1n \sum_{k=1}^n(\e^{\i t_k^{(n)}}\unit -T)^{-1}\e^{\i
  t_k^{(n)}}
\end{equation}

one has that
\[
\lim_{n\rightarrow\infty}\|E-E_n\|_p=\lim_{n\rightarrow\infty}\|E-F_n\|_p
=0.
\]

Moreover,
\begin{eqnarray*}
  E_nF_n &=& \frac{1}{n^2} \sum_{k,l_1}^n \e^{\i s_k^{(n)}}\e^{\i
  t_l^{(n)}}(\e^{\i s_k^{(n)}}\unit -T)^{-1}(\e^{\i t_l^{(n)}}\unit
  -T)^{-1}\\
  & = & \frac{1}{n^2} \sum_{k,l_1}^n
  \frac{\e^{\i s_k^{(n)}}\e^{\i t_l^{(n)}}}{\e^{\i t_l^{(n)}}-\e^{\i s_k^{(n)}}}\Big((\e^{\i s_k^{(n)}}\unit -T)^{-1}-(\e^{\i t_l^{(n)}}\unit
  -T)^{-1}\Big)\\
  & = & \frac1n \sum_{k=1}^n a_k^{(n)}\e^{\i s_k^{(n)}}(\e^{\i
  s_k^{(n)}}\unit -T)^{-1} +  \frac1n \sum_{l=1}^n b_l^{(n)}\e^{\i
  t_l^{(n)}}(\e^{\i t_l^{(n)}}\unit -T)^{-1},
\end{eqnarray*}

where
\[
 a_k^{(n)} = \frac 1n \sum_{l=1}^n\frac{\e^{\i t_l^{(n)}}}{\e^{\i
 t_l^{(n)}}-\e^{\i s_k^{(n)}}} = \frac 1n \sum_{l=1}^n\frac{1}{1-\e^{-\i
 t_l^{(n)}}\e^{\i s_k^{(n)}}},
\]
and
\[
b_l^{(n)} = \frac 1n \sum_{k=1}^n\frac{\e^{\i s_k^{(n)}}}{\e^{\i
 s_k^{(n)}}-\e^{\i t_l^{(n)}}} =  \frac 1n \sum_{k=1}^n\frac{1}{1-\e^{-\i
 s_k^{(n)}}\e^{\i t_l^{(n)}}}.
\]

It follows that
\[
a_k^{(n)} = b_l^{(n)} = \frac 1n \sum_{m=1}^n \frac{1}{1-\theta_m},
\]

where $\theta_m =\e^{\i \frac{(2m-1)\pi}{n}}$, $m=1, \ldots, n$, are the
roots of the polynomial
\[
p(z)=z^n +1, \qquad (z\in\C).
\]

Consequently,
\[
p(z)=\prod_{m=1}^n (z-\theta_m),
\]
and
\[
\frac{p'(z)}{p(z)}=\sum_{m=1}^n\frac{1}{z-\theta_m}.
\]

This implies that
\[
a_k^{(n)} = b_l^{(n)} = \frac1n \frac{p'(1)}{p(1)}=\frac12,
\]

and we have thus shown that
\begin{equation}
  E_nF_n = {\textstyle \frac 12} (E_n+F_n).
\end{equation}

Now, $E_nF_n = {\textstyle \frac 12} (E_n+F_n)\rightarrow E$
w.r.t. $\|\cdot\|_p$ and hence w.r.t.  $\|\cdot\|_{\frac p2}$. Moreover,
\begin{eqnarray*}
\|E_nF_n-E^2\|_{\frac p2}^{\frac p2} &\leq & \|(E_n-E)F_n\|_{\frac
  p2}^{\frac p2} + \|E(F_n-E)\|_{\frac p2}^{\frac p2}\\
& \leq & \big(\|E_n-E\|_p\|\|F_n\|_p\big)^{\frac p2} +
  \big(\|E\|_p\|\|F_n-E\|_p\big)^{\frac p2}\\
  &\rightarrow 0 & {\rm as} \; n\rightarrow \infty, 
\end{eqnarray*}

and it follows that $E=E^2$.

Clearly, $E_n$ and $T$ commute for every $n\in\N$. Since the map
$(a,b)\mapsto ab$ is continuous w.r.t the measure topology on products of
sets which are bounded in measure, this implies that
\begin{equation}
  ET=TE\; {\rm in} \;\tilde\CM. 
\end{equation}
In particular, ${\rm ker}(E)$ and $R(E)$ (the range of $E$)
are $T$-invariant.

Let $P$ denote the projection onto $\overline{R(E)}$. Then $P\in\CM$, and
$PTP$ agrees with $TP$ on $R(E)+P(\CH)^\bot$ which is a dense subset of
$\CH$.  Hence $PTP=TP$ in $\tilde\CM$, i.e. $\CK=\overline{R(E)}$ is $T$-invariant. 

Next take a fixed $m\in\N$. As above we find that the map
\[
t\mapsto (\e^{\i t}\unit-T)^{-1}\e^{\i (m+1)t}
\]

is H\"older continuous with exponent $\alpha$, and thus we may define
$G\in\Lp$ by
\begin{equation}
  G=\frac{1}{2\pi\i}\int_{\partial B(0,1)}(\lambda\unit
  -T)^{-1}\lambda^m\,\d\lambda.
\end{equation}

We know that $G=\lim_{n\rightarrow\infty}G_n$ in $p$-norm, where
\[
G_n= \frac 1n \sum_{k=1}^n(\e^{\i s_k^{(n)}}\unit
-T)^{-1}\e^{\i(m+1)s_k^{(n)}},
\]
with $s_k^{(n)}=\frac{2\pi k}{n}$. With $E_n$ as in \eqref{eq3-6a} we want
to prove that for fixed $n>m$,
\begin{equation}\label{eq3-16}
  G_n=T^mE_n,
\end{equation}
i.e. that
\begin{equation}\label{eq3-7}
  \frac 1n \sum_{k=1}^n(\e^{\i
  s_k^{(n)}}\unit-T)^{-1}\e^{\i(m+1)s_k^{(n)}}= \frac 1n
  \sum_{k=1}^n(\e^{\i s_k^{(n)}}\unit -T)^{-1}T^m\e^{\i s_k^{(n)}}.
\end{equation}

Note that \eqref{eq3-7} holds if
\begin{equation}\label{eq3-10}
  \frac1n \sum_{k=1}^n(\e^{\i  s_k^{(n)}}-z)^{-1}\e^{\i(m+1)s_k^{(n)}}=\frac1n \sum_{k=1}^n(\e^{\i s_k^{(n)}}
  -z)^{-1}z^m\e^{\i s_k^{(n)}}
\end{equation}
as rational functions of $z\in\C$.

Now, with
\[
q(z)=z^n-1, \qquad (z\in\C),
\]
$\e^{\i s_k^{(n)}}$, $k=1,\ldots, n$ are the roots of $q$, and therefore
\[
{\rm Res}\Big(\frac{1}{q};\e^{\i s_k^{(n)}}\Big)=\frac{1}{q'(\e^{\i
    s_k^{(n)}})}= \frac{1}{n\cdot \e^{\i(n-1) s_k^{(n)}}}= \frac{\e^{\i
    s_k^{(n)}}}{n}.
\]

It follows that
\begin{equation}\label{eq3-8}
  \frac{1}{q(z)}= \frac1n \sum_{k=1}^n(z-\e^{\i
  s_k^{(n)}})^{-1}\e^{\i s_k^{(n)}}.
\end{equation}

Also, $\lim_{|z|\rightarrow\infty}\frac{z^m}{q(z)}=0$ and
\[
{\rm Res}\Big(\frac{z^m}{q(z)};\e^{\i s_k^{(n)}}\Big)=\frac{\e^{\i
    ms_k^{(n)}}}{q'(\e^{\i s_k^{(n)}})}
    = \frac1n \e^{\i (m+1)s_k^{(n)}}.
\]
Hence
\begin{equation}\label{eq3-9}
  \frac{z^m}{q(z)}= \frac1n \sum_{k=1}^n(z-\e^{\i
  s_k^{(n)}})^{-1}\e^{\i (m+1)s_k^{(n)}}.
\end{equation}

Comparing \eqref{eq3-8} with  \eqref{eq3-9} we find that  \eqref{eq3-10}
holds when $n>m$, that is $G_n=T^mE_n$. Take the limit as
$n\rightarrow\infty$ on both sides of \eqref{eq3-16} (w.r.t. the measure topology) and
conclude that
\begin{equation}\label{eq3-11}
  G=T^mE.
\end{equation}

\eqref{eq3-11} enables us to make an estimate of $\|T^mE\|_p$: 
As in \eqref{eq3-12} one can show that
\begin{eqnarray*}
\|(\e^{\i s}\unit - T)^{-1}\e^{\i(m+1) s}- (\e^{\i t}\unit - T)^{-1}\e^{\i
  (m+1)t}\|_p^p&\leq& 
C^p|\e^{\i s}-\e^{\i t}|^{\alpha p} + (C')^p|\e^{\i (m+1)s}-\e^{\i(m+1)
  t}|^p\\
&\leq & C^p|s-t|^{\alpha p} + (C')^p|\e^{\i (m+1)s}-\e^{\i(m+1)
  t}|^p,
\end{eqnarray*}

and since $\alpha\leq 1$ and $|\e^{\i(m+1) s}-\e^{\i(m+1) t}|\leq 2$,
\begin{eqnarray*}
|\e^{\i(m+1) s}-\e^{\i (m+1)t}|&\leq& 2\Big|\frac{\e^{\i (m+1)s}-\e^{\i
 (m+1)t}}{2}\Big|^\alpha\\
&=& 2^{1-\alpha}|\e^{\i (m+1)s}-\e^{\i (m+1)t}|^\alpha \\
&\leq& 2^{1-\alpha}(m+1)^\alpha |s-t|^\alpha\\
&\leq &  2^{1-\alpha}(m+1) |s-t|^\alpha.
\end{eqnarray*}

It follows that $t\mapsto  (\e^{\i t}\unit - T)^{-1}\e^{\i (m+1)t}$ is
H\"older continuous with exponent $\alpha$ and constant
\[
\big(C^p+ (C')^p(m+1)^{p}2^{p(1-\alpha)}\big)^{\frac 1p}.
\]

Since $0<p<1$,
\[
\Big(\frac{x+y}{2}\Big)\geq \frac12 (x^p+y^p),
\]
whence
\[
x^p+y^p \leq 2^{1-p}(x+y)^p \leq 2^p(x+y)^p,
\]
because $p\geq \frac 12$. Consequently, $t\mapsto  (\e^{\i t}\unit - T)^{-1}\e^{\i (m+1)t}$ is
H\"older continuous with exponent $\alpha$ and constant
\begin{equation}\label{eq3-13}
C''=2(C+ C'(m+1)2^{1-\alpha}).
\end{equation}

Then according to Lemma~\ref{c},
\begin{equation}\label{eq3-14}
\|G\|_p^p \leq \frac{(C'')^p(2\pi)^{\alpha p}}{\alpha p +p-1} + \|(\unit-T)^{-1}\|_p^p,
\end{equation}

and combining \eqref{eq3-13} with \eqref{eq3-14} we find that
there exist positive constants $k_1, k_2$
(which are independent of $m$) such that
\[
\|T^mE\|_p\leq k_1 + mk_2, \qquad (m\in\N).
\]

In particular,
\[
\lim_{m\rightarrow \infty}\Big(\tau(P)^{-\frac1p}\|T^mE\|_p\Big)^{\frac 1m}\leq 1.
\]

Since $\CK$ is $T$-invariant, with $T_{11}=T_{|_\CK}$, we have that
$T_{11}^m=T^m_{|_\CK}$. Moreover, $TE$ and $T_{|_\CK}$ agree on $R(E)$, and hence
they must agree as operators in $(P\CM P)\tilde{} $. It follows that
\[
T_{11}^m=T^m_{|_\CK} =T^mE_{|_\CK},
\]
and then, by Proposition~\ref{modspecrad},
\begin{eqnarray*}
  r'(T_{11}) & \leq  &  \lim_{m\rightarrow\infty}\|T_{11}^m\|_p^{\frac 1m}\\
  & = & \lim_{m\rightarrow\infty}\Big(\tau(P)^{-\frac1p}\|T^mE\|_p\Big)^{\frac 1m}\\
  &\leq & 1.
\end{eqnarray*}
Thus $\supp(\mu_{T_{11}})\subseteq \overline{B(0,1)}$.

Finally, to prove that $\supp(\mu_{T_{22}})\subseteq
\{z\in\C\,|\,|z|\geq 1\}$, define for fixed $m\in\N$ $H\in\Lp$ by
\begin{equation}
  H=-\frac{1}{2\pi\i}\int_{\partial B(0,1)}(\lambda\unit
  -T)^{-1}\lambda^{-m}\,\d\lambda
\end{equation}
As above one may prove that the right hand side makes sense as a Riemann
integral in $\Lp$. Moreover, $H=\lim_{n\rightarrow\infty}H_n$ in $p$-norm,
where
\[
H_n= -\frac 1n \sum_{k=1}^n(\e^{\i s_k^{(n)}}\unit
-T)^{-1}\e^{\i(1-m)s_k^{(n)}},
\]
with $s_k^{(n)}=\frac{2\pi k}{n}$. We are going to prove that
\begin{equation}\label{eq3-17}
  H_n=T^{-m}(1-E_n), \qquad(n>m),
\end{equation}
and as in the above, \eqref{eq3-17} holds, as soon as
\begin{equation}\label{eq3-21}
-\frac 1n \sum_{k=1}^n(\e^{\i s_k^{(n)}}
-z)^{-1}\e^{\i(1-m)s_k^{(n)}}= z^{-m}\Big(1-\frac 1n \sum_{k=1}^n(\e^{\i
s_k^{(n)}}\unit -z)^{-1}\e^{\i s_k^{(n)}}\Big), \qquad (z\in\C).
\end{equation}  

According to \eqref{eq3-8},
\begin{equation}\label{eq3-20}
  z^{-m}\Big(1-\frac 1n \sum_{k=1}^n(\e^{\i
-s_k^{(n)}}\unit -z)^{-1}\e^{\i s_k^{(n)}}\Big)=z^{-m}\Big(1+\frac{1}{q(z)}\Big)=\frac{z^{n-m}}{z^n-1}.
\end{equation}  

Since $\frac{z^{n-m}}{z^n-1}\rightarrow 0$ as $|z|\rightarrow \infty$ and
\[
{\rm Res}\Big(\frac{z^{n-m}}{z^n-1}; \e^{\i s_k^{(n)}}\Big)=
\frac{\e^{(n-m)\i s_k^{(n)}}}{q'(\e^{\i s_k^{(n)}})}= \frac1n \e^{\i
  (1-m)s_k^{(n)}},
\]

$\frac{z^{n-m}}{z^n-1}$ has the following partial fraction decomposition:
\begin{equation}\label{eq3-19}
\frac{z^{n-m}}{z^n-1} =\frac1n \sum_{k=1}^n \e^{\i
  (1-m)s_k^{(n)}} \frac{1}{z-\e^{\i s_k^{(n)}}}= -\frac1n \sum_{k=1}^n \e^{\i
  (1-m)s_k^{(n)}} \frac{1}{\e^{\i s_k^{(n)}}-z}
\end{equation}

Combining \eqref{eq3-20} with \eqref{eq3-19} we find that \eqref{eq3-21}
holds, and hence $H_n=T^{-m}(1-E_n)$ when $n>m$. Taking the limit as
$n\rightarrow\infty$ on both sides of \eqref{eq3-17} (w.r.t. the measure
topology) we arrive at the identity
\begin{equation}
  H=T^{-m}(\unit-E).
\end{equation}

As in the above this implies that there are constants $k_3, k_4>0$ such
that
\[
\|T^{-m}(\unit-E)\|_p \leq k_3 +k_4 m, \qquad (m\in\N),
\]
whence
\[
\|(T\cc)^{-m}(\unit-E\cc)\|_p = \|(\unit-E)T^{-m}\|_p =
\|T^{-m}(\unit-E)\|_p \leq k_3 +k_4 m,
\]
because $T$ and $E$ commute.

Clearly, $(\unit-E\cc)^2=\unit-E\cc$, and $(T\cc)^{-1}$ commutes with
$\unit-E\cc$. As in the above this implies that $\overline{R(\unit-E\cc)}$
is $(T\cc)^{-1}$-invariant and that
\begin{equation}\label{eq3-22}
\|(T\cc)^{-m}Q\|_p\leq k_3 +k_4 m,
\end{equation}
where $Q\in\CM$ denotes the projection onto $\overline{R(\unit-E\cc)}$. Now
\[
\overline{R(\unit-E\cc)}^\bot=\overline{{\rm
    ker}(\unit-E)}=\overline{\{x\in\CD(E)\,|\,Ex=x\}}=\overline{R(E)}=\CK,
\]
and hence
\begin{equation}\label{eq3-23}
(T\cc)^{-m}Q = Q(T\cc)^{-m}Q=(Q(T\cc)^{-1}Q)^m = (T_{22}\cc)^{-m}.
\end{equation}

As above it follows from \eqref{eq3-22} and \eqref{eq3-23} that
\[
r'(T_{22}^{-1})\leq \lim_{m\rightarrow\infty}\|T_{22}^{-m}\|_p^{\frac
  1m}=\lim_{m\rightarrow\infty}\|(T_{22}\cc)^{-m}\|_p^{\frac 1m}\leq 1,
\]
i.e. $\supp(\mu_{T_{22}^{-1}})\subseteq \{z\in\C\,|\,|z|\leq 1\}$. Since $T\in\CMD$, it follows that $T_{22}\in (P^\bot \CM P^\bot)^\Delta$, and thus, by \cite[Proposition~2.16]{HS}, $\supp(\mu_{T_{22}})\subseteq \{z\in\C\,|\,|z|\geq 1\}$. $\endproof$

%% file: invariant.tex
\section{The invariant subspace problem relative to a II$_1$-factor}

The purpose of the present section is to combine the results of the
previous sections to give a proof of our main theorem:

\begin{thm}\label{Borelsets} For every $T\in\CM$ and every Borel set
  $B\subseteq \C$ there is a largest closed $T$-invariant subspace $\CK =
  \CK_T(B)$ affiliated with $\CM$, such that the Brown measure of $T|_\CK$,
  $\mu_{T|_\CK}$, is concentrated on $B$.\footnote{If $\CK = \{0\}$, then
  we define $\mu_{T|_\CK}:=0$. If $\CK \neq \{0\}$,  then $\mu_{T|_\CK}$ is computed
  relative to the II$_1$-factor $P\CM P$, where $P\in\CM$ denotes the
  projection onto $\CK$.} Moreover, $\CK$ is hyperinvariant\footnote{$\CK$ is said to be 
  {\it hyperinvariant} for $T$ if it is invariant under every operator
  $S\in\{T\}'$.} for $T$, and if $P=P_T(B)\in\CM$ denotes the
  projection onto $\CK$, then
  \begin{itemize}
    \item[(i)] $\tau(P)=\mu_T(B)$,
    \item[(ii)] the Brown measure of $P^\bot TP^\bot$, considered as an
    element of $P^\bot \CM P^\bot$, is concentrated on $\C\setminus B$.
  \end{itemize}
\end{thm}

\vspace{.2cm}

\begin{cor}\label{partial solution} Let $T\in \CM$, and suppose that $\mu_T$ is not concentrated on
  a singleton. Then there is a non-trivial subspace
  affiliated with $\CM$ which is hyperinvariant for $T$.
\end{cor}

\proof Take a Borel set $B\subseteq \C$ such that $\mu_T(B)>0$ and
$\mu_T(B^c)>0$. Then with $P=P_T(B)\in\CM$ as in Theorem~\ref{Borelsets}, one has
that $P$ is hyperinvariant for $T$, and $\tau(P)=\mu_T(B)\in (0,1)$. Since
$\tau$ is faithful, $P$ must be non-trivial. $\endproof$

\vspace{.2cm}

\subsection*{First case: Spectral subspaces associated with the set $B(0,1))$.}

Consider a fixed operator $T\in\CM$. As in sections~4 and 5, we can choose a circular system $\{x,y\}$ in $L(\F_4)$ and embed $\CM$ and $L(\F_4)$ into the free product $\CN=\CM\ast L(\F_4)$. This way $\{x,y\}$ becomes a circular system in the II$_1$--factor $\CN$ such that $\CM$ is free from $\{x,y\}$. We will denote the trace on $\CN$ by $\tau$ as
well. For each $n\in\N$ we define an
operator $T_n\in \LpN$, ($0<p<1$), by
\begin{equation}
  T_n = T+{\textstyle \frac 1n} xy^{-1}.
\end{equation}

Note that $\|T-T_n\|_p\rightarrow 0$ as $n\rightarrow\infty$. According to
Theorem~\ref{projection} and Theorem~\ref{Lipschitz} we have:

\begin{thm} \label{P_n} For each $n\in\N$ there exists a projection
  $P_n\in\CN$ such that
  \begin{itemize}
    \item[(i)]$P_nT_nP_n = T_nP_n$,
    \item[(ii)]  $\mu_{P_nT_nP_n}$ (computed relative to $P_n\CN P_n$) is
    concentrated on $\overline{B(0,1)}$, 
    \item[(iii)]  $\mu_{P_n^\bot T_nP_n^\bot}$ (computed relative to
    $P_n^\bot \CN P_n^\bot$) is concentrated on $\C\setminus B(0,1)$. 
  \end{itemize}
\end{thm}

\proof We will show that $T_n$ satisfies the assumptions of Theorem~\ref{projection} with $\alpha=1$ and $p\in (\frac12, \frac23)$. Clearly, $T_n\in \CN^\Delta$, and by Theorem~\ref{Lipschitz}~(i), 
\[
(T_n-\lambda\unit)^{-1} = n(nT-n\lambda\unit +z)^{-1}\in L^p(\CN,\tau)
\]
for $p\in (0,1)$, $\lambda\in\C$. Moreover, by Theorem~\ref{Lipschitz}~(ii),
\[
\begin{split}
\|(T_n-\lambda\unit)^{-1}-(T_n-\mu\unit)^{-1}\|_p &= n \|(nT-n\lambda\unit +z)^{-1}- (nT-n\mu\unit +z)^{-1}\|_p\\
&\leq n^2C_p^{(1)}|\lambda-\mu|
\end{split}
\]
for $p\in (0,\frac23)$ and $\lambda,\mu\in\C$. Thus, Theorem~\ref{P_n} follows from Theorem~\ref{projection}. $\endproof$

\vspace{.5cm}

\begin{prop}\label{prop20.7.4} Let $n\in\N$, and let $P_n$ be as in Theorem~\ref{P_n}. Then
\begin{equation}\label{20.7.2}
0<\tau(P_n)<1,
\end{equation}
\begin{equation}\label{20.7.3}
\mu_{T_n}(\partial B(0,1))=0,
\end{equation}
\begin{equation}
\mu_{T_n}(B(0,1))=\tau(P_n).
\end{equation}
Moreover, for every Borel set $A\in\B(\C)$, 
\begin{equation}\label{20.7.5}
\mu_{T_n}(A\cap B(0,1))=\tau(P_n)\mu_{P_nT_nP_n}(A),
\end{equation}
\begin{equation}\label{20.7.6}
\mu_{T_n}(A\setminus B(0,1))=\tau(P_n^\bot)\mu_{P_n^\bot T_nP_n^\bot}(A),
\end{equation}
where the Brown measures of $P_nT_nP_n$ and $P_n^\bot T_nP_n^\bot$ are computed relative to $P_n\CN P_n$ and $P_n^\bot \CN P_n^\bot$, respectively. 
\end{prop}

\proof According to Corollary~\ref{cor4.6}, $\mu_{T_n}$ has a density $\varphi_n$ w.r.t. Lebesgue measure on $\C$, and $\varphi_n(z)>0$ for all $z\in\C$. Therefore $\supp(\mu_{T_n})=\C$. Moreover,
\begin{equation}\label{20.7.7}
\mu_{T_n}=\tau(P_n)\mu_{P_nT_nP_n} + \tau(P_n^\bot)\mu_{P_n^\bot T_nP_n^\bot}
\end{equation}
(cf. Proposition~\ref{prop6.9}). It then follows from (ii) and (iii) of Theorem~\ref{P_n} that \eqref{20.7.2} holds. Moreover, since $\partial B(0,1)$ is a nullset w.r.t. Lebesgue measure on $\C$, \eqref{20.7.3} holds. Then by (ii) of Theorem~\ref{P_n},
\begin{eqnarray*}
\tau(P_n) &=&  \tau(P_n)\mu_{P_nT_nP_n}(B(0,1))\\
&=& \mu_{T_n}(B(0,1))- \tau(P_n^\bot) \mu_{P_n^\bot T_nP_n^\bot}(B(0,1))\\
&\overset{{\rm Theorem~\ref{P_n}~(iii)}}{=}& \mu_{T_n}(B(0,1)).
\end{eqnarray*}
\eqref{20.7.5} and \eqref{20.7.6} now follow from \eqref{20.7.7} and the fact that $\mu_{P_nT_nP_n}$ and $\mu_{P_n^\bot T_nP_n^\bot}$ are concentrated on $B(0,1)$ and $\C\setminus B(0,1)$, respectively. $\endproof$

\vspace{.5cm}

\noindent Now, take a free ultrafilter $\omega$ on $\N$, and let
\begin{eqnarray*}
\CJ_\omega &=& \{(x_n)_{n=1}^\infty \in \ell^\infty(\CN)\,|\, \lim_\omega
\|x_n\|_2 = 0\},\\
\CN^\omega &=& \ell^\infty(\CN)/\CJ_\omega.
\end{eqnarray*}

Moreover, let $\rho: \ell^\infty(\CN)\rightarrow \CN^\omega$ denote the
quotient mapping, let $\tilde{T}=\rho((T)_{n=1}^\infty)$ denote the copy of $T$ in $\CN^\omega$, and
define a projection $P\in\CN^\omega$ by
\[
P=\rho((P_n)_{n=1}^\infty).
\]

Recall that the ultrapower $\CN^\omega$ is a II$_1$-factor 
equipped with a faithful tracial state $\tau_\omega$ given by
\[
\tau_\omega(\rho(x))=\lim_{n\rightarrow\omega} \tau(x_n), \qquad
x=(x_n)_{n=1}^\infty\in \ell^\infty(\CN).
\]
For $0<p<\infty$,
\[
\|\rho(x)\|_p = \lim_{n\rightarrow\omega}\|x_n\|_p, \qquad x=(x_n)_{n=1}^\infty\in \ell^\infty(\CN),
\]
where the $p$--norm on the left--hand side is computed w.r.t. the trace $\tau_\omega$.

\begin{prop}\label{T-inv} With $P$ and $\tilde{T}$ as defined above,
  $P\tilde{T}P=\tilde{T}P$.
\end{prop}

\proof For all $p\in (0,1)$,
\begin{eqnarray*}
  \|P\tilde{T}P-\tilde{T}P\|_p^p & = & \lim_{n\rightarrow\omega}
  \|P_nTP_n-TP_n\|_p^p \\
  & = & \lim_{n\rightarrow\omega} \|(T-T_n)P_n +T_nP_n -P_nT_nP_n
  -P_n(T-T_n)P_n\|_p^p,
  \end{eqnarray*}
where
\begin{eqnarray*}
\|(T-T_n)P_n +T_nP_n -P_nT_nP_n
  -P_n(T-T_n)P_n\|_p^p &\leq &2\|T-T_n\|_p^p\\
   & \rightarrow& 0 \quad {\rm as \;} n\rightarrow \infty.
\end{eqnarray*}
Since $\omega$ is a free ultrafilter, any convergent sequence converges
along $\omega$ as well. Hence,
\[
\|P\tilde{T}P-\tilde{T}P\|_p=\lim_{n\rightarrow\omega}\|P_nTP_n-TP_n\|_p= 0, \qquad 0<p<1,
\]
and therefore $\tilde T P =P\tilde T P$. $\endproof$

\vspace{.2cm}

\begin{lemma}\label{MSRI4.5}
  Let $\mu, \mu_1, \mu_2, \ldots\in \Prob(\C)$, and suppose that $\mu_n
  \rightarrow \mu$ weakly as $n\rightarrow \infty$. Then
  \begin{itemize}
    \item[(i)] For every open set $\CU\subseteq \C$,
      \[
      \mu(\CU)\leq \liminf_{n\rightarrow\infty}\mu_n(\CU).
      \]
    \item[(ii)] For every closed set $F\subseteq \C$,
      \[
      \mu(F)\geq \limsup_{n\rightarrow\infty}\mu_n(F).
      \]
    \item[(iii)] For every Borel set $B\subseteq \C$ with $\mu(\partial
    B)=0$,
    \[
    \mu(B)=\lim_{n\rightarrow\infty}\mu_n(B).
    \]
  \end{itemize}
\end{lemma}

\proof This is standard. In order to prove (i), write $1_\CU$ as the
pointwise limit of an increasing sequence of non-negative functions
$\phi_n\in C_c(\C)$. To prove (ii), use (i) and the fact that $\mu(F)=
1-\mu(F^c)$, where $F^c$ is open. (iii) follows from (i) and (ii) with
$\CU= {\rm int}B$ (the interior of $B$) and $F=\overline B$. Note that
$\mu(F)=\mu(\CU)= \mu(B)$. $\endproof$

\vspace{.2cm}

\begin{prop}\label{measureofball} If $\mu_T(\partial B(0,1))=0$, then the
sequence $(\tau(P_n))_{n=1}^\infty$ converges as $n\rightarrow\infty$. Moreover,
\begin{equation}
\mu_T(B(0,1))=\lim_{n\rightarrow\infty}\tau(P_n)=\tau_\omega(P).
\end{equation}
\end{prop}

\proof According to Corollary~\ref{cor4.6}, $\mu_{T_n}$converges weakly to $\mu_T$ as $n$
tends to infinity. Then by Proposition~\ref{prop20.7.4} and Lemma~\ref{MSRI4.5},
\begin{eqnarray*}
\mu_T(B(0,1)) &=&
\lim_{n\rightarrow\infty}\mu_{T_n}(B(0,1))\\
&=& \lim_{n\rightarrow\infty}\tau(P_n),
\end{eqnarray*}
and since $\omega$ is a free ultrafilter on $\N$, $\tau(P_n)$ converges along $\omega$ to $\mu_T(B(0,1))$ as well. Thus,
\[
\tau_\omega(P)= \lim_{n\rightarrow\omega}\tau(P_n)=\mu_T(B(0,1))  . \quad \endproof
\]

\vspace{.2cm}

\begin{lemma}\label{infimum}\begin{itemize}
\item[(i)] Consider a (classical) probability space $(X,
  \CE, \mu)$. For every $f\in \bigcup_{p>0}L^p(X, \CE, \mu)$ one
  has that
  \begin{equation*}
    \exp\Big\{\int_X \log|f|\,\d\mu\Big\}=\inf_{p>0}\|f\|_p =
    \lim_{p\rightarrow 0+}\|f\|_p,
  \end{equation*}
where $\exp(-\infty):=0$. 
\item[(ii)] For every $S\in\bigcup_{p>0}L^p(\CN,\tau)$,
\[
\Delta(S)=\inf_{p>0}\|S\|_p=\lim_{p\rightarrow 0+}\|S\|_p.
\]
\end{itemize}
\end{lemma}

\proof (i) follows from \cite[Section~3, exercise~5(d)]{Ru}. (ii) now follows by application of (i) to $(\R,\B, \mu_{|S|})$ and $f={\rm id}_{\R}$. $\endproof$

\vspace{.2cm}

\begin{lemma}\label{limsupleq} Suppose $\mu_T(\partial B(0,1))=0$ and that $\tau_\omega(P)\in (0,1)$. Then for every complex number $\lambda$,
  \begin{equation}
    \lim_{n\rightarrow\omega}\Delta_{P_n\CN P_n}(P_nT_nP_n-\lambda P_n)\leq
    \Delta_{P\CN^\omega P}(P\tilde TP-\lambda P),
  \end{equation}
  and
  \begin{equation}\label{eq2-1}
    \lim_{n\rightarrow\omega}\Delta_{P_n^\bot \CN P_n^\bot}(P_n^\bot T_nP_n^\bot-\lambda P_n^\bot)\leq \Delta_{P^\bot\CN^\omega P^\bot}(P^\bot\tilde T P^\bot-\lambda P^\bot).
  \end{equation}
\end{lemma}

\proof We consider the case $\lambda =0$ only. The general case can be taken care of in the same way. According to
Lemma~\ref{infimum}
\[
\Delta_{P_n\CN P_n}(P_nT_nP_n) =
\inf_{0<p<1}\|P_nT_nP_n\|_{L^p(P_n\CN P_n)},
\]
where
\[
\|P_nT_nP_n\|_{L^p(P_n\CN P_n)}^p = \frac{1}{\tau(P_n)} \|P_nT_nP_n\|_p^p.
\]
It follows that
\begin{equation}\label{eq2-21}
\Delta_{P_n\CN P_n}(P_nT_nP_n) = \inf_{0<p<1}
(\tau(P_n)^{-\frac1p}\|P_nT_nP_n\|_p),
\end{equation}
and by the same argument, 
\begin{equation}\label{eq2-18}
\Delta_{P\CN^\omega P}(P\tilde T P) =
\inf_{0<p<1}(\tau_\omega(P)^{-\frac1p}\|P\tilde T P\|_p),
\end{equation}
where
\begin{equation}\label{eq2-19}
\tau_\omega(P)^{-\frac1p}=\lim_{n\rightarrow\omega}\tau(P_n)^{-\frac1p},
\end{equation}
and
\[
\|P\tilde T P\|_p= \lim_{n\rightarrow\omega}\|P_n T P_n\|_p, \qquad 0<p<1.
\]
Now,
\[
\|P_nTP_n-P_nT_nP_n\|_p \leq \|T-T_n\|_p \rightarrow 0 \;{\rm as }\;
n\rightarrow\omega,
\]
and hence
\begin{equation}\label{eq2-20}
\|P\tilde T P\|_p=\lim_{n\rightarrow\omega}\|P_n T P_n\|_p = \lim_{n\rightarrow\omega}\|P_n T_n P_n\|_p. 
\end{equation}
Combining now \eqref{eq2-18}, \eqref{eq2-19} and \eqref{eq2-20}, we obtain:
\[
\Delta_{P\CN^\omega P}(P\tilde T P) =
\inf_{0<p<1}(\lim_{n\rightarrow\omega}\tau(P_n)^{-\frac1p}\|P_n T_n
P_n\|_p).
\]
According to \eqref{eq2-21}, for every $p\in(0,1)$,
\[
\lim_{n\rightarrow\omega}(\tau(P_n)^{-\frac1p}\|P_n T_n
P_n\|_p) \geq \lim_{n\rightarrow\omega}\Delta_{P_n\CN P_n}(P_nT_nP_n),
\]
and it follows that
\[
\Delta_{P\CN^\omega P}(P\tilde T P) \geq \lim_{n\rightarrow\omega}\Delta_{P_n\CN P_n}(P_nT_nP_n).
\]

By similar arguments one can prove that \eqref{eq2-1} holds. $\endproof$

\vspace{.2cm}

\begin{lemma}\label{limsup} Let $(a_n)_{n=1}^\infty$ and  $(b_n)_{n=1}^\infty$ be
  sequences in $\R\cup \{-\infty\}$ and let $a,b\in\R$. 
  \begin{itemize}
  \item[(i)] If 
  \begin{eqnarray*}
  \limsup_{n\rightarrow\infty}a_n&\leq & a,\\
  \limsup_{n\rightarrow\infty}b_n&\leq & b,\\
  \liminf_{n\rightarrow\infty}(a_n+b_n)&\geq &a+b,
  \end{eqnarray*}
  then
  \[
  \lim_{n\rightarrow\infty}a_n = a \qquad {\rm and} \qquad 
  \lim_{n\rightarrow\infty}b_n = b.
  \]
  \item[(ii)] (i) holds for convergence along a free ultrafilter $\omega$ on $\N$ as well. 
  \end{itemize}
\end{lemma}

\proof (i) Let $\eps>0$. Then, eventually as $n\rightarrow \infty$, $a_n\leq a+\eps$, $b_n\leq b+\eps$, and $a_n+b_n\geq a+b-\eps$. Hence,
\[
a_n=(a_n+b_n)-b_n\geq a-2\eps
\]
and 
\[
b_n=(a_n+b_n)-a_n\geq b-2\eps,
\]
eventually as $n\rightarrow\infty$. It follows that $ \lim_{n\rightarrow\infty}a_n = a$ and $\lim_{n\rightarrow\infty}b_n = b$. (ii) follows in a similar way. $\endproof$

\vspace{.2cm}

\begin{prop}\label{lim} Suppose $\mu_T(\partial B(0,1))=0$ and $\tau_\omega(P)\in(0,1)$. Let $\lambda\in\C$ with $\Delta(T-\lambda\unit)>0$. Then
 \begin{equation}\label{eq2-11}
    \lim_{n\rightarrow\omega}\Delta_{P_n\CN P_n}(P_nT_nP_n-\lambda P_n)= \Delta_{P\CN^\omega P}(P\tilde T
    P-\lambda P),
  \end{equation}
  and
  \begin{equation}\label{eq2-12}
    \lim_{n\rightarrow\omega}\Delta_{P_n^\bot \CN P_n^\bot}(P_n^\bot T_nP_n^\bot-\lambda
    P_n^\bot)= \Delta_{P^\bot \CN^\omega P^\bot}(P^\bot\tilde T P^\bot-\lambda P^\bot).
  \end{equation}
\end{prop}

\proof We will apply Lemma~\ref{limsup} with
\begin{eqnarray*}
a_n&=&\tau(P_n)\log \Delta_{P_n\CN P_n}(P_nT_nP_n-\lambda P_n),\\
b_n&=&\tau(P_n^\bot)\log \Delta_{P_n^\bot\CN P_n^\bot}(P_n^\bot T_nP_n^\bot-\lambda P_n^\bot),\\
a&=& \tau_\omega(P)\log \Delta_{P\CN^\omega P}(P\tilde T P-\lambda P),\\
b&=& \tau_\omega(P^\bot)\log \Delta_{P^\bot\CN^\omega P^\bot}(P^\bot\tilde T P^\bot-\lambda P^\bot).
\end{eqnarray*}
Since $\tau(P_n)$ converges to $\tau_\omega(P)$ as $n\rightarrow\omega$, we get from Lemma~\ref{limsupleq} that
\[
\lim_{n\rightarrow\omega}a_n\leq a \qquad {\rm and} \qquad \lim_{n\rightarrow\omega}b_n\leq b.
\]
Moreover, by Proposition~\ref{prop6.9},
\[
a_n+b_n=\log\Delta(T_n-\lambda\unit) \qquad {\rm and} \qquad a+b=\log\Delta(T-\lambda\unit).
\]
Hence, by the proof of Corollary~\ref{cor4.8},
\[
\lim_{n\rightarrow\omega}(a_n+b_n)=a+b.
\]
$a,b\in\R$ because $\Delta(T-\lambda\unit)>0$. Then by Lemma~\ref{limsup}~(ii),
\[
  \lim_{n\rightarrow\omega}a_n = a \qquad {\rm and} \qquad 
  \lim_{n\rightarrow\omega}b_n = b.
  \]
Since $\lim_{n\rightarrow\omega}\tau(P_n)=\tau_\omega(P)>0$ and $\lim_{n\rightarrow\omega}\tau(P_n^\bot)=\tau_\omega(P^\bot)>0$, this completes the proof. $\endproof$

\vspace{.2cm}

We will now prove that under the assumptions $\mu_T(\partial B(0,1))=0$ and $0<\tau_\omega(P)<1$,
\[
\lim_{n\rightarrow\omega}\mu_{P_nT_nP_n}=\mu_{P\tilde TP}
\]
and
\[
\lim_{n\rightarrow\omega}\mu_{P_n^\bot T_nP_n^\bot}=\mu_{P^\bot\tilde TP^\bot}
\]
(weak convergence in ${\rm Prop}(\C)$). 

\vspace{.2cm}

\begin{lemma}\label{lemma20.7.13} Suppose $\mu_T(\partial B(0,1))=0$ and $0<\tau_\omega(P)<1$. Define maps $\rho,\sigma:\B(\C)\rightarrow [0,\infty[$ by
\begin{equation}\label{20.7.11}
\rho(A)=\frac{1}{\tau_\omega(P)}\mu_T(A\cap B(0,1)),
\end{equation}
\begin{equation}\label{20.7.12}
\sigma(A)=\frac{1}{\tau_\omega(P^\bot)}\mu_T(A\setminus B(0,1)).
\end{equation}
Then $\rho,\sigma\in {\rm Prob}(\C)$, and
\[
\lim_{n\rightarrow\infty}\mu_{P_nT_nP_n}=\rho \qquad {\rm and}\qquad \lim_{n\rightarrow\infty}\mu_{P_n^\bot T_nP_n^\bot}=\sigma
\]
(weak convergence in ${\rm Prob}(\C)$).
\end{lemma}

\proof $\rho,\sigma\in {\rm Prob}(\C)$ because $\tau_\omega(P)=\mu_T(B(0,1))$ and $\tau_\omega(P^\bot)=\mu_T(\C\setminus B(0,1))$ (cf. Proposition~\ref{measureofball}). Put
\[
\rho_n=\mu_{P_nT_nP_n} \qquad {\rm and} \qquad \sigma_n= \mu_{P_n^\bot T_nP_n^\bot}.
\]
Then by \eqref{20.7.5} and \eqref{20.7.6},
\[
\rho_n(A)=\frac{1}{\tau(P_n)}\mu_{T_n}(A\cap B(0,1)) \qquad {\rm and} \qquad \rho_n(A)=\frac{1}{\tau(P_n^\bot)}\mu_{T_n}(A\setminus B(0,1)).
\]
Hence, for every continuous function $\varphi:\C\rightarrow [0,1]$,
\[
\int_\C \varphi\,\d\rho_n = \frac{1}{\tau(P_n)}\int_{B(0,1)} \varphi\,\d\mu_{T_n},
\]
and by the definition of $\rho$,
\[
\int_\C \varphi\,\d\rho = \frac{1}{\tau_\omega(P)}\int_{B(0,1)} \varphi\,\d\mu_T.
\]
Choose $f_k\in C_c(\C)$, $0\leq f_k\leq 1$, such that $f_k \nearrow 1_{B(0,1)}$ as $k\rightarrow\infty$. Since $\mu_{T_n}\rightarrow\mu_T$ weakly as $n\rightarrow \infty$,
\begin{eqnarray*}
\liminf_{n\rightarrow\infty}\int_\C\varphi\,\d\rho_n &\geq & \lim_{n\rightarrow\infty}\Big(\frac{1}{\tau(P_n)}\int_\C\varphi \cdot f_k\,\d\mu_{T_n}\Big)\\
&=& \frac{1}{\tau_\omega(P)}\int_\C\varphi\cdot f_k\,\d\mu_T, \qquad k\in\N.
\end{eqnarray*}
Letting $k\rightarrow\infty$, we get that
\begin{equation}\label{20.7.13}
\liminf_{n\rightarrow\infty}\int_\C\varphi\,\d\rho_n \geq \frac{1}{\tau_\omega(P)}\int_\C\varphi\,\d\mu_T=\int_\C \varphi\,\d\rho.
\end{equation}
The same argument applied to the function $1-\varphi$ gives
\[
\liminf_{n\rightarrow\infty}\int_\C(1-\varphi)\,\d\rho_n \geq \int_\C (1-\varphi)\,\d\rho,
\]
and hence
\begin{equation}\label{20.7.14}
\limsup_{n\rightarrow\infty}\int_\C\varphi\,\d\rho_n\leq \int_\C\varphi\,\d\rho.
\end{equation}
Combining \eqref{20.7.13} and \eqref{20.7.14}, we find that for every continuous function $\varphi:\C\rightarrow [0,1]$,
\[
\lim_{n\rightarrow\infty}\int_\C\varphi\,\d\rho_n = \int_\C\varphi\,\d\rho,
\]
proving that $\rho_n\rightarrow\rho$ weakly as $n\rightarrow\infty$.  Since $\mu_{T_n}(\partial B(0,1))=\mu_T(\partial B(0,1))=0$, a similar proof gives that $\sigma_n\rightarrow \sigma$ weakly as $n\rightarrow\infty$: Simply replace $f_k$ in the above by $g_k\in C_c(\C)$ such that $0\leq g_k\leq 1$ and $g_k\nearrow 1_{\C\setminus \overline{B(0,1)}}$ as $k\rightarrow\infty$. $\endproof$

\vspace{.2cm}

\begin{prop}\label{prop20.7.14} Suppose $\mu_T(\partial B(0,1))=0$ and $0<\tau_\omega(P)<1$, and let $\rho,\sigma$ be as in Lemma~\ref{lemma20.7.13}. Then for all $\lambda\in\C$,
\begin{equation}\label{20.7.15}
\limsup_{n\rightarrow\infty} \log\Delta_{P_n\CN P_n}(P_nT_nP_n-\lambda P_n)\leq \int_\C \log|z-\lambda|\,\d\rho(z),
\end{equation}
and
\begin{equation}\label{20.7.16}
\limsup_{n\rightarrow\infty} \log\Delta_{P_n^\bot\CN P_n^\bot}(P_n^\bot T_nP_n^\bot-\lambda P_n^\bot)\leq \int_\C \log|z-\lambda|\,\d\sigma(z).
\end{equation}
Moreover, if $\Delta(T-\lambda\unit)>0$, then
\begin{equation}\label{20.7.17}
\lim_{n\rightarrow\infty} \log\Delta_{P_n\CN P_n}(P_nT_nP_n-\lambda P_n)= \int_\C \log|z-\lambda|\,\d\rho(z),
\end{equation}
and
\begin{equation}\label{20.7.18}
\lim_{n\rightarrow\infty} \log\Delta_{P_n^\bot\CN P_n^\bot}(P_n^\bot T_nP_n^\bot-\lambda P_n^\bot)= \int_\C \log|z-\lambda|\,\d\sigma(z).
\end{equation}
\end{prop}

\proof At first note that by \eqref{20.7.11} and \eqref{20.7.12}, $\supp(\rho),\supp(\sigma)\subseteq\supp(\mu_T)\subseteq\sigma(T)$, whence $\supp(\rho)$ and $\supp(\sigma)$ are compact. Therefore, the right--hand sides of \eqref{20.7.15} and \eqref{20.7.16} are well--defined. Let $\rho_n=\mu_{P_nT_nP_n}$ and $\sigma_n=\mu_{P_n^\bot T_n P_n^\bot}$. Then by Lemma~\ref{lemma20.7.13},
\begin{equation}\label{20.7.19}
\lim_{n\rightarrow\infty}\rho_n=\rho \qquad {\rm and} \qquad \lim_{n\rightarrow\infty}\sigma_n=\sigma
\end{equation}
(weak convergence in ${\rm Prob}(\C))$. Note that \eqref{20.7.15} and \eqref{20.7.16} are equivalent to
\begin{equation}\label{20.7.20}
\limsup_{n\rightarrow\infty}\int_\C\log|z-\lambda|\,\d\rho_n(z) \leq \int_\C\log|z-\lambda|\,\d\rho(z)
\end{equation}
and
\begin{equation}\label{20.7.21}
\limsup_{n\rightarrow\infty}\int_\C\log|z-\lambda|\,\d\sigma_n(z) \leq \int_\C\log|z-\lambda|\,\d\sigma(z),
\end{equation}
respectively. Let $p\in (0,1)$. By \cite[Theorem~2.19]{HS},
\begin{eqnarray*}
\int|z|^p\,\d\rho_n(z) &\leq& \|P_nT_nP_n\|_{L^p(P_n\CN P_n)}^p\\
&=&\frac{1}{\tau(P_n)}\|P_nT_nP_n\|_p^p\\
&\leq& \frac{1}{\tau(P_n)}\|T + \textstyle{\frac1n}xy^{-1}\|_p^p\\
&\leq& \frac{1}{\tau(P_n)}\big(\|T\|_p^p + \|xy^{-1}\|_p^p\big),
\end{eqnarray*}
and since $\rho$ has compact support, 
\[
\int|z|^p\,\d\rho(z)<\infty.
\]
Now, $\tau(P_n)\rightarrow\tau_\omega(P)$ as $n\rightarrow\infty$, and we can therefore choose a constant $C_p$ depending only on $p$ such that
\[
\int|z|^p\,\d\rho(z)<C_p \qquad {\rm and} \qquad \int|z|^p\,\d\rho_n(z)<C_p, \qquad n\in\N.
\]
Then since $|z-\lambda|^p\leq |z|^p+|\lambda|^p$, 
\[
\int|z-\lambda|^p\,\d\rho(z)\leq C_p(\lambda)  \qquad {\rm and} \qquad 
\int|z-\lambda|^p\,\d\rho_n(z)\leq C_p(\lambda),
\]
where $C_p(\lambda):=C_p+|\lambda|^p$. We will prove \eqref{20.7.20} and \eqref{20.7.21} in the case $\lambda=0$ only. For general $\lambda$ the proof is essentially the same. For $1<R<\infty$ define
\[
\varphi_R(z)=\left\{\begin{array}{ll}
0 & |z|\leq 1,\\
\log|z|, &1<|z|\leq R,\\
 \log R, &|z|> R.
 \end{array}\right.
\]
Then
\begin{equation}\label{20.7.22}
\log^+|z|=\varphi_R(z)+ \log^+\Big(\frac{|z|}{R}\Big).
\end{equation}
$\varphi_R$ is continuous and bounded. Hence, \eqref{20.7.19} implies that
\begin{equation}\label{20.7.23}
\lim_{n\rightarrow\infty}\int_\C\varphi_R(z)\,\d\rho_n(z)= \int_\C \varphi_R(z)\,\d\rho(z).
\end{equation}
Moreover, for $0<p<1$,
\begin{eqnarray*}
\int_\C\log^+\Big(\frac{|z|}{R}\Big)\,\d\rho(z)&\leq &\frac1p \int_\C\Big(\frac{|z|}{R}\Big)^p\,\d\rho(z)\\
&\leq & \frac{R^{-p}C_p}{p}.
\end{eqnarray*}
Similarly,
\[
\int_\C\log^+\Big(\frac{|z|}{R}\Big)\,\d\rho_n(z)\leq \frac{R^{-p}C_p}{p}.
\]
Then by \eqref{20.7.22} and \eqref{20.7.23},
\[
\limsup_{n\rightarrow\infty}\Big|\int_\C \log^+|z|\,\d\rho_n(z) -  \int_\C\log^+|z|\,\d\rho(z)\Big| \leq \frac{2R^{-p}C_p}{p}
\]
for all $R>1$, implying that
\begin{equation}\label{20.7.24}
\lim_{n\rightarrow\infty} \int_\C \log^+|z|\,\d\rho_n(z) = \int_\C\log^+|z|\,\d\rho(z).
\end{equation}
Now choose a sequence $(f_k)_{k=1}^\infty$ of compactly supported continuous functions on $\C$ such that $f_k\geq 0$ and $f_k\nearrow \log^-|z|$ as $k\rightarrow\infty$. Then as in the proof of Lemma~\ref{lemma20.7.13}, 
\begin{eqnarray*}
\liminf_{n\rightarrow\infty}\int_\C\log^-|z|\,\d\rho_n(z) & \geq & \lim_{n\rightarrow\infty}\int_\C f_k\,\d\rho_n\\
&=& \int_\C f_k\,\d\rho.
\end{eqnarray*}
Letting $k\rightarrow\infty$, we find that
\begin{equation}\label{20.7.25}
\liminf_{n\rightarrow\infty}\int_\C\log^-|z|\,\d\rho_n(z)\geq \int_\C \log^-|z|\,\d\rho(z).
\end{equation}
Since $\log|z|=\log^+|z|-\log^-|z|$, \eqref{20.7.20} (with $\lambda=0$) now follows from \eqref{20.7.24} and \eqref{20.7.25}. \eqref{20.7.21} follows from a similar proof. This finishes the proof of \eqref{20.7.15} and \eqref{20.7.16}. 

Now, assume that $\Delta(T-\lambda\unit)>0$. We will apply Lemma~\ref{limsup} with
\begin{eqnarray*}
a_n&=& \tau(P_n)\log\Delta_{P_n\CN P_n}(P_nT_nP_n-\lambda P_n),\\
b_n&=& \tau(P_n^\bot)\log\Delta_{P_n^\bot\CN P_n^\bot}(P_n^\bot T_nP_n^\bot-\lambda P_n^\bot),\\
a&=& \tau_\omega(P)\int_\C\log|z-\lambda|\,\d\rho(z),\\
b&=& \tau_\omega(P^\bot)\int_\C\log|z-\lambda|\,\d\sigma(z).
\end{eqnarray*}
Then by \eqref{20.7.15}, \eqref{20.7.16} and Proposition~\ref{measureofball},
\[
\limsup_{n\rightarrow\infty}a_n\leq a \qquad {\rm and} \qquad \limsup_{n\rightarrow\infty}b_n\leq b,
\]
and by Proposition~\ref{prop6.9},
\[
a_n+b_n=\log\Delta(T_n-\lambda\unit).
\]
Since $\mu_T=\tau_\omega(P)\rho + \tau_\omega(P^\bot)\sigma$, we also have
\[
a+b=\int_\C\log|z-\lambda|\,\d\mu_T(z) = \log\Delta(T-\lambda\unit).
\]
Hence,
\[
\lim_{n\rightarrow\infty}(a_n+b_n)=a+b>-\infty.
\]
In particular, $a,b\in\R$. Then by Lemma~\ref{limsup},
\[
\lim_{n\rightarrow\infty}a_n=a \qquad {\rm and} \qquad \lim_{n\rightarrow\infty}b_n=b.
\]
Since $\lim_{n\rightarrow\infty}\tau(P_n)=\tau_\omega(P)\in (0,1)$, \eqref{20.7.17} now follows. A similar argument shows that \eqref{20.7.18} holds. $\endproof$

\vspace{.2cm}

\begin{prop}\label{prop20.7.15} If $\mu_T(\partial B(0,1))=0$ and $0<\tau_\omega(P)<1$, then w.r.t. weak convergence in ${\rm Prob}(\C)$,
\begin{equation}\label{20.7.26}
\lim_{n\rightarrow\infty}\mu_{P_nT_nP_n}=\mu_{P\tilde TP},
\end{equation}
and
\begin{equation}\label{20.7.27}
\lim_{n\rightarrow\infty}\mu_{P_n^\bot T_nP_n^\bot}=\mu_{P^\bot \tilde TP^\bot}.
\end{equation}
Moreover,
\[
\supp(\mu_{P\tilde TP})\subseteq \overline{B(0,1)}
\]
and
\[
\supp(\mu_{P^\bot\tilde TP^\bot})\subseteq \C\setminus B(0,1).
\]
\end{prop}

\proof Combining \eqref{20.7.17}, \eqref{20.7.18}, and Proposition~\ref{lim}, we find that when $\lambda\in\C$ and $\Delta(T-\lambda\unit)>0$, then
\begin{equation}\label{20.7.28}
\log\Delta_{P\CN^\omega P}(P\tilde TP-\lambda P)= \int_\C \log|z-\lambda|\,\d\rho(z),
\end{equation}
and
\begin{equation}\label{20.7.29}
\log\Delta_{P^\bot\CN^\omega P^\bot}(P^\bot\tilde TP^\bot-\lambda P^\bot)= \int_\C \log|z-\lambda|\,\d\sigma(z).
\end{equation}
The map $\lambda\mapsto \log\Delta(T-\lambda\unit)$ is subharmonic, and hence $\Delta(T-\lambda\unit)>0$ for a.e. $\lambda\in\C$ w.r.t. Lebesgue measure. Both sides of \eqref{20.7.28} and \eqref{20.7.29} are subharmonic functions of $\lambda$, and since two subharmonic functions which agree almost everywhere are identical, \eqref{20.7.28} and \eqref{20.7.29} hold for all $\lambda\in\C$, showing that
\[
\rho=\mu_{P\tilde T P} \qquad {\rm and} \qquad \sigma=\mu_{P^\bot \tilde T P^\bot}.
\]
Then \eqref{20.7.26} and \eqref{20.7.27} follow from Lemma~\ref{lemma20.7.13}. Moreover, by the definition of $\rho$ and $\sigma$ in Lemma~\ref{lemma20.7.13}, it is clear that
\[
\supp(\rho)\subseteq \overline{B(0,1)} \qquad {\rm and} \qquad \supp(\sigma)\subseteq \C\setminus B(0,1). \qquad \endproof
\]

\vspace{.2cm}

\begin{prop}\label{generalballs} Let $T\in\CM$, and suppose that for some
  $r>0$, $\mu_T(\partial B(0,r))=0$. Then there is a $\tilde T$-invariant
  projection $P_r\in \CN^\omega$, such that $\supp(\mu_{P_r \tilde T P_r})\subseteq
  \overline{B(0,r)}$, $\supp(\mu_{P_r^\bot \tilde T P_r^\bot})\subseteq \C\setminus
  B(0,r)$, and
  \[
  \tau_\omega(P_r)=\mu_T(\overline{B(0,r)}).
  \]
\end{prop}

\proof If $\mu_T(B(0,r))=1$, take $P_r=\unit$ and define $\mu_{P_r^\bot \tilde TP_r^\bot}$ to be 0. If $\mu_T(B(0,r))=0$, take $P_r=0$ and define $\mu_{P_r \tilde TP_r}$ to be 0. If $0<\mu_T(B(0,r))<1$, then $0<\mu_{{\frac1r} T}(B(0,1))<1$, and we can take $P_r\in\CN^\omega$ to be the ${\frac1r} \tilde T$--invariant projection found in the above with
\[
\supp(\mu_{P_r{\frac1r} \tilde T P_r})\subseteq \overline{B(0,1)} \qquad {\rm and} \qquad
\supp(\mu_{P_r^\bot {\frac1r} \tilde T P_r^\bot})\subseteq \C\setminus B(0,1)
\]
and with
\[
\tau_\omega(P_r)=\mu_{{\frac1r} T}(B(0,1))=\mu_{{\frac1r} T}(\overline{B(0,1)}).
\]
Clearly, $P_r$ satisfies the conditions listed in Proposition~\ref{generalballs}. $\endproof$

\subsection*{The $T$-invariant subspaces $(E(T,r))_{r>0}$ and $(F(T,r))_{r>0}$.}

Recall from Section~\ref{sec6} that for each $r>0$ one can define
$T$-hyperinvariant subspaces $E(T,r)$ and $F(T,r)$, for which the projections $P_{E(T,r)}$ and $P_{F(T,r)}$ are independent of
the representation of $\CM$ on a Hilbert space. In particular, we may
regard $E(T,r)$ and $F(T,r)$ as subspaces of the Hilbert space $\CH = L^2(\CN^\omega, \tau_\omega)$, on which $\CN^\omega$ acts. 

\begin{lemma}\label{Pr(H)=E(T,r)} Let $T\in\CM$ and let $r>0$. If $\mu_T(\partial B(0,r))=0$,
  then 
\begin{itemize}
\item[(i)] $\supp(\mu_{T|_{E(T,r)}})\subseteq \overline{B(0,r)}$,
\item[(ii)]   $\supp(\mu_{T|_{F(T,r)}})\subseteq \C\setminus B(0,r)$,
\item[(iii)] $\tau(P_{E(T,r)}) = \mu_T(\overline{B(0,r)})$,
\item[(iv)] $\tau(P_{F(T,r)})=\mu_T(\C\setminus B(0,r))$,
\item[(v)] $E(T,r)=F(T\cc,r)^\bot$,
\item[(vi)] $F(T,r)=E(T\cc,r)^\bot$.
\end{itemize}
\end{lemma}

\proof In the proof of (i)--(vi) we will consider $\CM$ as a von Neumann algebra acting on $L^2(\CN^\omega,\tau_\omega)$ and thus identify $T\in\CM$ with $\tilde{T}\in\CN^\omega$. Let $P_r$ denote the projection from Proposition~\ref{generalballs}. Then
\[
\supp(\mu_{P_rTP_r})\subseteq \overline{B(0,r)},
\]
\[
\supp(\mu_{P_r^\bot T P_r^\bot})\subseteq \C\setminus B(0,r),
\]
and
\[
\tau_\omega(P_r)=\mu_T(B(0,r))=\mu_T(\overline{B(0,r)}).
\]
For $A\subseteq \C$, let $A\cc= \{\overline{z}\,|\,z\in A\}$. Then
\[
\supp(\mu_{P_r^\bot T\cc P_r^\bot})\subseteq (\C\setminus B(0,r))\cc =\C\setminus B(0,r).
\]
Since $P_r(\CH)$ is $T$--invariant and $P_r(\CH)^\bot$ is $T\cc$--invariant, it follows from Corollary~\ref{Cor3.5} that 
\begin{equation}\label{25/11b}
P_r(\CH)\subseteq E(T,r) \qquad {\rm and}\qquad P_r(\CH)^\bot \subseteq F(T\cc,r).
\end{equation}
Hence,
\begin{equation}\label{25/11c}
\tau(P_{E(T,r)}) \geq \tau_\omega(P_r)= \mu_T(\overline{B(0,r)}),
\end{equation}
\begin{equation}\label{25/11d}
\tau(P_{F(T\cc,r)}) \geq 1- \tau_\omega(P_r)= \mu_T(\C\setminus B(0,r)).
\end{equation}
Now choose a sequence $(s_n)_{n=1}^\infty$ from $]r,\infty[$ such that
$\mu_T(\partial B(0,s_n))=0$ for all $n\in\N$ and $s_n\searrow r$ as
$n\rightarrow \infty$. Applying \eqref{25/11d} to $s_n$, we find that
\[
\tau(F(T\cc,s_n))\geq \mu_T(\C\setminus B(0,s_n)), \qquad n\in\N.
\]
Moreover, $E(T,r)\bot F(T\cc,s_n)$ (cf. Lemma~\ref{perpendicular}). Hence,
\[
\tau(P_{E(T,r)})\leq 1-\tau(P_{F(T\cc,s_n)})\leq \mu_T(\C\setminus B(0,s_n)).
\]
Letting $n\rightarrow \infty$, we see that
\[
\tau(P_{E(T,r)})\leq \mu_T(\C\setminus B(0,r)).
\]
Thus the inequality \eqref{25/11c} is an equality. Next, choose a sequence $(t_n)_{n=1}^\infty$ from $(0,r)$ such that $t_n\nearrow r$ as $n\rightarrow\infty$ and with $\mu_T(\partial B(0,t_n))=0$ for all $n$. Arguing as above, we find that
\[
\tau(P_{F(T\cc,r)})\leq 1-\tau(P_{E(T,t_n)})\leq \mu_T(\C\setminus\overline{B(0,t_n)}).
\]
Letting $n\rightarrow \infty$, we obtain
\[
\tau(P_{F(T\cc,r)})\leq \mu_T(\C\setminus B(0,r)).
\] 
Hence, \eqref{25/11d} is an equality too. With \eqref{25/11c} and \eqref{25/11d} being equalities it follows from \eqref{25/11b} that
\[
P_r(\CH)=E(T,r) \qquad {\rm and} \qquad P_r(\CH)^\bot =F(T\cc,r).
\]
Altogether, we have proven (i), (ii), and (v) for $T$ and (ii), (iv), and (vi) with $T$ replaced by $T\cc$. Since $\mu_{T\cc}(\partial B(0,r))=\mu_T(\partial B(0,r))=0$, and since $T$ was arbitrary, this finishes the proof. $\endproof$
  
\vspace{.2cm}

\begin{remark}\label{remark on supports} If $T\in\CM$, and if $E_1$ and
  $E_2$ are closed, $T$-invariant subspaces affiliated with $\CM$ with
  $E_1\subseteq E_2\neq \{0\}$, then $\supp(\mu_{T|_{E_1}})\subseteq
  \supp(\mu_{T|_{E_2}})$. Indeed, $E_1$ is invariant under $T|_{E_2}$, and
  if $P_i$ denotes the projection onto $E_i$, then $E_1$ is affiliated with
  $P_2\CM P_2$, whence
  \begin{eqnarray*}
  \mu_{T|_{E_2}}&=& \tau_2(P_1)\cdot\mu_{(T|_{E_2})|_{E_1}}) +
  \tau_2(P_2-P_1)\cdot \mu_{(P_2-P_1)T|_{E_2}(P_2-P_1)}\\
  & = & \tau_2(P_1)\cdot\mu_{T|_{E_1}} +  \tau_2(P_2-P_1)\cdot
  \mu_{(P_2-P_1)T|_{E_2}(P_2-P_1)},
  \end{eqnarray*}
  where $\tau_2= \frac{1}{\tau(P_2)}\tau$ denotes the trace on $P_2\CM
  P_2$. It follows that $\supp(\mu_{T|_{E_1}})\subseteq
  \supp(\mu_{T|_{E_2}})$. 

\end{remark}

\vspace{.2cm}

\begin{lemma}\label{general r} For $T\in\CM$ and for arbitrary $r>0$, (i), (ii), (iii), and (iv) of Lemma~\ref{Pr(H)=E(T,r)} hold. Moreover,
\begin{itemize}
\item[(v')] $E(T,r)=\bigcap_{s>r}F(T\cc,s)^\bot \supseteq F(T\cc,r)^\bot$,
\item[(vi')] $F(T,r)=\bigcap_{s<r} E(T\cc,s)^\bot \supseteq E(T\cc,r)^\bot$.
\end{itemize}
\end{lemma}

\proof Choose sequences $(s_n)_{n=1}^\infty$ from $(r,\infty)$ and $(t_n)_{n=1}^\infty$ from $(0,r)$ such that $s_n\searrow r$, $t_n\nearrow r$, and with $\mu_T(\partial B(0,s_n))=\mu_T(\partial B(0,t_n))=0$ for all $n$. According to Remark~\ref{remark on supports} and Lemma~\ref{Pr(H)=E(T,r)}~(i) and (ii),
\begin{eqnarray*}
\supp(\mu_{T|_{E(T,r)}}) &=& \bigcap_{n=1}^\infty \supp(\mu_{T|_{E(T,s_n)}})\\
&\subseteq & \bigcap_{n=1}^\infty \overline{B(0,s_n)}\\
&=& \overline{B(0,r)}.
\end{eqnarray*}
Similarly, with the aid of the sequence $(t_n)_{n=1}^\infty$ we find that
\[
\supp(\mu_{T|_{F(T,r)}})\subseteq \C\setminus B(0,r).
\]
According to Lemma~\ref{Pr(H)=E(T,r)}~(iii) and (iv),
\begin{eqnarray*}
\tau(P_{E(T,s_n)}) &=& \mu_T(\overline{B(0,s_n)}),\\
\tau(P_{F(T,t_n)}) &=& \mu_T(\C\setminus B(0,t_n)).
\end{eqnarray*}
Letting $n\rightarrow\infty$ and applying Lemma~\ref{Properties}~(b) we then have
\begin{eqnarray*}
\tau(P_{E(T,r)}) &=& \mu_T(\overline{B(0,r)}),\\
\tau(P_{F(T,r)}) &=& \mu_T(\C\setminus B(0,0)).
\end{eqnarray*}
By application of Lemma~\ref{Pr(H)=E(T,r)}~(v),
\[
E(T,r)=\bigcap_{n=1}^\infty E(T,s_n)= \bigcap_{n=1}^\infty F(T\cc,s_n)^\bot \supseteq \bigcap_{s>r} F(T\cc,s)^\bot,
\]
and then by Lemma~\ref{perpendicular},
\[
E(T,r)=\bigcap_{s>r} F(T\cc,s)^\bot.
\]
Since $F(T\cc,s)\subseteq F(T\cc,r)$, $s>r$ (cf. Lemma~\ref{Properties}~(b)), this proves (v'). (vi') is proven in a similar way using the sequence $(t_n)_{n=1}^\infty$. $\endproof$

\vspace{.2cm}

\begin{cor}\label{maximality of E(T,r)} For every $T\in\CM$, every $\lambda\in\C$ and every $r>0$ one
  has:
  \begin{itemize}
    \item[(i)] $E:=E(T-\lambda\unit,r)$ is the largest, closed $T$-invariant
    subspace affiliated with $\CM$, such that $\supp(\mu_{T|_E})\subseteq
    \overline{B(\lambda,r)}$.
    \item[(ii)] $F:=F(T-\lambda\unit,r)$ is the largest, closed $T$-invariant
    subspace affiliated with $\CM$, such that $\supp(\mu_{T|_F})\subseteq
    \C\setminus B(\lambda,r)$.
  \end{itemize}
\end{cor}

\proof It suffices to consider the case $\lambda=0$, because
$\mu_{T-\lambda\unit|_E}$ and $\mu_{T-\lambda\unit|_F}$ are the
push-forward measures of $\mu_{T|_E}$ and $\mu_{T|_F}$, respectively, under
the map $z\mapsto z-\lambda$. Moreover, because of Corollary~\ref{Cor3.5}
we only have to prove that $\supp(\mu_{T|_E})\subseteq \overline{B(0,r)}$
and $\supp(\mu_{T|_F})\subseteq \C\setminus B(0,r)$, and these properties follow from Lemma~\ref{general r}. $\endproof$

\vspace{.2cm}

\subsection*{Spectral subspaces corresponding to closed sets.}

\begin{prop}\label{closedsets} For every $T\in\CM$ and every closed set $F\subseteq \C$ there
  is a largest, closed $T$-invariant subspace $\CK=\CK_T(F)$ affiliated with $\CM$,
  such that $\supp(\mu_{T|_\CK})\subseteq F$. Moreover, $\CK$ is
  hyperinvariant for $T$.
\end{prop}

\proof We may write $\C\setminus F$ as a union of countably many open balls
$(B(\lambda_k, r_k))_{k=1}^\infty$:
\[
\C\setminus F = \bigcup_{k=1}^\infty B(\lambda_k,r_k).
\]
With
\begin{equation}\label{defineK}
  \CK:= \bigcap_{k=1}^\infty F(T-\lambda_k\unit,r_k),
\end{equation}
$\CK$ is hyperinvariant for $T$, and according to Remark~\ref{remark on
  supports},
\[
\supp(\mu_{T|_\CK})\subseteq \C\setminus B(\lambda_k,r_k)
\]
for every $k\in\N$. Hence, $\supp(\mu_{T|_\CK})\subseteq F$.

In order to prove that $\CK$ is the largest closed subspace of $\CH$ having this property, assume that $\CK'$ is a closed
$T$-invariant subspace affiliated with $\CM$, such that
$\supp(\mu_{T|_{\CK'}})\subseteq F$. Then, for every $k\in\N$,
$\supp(\mu_{T|_{\CK'}})\subseteq\C\setminus B(\lambda_k,r_k)$. Therefore,
by Corollary~\ref{Cor3.5}, $\CK'\subseteq F(T-\lambda_k\unit,r_k)$, and it
follows that $\CK'\subseteq \CK$. $\endproof$

\vspace{.2cm}

\begin{definition} For $T\in\CM$ and $F$ a closed subset of $\C$ we denote
  by $P_T(F)\in W\cc(T)$ the projection onto the subspace $\CK_T(F)$ found
  in Proposition~\ref{closedsets}.
\end{definition}

The following proposition is a trivial consequence of
Proposition~\ref{closedsets} and Remark~\ref{remark on supports}:

\begin{prop}\label{26/10a} For every $T\in\CM$ one has that
  \begin{itemize}
    \item[(i)] $P_T(\emptyset)=0$ and $P_T(\C)=\unit$
    \item[(ii)] If $F_1$ and $F_2$ are closed subsets of $\C$ with
    $F_1\subseteq F_2$, then $P_T(F_1)\leq P_T(F_2)$
    \item[(iii)] If $(F_i)_{i\in I}$ is a family of closed subsets of $\C$,
    then $P_T\big(\bigcap_{i\in I}F_i\big)= \bigwedge_{i\in I}P_T(F_i)$
  \end{itemize}
\end{prop}
    
\vspace{.2cm}

\begin{lemma}\label{Oct5} For every $T\in\CM$ and every closed subset $F$ of $\C$,
  \begin{equation}\label{mindreend}
    \tau(P_T(F))\leq \mu_T(F).
  \end{equation}
\end{lemma}

\proof This is an easy consequence of the fact that
\[
\mu_T= \tau(P_T(F))\cdot \mu_{P_T(F)TP_T(F)}+\tau(P_T(F)^\bot)\cdot
\mu_{P_T(F)^\bot T P_T(F)^\bot}
\]
with $ \mu_{P_T(F)TP_T(F)}(F)=1$. $\endproof$

\vspace{.2cm}

\begin{prop}\label{boundary, F} Let $T\in\CM$. Then for every closed subset $F$ of $\C$ with
  $\mu_T(\partial F)=0$,
\begin{itemize}
\item[(i)] $\tau(P_T(F))=\mu_T(F)$,
\item[(ii)] $\CK_T(F)=\CK_{T\cc}(\overline{\C\setminus F\cc})^\bot$.
\end{itemize}
\end{prop}

\proof As in the proof of Proposition~\ref{closedsets}, write $\C\setminus F$ as a union of countably many open balls
$(B(\lambda_k, r_k))_{k=1}^\infty$:
\[
\C\setminus F = \bigcup_{k=1}^\infty B(\lambda_k,r_k).
\]
Then $\CK_T(F)$ is given by
\[
\CK_T(F)=\bigcap_{k=1}^\infty F(T-\lambda_k\unit,r_k).
\]
By Proposition~\ref{26/10a}~(ii),
\[
\CK_{T\cc}(\overline{\C\setminus F\cc}) \supseteq \CK_{T\cc}(\overline{B(\lambda_k,r_k)})=E(T\cc-\lambda_k\unit,r_k),
\]
and then by Lemma~\ref{Pr(H)=E(T,r)}~(v'),
\begin{equation*}
\CK_{T\cc}(\overline{\C\setminus F\cc}) \supseteq F(T-\lambda_k\unit,r_k)^\bot.
\end{equation*}
It follows that
\begin{equation}\label{26/10b}
\CK_T(F)=\bigcap_{k=1}^\infty F(T-\lambda_k\unit,r_k)\supseteq \CK_{T\cc}(\overline{\C\setminus F\cc})^\bot.
\end{equation}
Then for the corresponding projections we have:
\begin{equation}\label{26/10c}
\tau(P_T(F))\geq 1- \tau(P_{T\cc}(\overline{\C\setminus F\cc}).
\end{equation}
According to Lemma~\ref{Oct5},
\begin{equation}\label{26/10d}
\tau(P_T(F))\leq \mu_T(F),
\end{equation}
and since $\mu_T(\partial F)=0$, 
\begin{equation}\label{26/10f}
1-\tau(P_{T\cc}(\overline{\C\setminus F\cc})) \geq 1- \mu_{T\cc}(\overline{\C\setminus F\cc}) = 1- \mu_{T}(\C\setminus F)= \mu_T(F).
\end{equation}
We deduce from \eqref{26/10c}, \eqref{26/10d}, and \eqref{26/10f} that 
\[
\tau(P_T(F)) = \mu_T(F) = \tau(\unit-P_{T\cc}(\overline{\C\setminus F\cc})).
\]
Then by \eqref{26/10b},
\[
\CK_T(F)=\CK_{T\cc}(\overline{\C\setminus F\cc})^\bot. \qquad \endproof
\]

\vspace{.2cm}

\begin{prop}\label{ligmed} For every $T\in\CM$ and every closed subset $F$ of $\C$,
  $\tau(P_T(F))=\mu_T(F)$, and hence $\mu_{P_T(F)^\bot
  TP_T(F)^\bot}$ is concentrated on $ \C\setminus F$. 
\end{prop}

\proof For $t>0$ define
\[
F_t=\Big\{z\in\C\,|\,{\rm dist}(z,F)\leq \frac1t\Big\}.
\]
The map $t\mapsto \mu_T(F_t)\in [0,1]$ is decreasing, and therefore it has
at most countably many points of discontinuity. This entails that
$\mu_T(\partial F_t)=0$ for all but countably many $t>0$. Choose
$t_1>t_2>t_3>\ldots$, such that $t_n \searrow 0$ as $n\rightarrow \infty$
and $\mu_T(\partial F_{t_n})=0$ for all $n\in\N$. Since $F_{t_n}\searrow
\bigcap_{m=1}^\infty F_{t_m}=F$, we have that $P_T(F_{t_n})\searrow
P_T(F)$, and hence
\begin{eqnarray*}
  \tau(P_T(F))&=& \lim_{n\rightarrow\infty} \tau(P_T(F_{t_n}))\\
  &=& \lim_{n\rightarrow\infty}\mu_T(F_{t_n})\\
  & = & \mu_T(F),
\end{eqnarray*}
as claimed. Finally, since
\[
\mu_T = \tau(P_T(F))\cdot \mu_{P_T(F)TP_T(F)} + \tau(P_T(F)^\bot)\cdot
\mu_{P_T(F)^\bot T P_T(F)^\bot},
\]
where $ \mu_{P_T(F)TP_T(F)}(F)=1$, we conclude that
\[
\tau(P_T(F)^\bot)\cdot
\mu_{P_T(F)^\bot T P_T(F)^\bot}(F)=0.
\]
Hence, if $P_T(F)\neq \unit$, then $\mu_{P_T(F)^\bot T
  P_T(F)^\bot}(F)=0$. If $P_T(F)=\unit$, then, by definition,
  $\mu_{P_T(F)^\bot T P_T(F)^\bot}=0$, and this measure is trivially
  concentrated on  $ \C\setminus F$. $\endproof$

\vspace{.2cm}

\begin{lemma}\label{intersect-inv} Let $T\in\CM$, and let $P\in\CM$ be a
  $T$-invariant projection. Then for every closed subset $F$ of $\C$, 
  \begin{equation}
    \CK_{T|_{P(\CH)}}(F)=\CK_T(F)\cap P(\CH).
  \end{equation}
\end{lemma}

\proof Let $Q\in P\CM P$ denote the projection onto $\CK_{T|_{P(\CH)}}(F)$, and
let $R = P_T(F)\wedge P$. We will prove that $Q\leq R$ and $R\leq Q$.

Clearly, $Q\leq P$. In order to see that $Q\leq P_T(F)$, recall that
$P_T(F)$ is the largest projection with the properties
\begin{itemize}
  \item[(i)] $P_T(F)TP_T(F)=TP_T(F)$,
  \item[(ii)] $\mu_{P_T(F)TP_T(F)}$ (computed relative to $P_T(F)\CM P_T(F)$) is concentrated on $F$.
\end{itemize}

Since
\begin{equation}
  QTQ = QTPQ = TPQ = TQ,
\end{equation}
and $\mu_{QTQ}= \mu_{QTPQ}$ (computed relative to $Q\CM Q$) is concentrated
on $F$, we get that $Q\leq P_T(F)$, and hence $Q\leq R$.

Similarly, to prove that $R\leq Q$, prove that
\begin{itemize}
  \item[(i')] $RTPR = TPR$, i.e. $RTR = TR$,
  \item[(ii')] $\mu_{RTPR}= \mu_{RTR}$  (computed relative to $R\CM R$) is concentrated on $F$.
\end{itemize}

Note that if $P_T(F)=0$, then $R\leq Q$, so we may assume that $P_T(F)\neq 0$.
(i') holds, because $R(\CH)= P(\CH)\cap P_T(F)(\CH)$ is $T$-invariant when
$P(\CH)$ and $P_T(F)(\CH)$ are $T$-invariant. In order to prove (ii'), at
first note that $R(\CH)$ is $TP_T(F)$-invariant. Hence
\begin{equation}
  \mu_{TP_T(F)}= \tau_1(R)\cdot \mu_{RTR} + \tau_1(R^\bot)\cdot
  \mu_{R^\bot TR^\bot},
\end{equation}
where $\tau_1=\frac{1}{\tau(P_T(F))}\cdot\tau|_{P_T(F)\CM P_T(F)}$. It follows that
\begin{equation}
  \tau_1(R)\cdot\mu_{RTR}(F^c)\leq \mu_{TP_T(F)}(F^c)=0,
\end{equation}
and thus, if $R\neq 0$, then $\mu_{RTR}(F^c)=0$, and (ii') holds. If $R=0$,
then $R\leq Q$ is trivially fulfilled. $\endproof$

\subsection*{The general case.}

{\it Proof of Theorem~\ref{Borelsets}.} Define
\begin{equation}
  P_T(B):= \bigvee_{K\;{\rm compact},\; K\subseteq B} P_T(K)
\end{equation}
and
\begin{equation}
  \CK_T(B)=P_T(B)(\CH).
\end{equation}
Then $\CK_T(B)$ is $T$-hyperinvariant. Since $\mu_{T|_{\CK_T(B)}}$ is a
regular measure (cf. \cite[Theorem~7.8]{Fo}),
\[
\mu_{T|_{\CK_T(B)}}(B)=\sup\{\mu_{T|_{\CK_T(B)}}(K)\,|\,K\;{\rm compact},
\; K\subseteq B\}.
\]
For every compact set $K\subseteq B$, $\CK_T(K)\subseteq \CK_T(B)$ is
$T|_{\CK_T(B)}$-invariant, so with $P=P_T(B)$,
\[
\mu_{T|_{\CK_T(B)}}= \tau_{P\CM P}(P_T(K))\cdot \mu_{T|_{\CK_T(K)}} +
\tau_{P\CM P}(P-P_T(K))\cdot \mu_{(P-P_T(K))T(P-P_T(K))}.
\]
Therefore, by Proposition~\ref{ligmed},
\[
\mu_{T|_{\CK_T(B)}}(K)=\frac{1}{\tau(P)}\cdot 
\tau(P_T(K))\cdot \mu_{T|_{\CK_T(K)}}(K)+0=\frac{1}{\tau(P)}\cdot
\tau(P_T(K)).
\]
$(P_T(K))_{K\;{\rm compact},\; K\subseteq B}$ is an increasing net of
projections with SO-limit $P_T(B)$. Therefore,
\[
\mu_{T|_{\CK_T(B)}}(B)=\sup\Big\{\frac{1}{\tau(P)}\cdot
\tau(P_T(K))\,|\,K\;{\rm compact},
\; K\subseteq B\Big\}= 1.
\]
This shows that $\mu_{T|_{\CK_T(B)}}$ is concentrated on $B$. Moreover,
by similar arguments,
\begin{eqnarray*}
\mu_T(B) &=& \sup\{\mu_{T}(K)\,|\,K\;{\rm compact},
\; K\subseteq B\}\\
&= & \sup\{\tau(P_T(K))\,|\,K\;{\rm compact},
\; K\subseteq B\}\\
&=& \tau(P_T(B)),
\end{eqnarray*}
proving that (i) of Theorem~\ref{Borelsets} holds. (ii) then follows as in the foregoing proof.

Finally, suppose that $Q\in\CM$ is a $T$-invariant projection and that
$\mu_{T|_{Q(\CH)}}$ is concentrated on $B$. Then by Lemma~\ref{intersect-inv} and Proposition~\ref{ligmed},
\begin{eqnarray*}
  \tau_{Q\CM Q}(P\wedge Q)&=& \sup\{\tau_{Q\CM Q}(P_T(K)\wedge Q)\,|\,K\;{\rm compact},
\; K\subseteq B\}\\
&=& \sup\{\tau_{Q\CM Q}(P_{T|_{Q(\CH)}}(K))\,|\,K\;{\rm compact},
\; K\subseteq B\}\\
&=& \sup\{\mu_{T|_{Q(\CH)}}(K)\,|\,K\;{\rm compact},
\; K\subseteq B\}\\
&=& \mu_{T|_{Q(\CH)}}(B)\\
&=& 1.
\end{eqnarray*}
Hence, $P\wedge Q= Q$, and we get that $Q\leq P$. $\endproof$

\vspace{.2cm}

\begin{cor} For every Borelset $B\subseteq\C$, 
\[
P_{T\cc}(\C\setminus B\cc)=\unit-P_T(B),
\]
where $B\cc =\{\overline{z}\,|\,z\in B\}$.
\end{cor}

\proof Let $\CK=\CK_T(B)$ and $P=P_T(B)$ be as in Theorem~\ref{Borelsets}. Then by Theorem~\ref{Borelsets}~(ii), $\CK^\bot$ is a closed $T\cc$--invariant subspace, and the Brown measure of $P^\bot T\cc P^\bot =(P^\bot T P^\bot)\cc$ is concentrated on $B\cc$. Hence,
\begin{equation}\label{Oct5a}
P_{T\cc}(\C\setminus B\cc)\geq P^\bot = P_T(B)^\bot.
\end{equation}
But
\begin{eqnarray*}
\tau(P_{T\cc}(\C\setminus B\cc)) &=& \mu_{T\cc}(\C\setminus B\cc)\\
&=& \mu_T(\C\setminus B)\\
&=& 1-\mu_T(B)\\
&=& \tau(P_T(B)^\bot).
\end{eqnarray*}
Hence, equality must hold in \eqref{Oct5a}. $\endproof$  

\vspace{.2cm}

\begin{cor} Let $T\in\CM$, let $B\subseteq \C$ be a Borelset, and let $\CK$ be a closed, $T$--invariant subspace of $\CH$ which is affiliated with $\CM$. Then the following two conditions are equivalent:
\begin{itemize}
\item[(i)] $\CK=\CK_T(B)$,
\item[(ii)] $\mu_{P_\CK TP_\CK}$ is concentrated on $B$ and $\mu_{P_\CK^\bot TP_\CK^\bot}$ is concentrated on $\C\setminus B$.
\end{itemize}
\end{cor}

\proof That (i) implies (ii) is a consequence of Theorem~\ref{Borelsets}. Now, suppose that (ii) holds. Then $\CK\subseteq \CK_T(B)$. Moreover, $\CK^\bot$ is $T\cc$--invariant, and $\mu_{P_\CK^\bot T\cc P_\CK^\bot}$ is concentrated on $(\C\setminus B)\cc = \C\setminus B\cc$. Therefore, $\CK^\bot\subseteq \CK_{T\cc}(\C\setminus B\cc)=\CK_T(B)^\bot$. Hence, $\CK=\CK_T(B)$. $\endproof$ 

%% file: stronglimit.tex
\section{Realizing $P_{E(T,r)}$ and $P_{F(T,r)}$ as spectral projections}

Recall from Section~\ref{sec6} that for every $T\in\CM$ and every $r>0$ we defined $T$-hyperinvariant subspaces $E(T,r)$ and $F(T,r)$. The aim of the present section is to show that the
corresponding projections, $P_{E(T,r)}$ and $P_{F(T,r)}$, have the following property:

\begin{thm}\label{existence of A and B} For every $T\in\CM$ we have:
  \begin{itemize}
    \item[(a)] There is a unique operator $A\in\CM^+$, such that for every
    $r>0$,
    \begin{equation}\label{eq8-2}
      P_{E(T,r)}=1_{[0,r]}(A).
    \end{equation}
    Moreover,
    \begin{equation}
      A= {\rm SO}-\lim_{n\rightarrow\infty}((T\cc)^n T^n)^{\frac{1}{2n}}.
    \end{equation}
    \item[(b)] There is a unique operator $B\in\CM^+$, such that for every
    $r>0$,
    \begin{equation}
      P_{F(T,r)}=1_{[r,\infty[}(B).
    \end{equation}
    Moreover,
    \begin{equation}
      B= {\rm SO}-\lim_{n\rightarrow\infty}(T^n (T\cc)^n)^{\frac{1}{2n}}.
    \end{equation}
  \end{itemize}
\end{thm}

\vspace{.2cm}

In the proof of Theorem~\ref{existence of A and B} we shall need the
following two lemmas. The first one of them is elementary, and we omit the
proof of it.

\begin{lemma}\label{lemma a} Let $P, P_1, P_2, \ldots$ be projections in $B(\CH)$. Then
  the following are equivalent:
  \begin{itemize}
    \item[(i)] $P_n\rightarrow P$ in the strong operator topology,
     \item[(ii)] $\|P_n\xi -\xi\|\rightarrow 0$ for every $\xi\in P(\CH)$
     and $\|P_n^\bot\eta -\eta\|\rightarrow 0$ for every $\eta\in
     P(\CH)^\bot$.
  \end{itemize}
\end{lemma}

\vspace{.2cm}

\begin{lemma}\label{B} Let $A, A_1, A_2, \ldots$ be operators from
  $B(\CH)^+$, and assume that
  \[
  M:= \max\{\|A\|, \, \sup_{n\in\N}\|A_n\|\}<\infty.
  \]
  If
  \begin{equation}\label{eq8-1}
  1_{[0,r]}(A_n)\rightarrow 1_{[0,r]}(A)
  \end{equation}
  in the strong operator topology for all but countably many $r\in [0,M]$,
  then $A_n\rightarrow A$ in the strong operator topology.
\end{lemma}

\proof We prove this, approximating $A_n$ and $A$ with linear combinations
of projections of the form $1_{[0,r]}(A_n)$ and $ 1_{[0,r]}(A)$,
respectively. Choose $\alpha>0$ such that \eqref{eq8-1} holds for all $r\in
\alpha \Q_+ \cap [0,M]$. For each $k\in\N$ define
\[
B_k = \frac{\alpha}{k}\sum_{k=1}^{[kM\alpha^{-1}]} 1_{]\frac{n\alpha}{k},
  \infty[}(A),
\]
\[
B_{n,k} = \frac{\alpha}{k}\sum_{k=1}^{[kM\alpha^{-1}]} 1_{]\frac{n\alpha}{k},
  \infty[}(A_n).
\]
Then, by the Borel functional calculus for normal operators,
\[
\|A-B_k\|\leq \frac{\alpha}{k},
\]
and
\[
\|A_n-B_{n,k}\|\leq \frac{\alpha}{k}.
\]
By assumption, for fixed $k\in\N$, $B_k= {\rm
  SO}-\lim_{n\rightarrow\infty}B_{n,k}$, and it follows that
  $A_n\rightarrow A$ in the strong operator topology. $\endproof$

{\it Proof of Theorem~\ref{existence of A and B}.} According to
Lemma~\ref{Properties}~(b), if we define $E(T,0):= \bigcap_{r>0}E(T,r)$,
then $r\mapsto P_{E(T,r)}$ is increasing and SO-continuous from the
right. Moreover, according to Lemma~\ref{Lemma3.4}, $P_{E(T,r)}=\unit$ for
every $r\geq r'(T)$. \cite[Theorem~5.2.4]{KR} then implies that there is
one and only one operator $A\in\CM^+$, such that \eqref{eq8-2} holds. Moreover,
$\|A\|\leq r'(T)\leq \|T\|$.

Now, take $r>0$ such that $\mu_T(\partial B(0,r))=0$. Then
\begin{itemize}
  \item[(i)] $E(T,r)= F(T\cc, r)^\bot$.
\end{itemize}
Moreover, we claim that
\begin{itemize}
  \item[(ii)] $\overline{\bigcup_{0<s<r}E(T,s)}= E(T,r)$,
  \item[(iii)] $\overline{\bigcup_{r<t<\infty}F(T\cc,t)}= F(T\cc,r)$.
\end{itemize}
The one inclusion $\subseteq$ in (ii) is obvious. On the other hand,
\begin{eqnarray*}
\tau(P_{\overline{\bigcup_{0<s<r}E(T,s)}})&=& \lim_{s\rightarrow
  r-}\tau(P_{E(T,s)})\\
&= & \lim_{s\rightarrow r-}\mu_T(\overline{B(0,s)})\\
&=& \mu_T(B(0,r))\\
&=& \tau(P_{E(T,r)}).
\end{eqnarray*}
Hence, ``$=$'' holds in (ii). Similarly, ``$\subseteq$'' holds in (iii). On
the other hand,
\begin{eqnarray*}
\tau(P_{\overline{\bigcup_{r<t<\infty}F(T\cc,t)}})&=& \lim_{t\rightarrow
  r+}\tau(P_{F(T\cc,t)})\\
&=& \lim_{t\rightarrow
  r+}\tau(\unit - P_{E(T,t)})\\
&= & \lim_{t\rightarrow r+}\mu_T(\C\setminus {B(0,t)})\\
&=& \mu_T(\C\setminus B(0,r))\\
&=& \tau(\unit - P_{E(T,r)})\\
&=& \tau(P_{F(T\cc, r)}).
\end{eqnarray*}
Thus, (iii) holds.

Now, let $\xi\in E(T,r)$ and $\eta\in E(T,r)^\bot = F(T\cc,r)$ with
$\|\xi\|=\|\eta\|=1$. Let $\eps>0$. According to (ii) and (iii), we may take $s\in
(0,r)$, $t\in (r,\infty)$, $\xi'\in E(T,s)$ and $\eta'\in F(T\cc,t)$ with
$\|\xi'\|= \|\eta'\|=1$, $\|\xi-\xi'\|<\frac{\eps}{2}$ and
$\|\eta-\eta'\|<\frac{\eps}{2}$ . Next choose $\xi_n'$ and $\eta_n'$ in
$\CH$ such that
\[
\lim_{n\rightarrow\infty}\|\xi_n'-\xi'\|=0 \quad {\rm and} \quad
\limsup_{n\rightarrow\infty}\|T^n\xi_n'\|^{\frac 1n}\leq s,
\]
and
\[
\lim_{n\rightarrow\infty}\|(T\cc)^n\eta_n'-\eta'\|=0 \quad {\rm and} \quad
\limsup_{n\rightarrow\infty}\|\xi_n'\|^{\frac 1n}\leq \frac 1t.
\]
We can, without loss of generality, assume that
$\|\xi_n'\|=\|(T\cc)^n\eta_n'\|=1$ for all $n\in\N$. Then let $\rho_n$
(respectively $\sigma_n$) denote the distribution of $(T\cc)^nT^n$ w.r.t
the vector state on $\CM$ induced by $\xi_n'$ (respectively
$(T\cc)^n\eta_n'$). Arguing as in the proof of Lemma~\ref{Lemma3.4}, we get that
\[
\rho_n(]r^{2n},\infty[)\rightarrow 0 \quad {\rm as} \quad
n\rightarrow\infty,
\]
and
\[
\sigma_n([0,r^{2n}])\rightarrow 0 \quad {\rm as} \quad
n\rightarrow\infty.
\]
Put $A_n=((T\cc)^nT^n)^{\frac{1}{2n}}$. Then
\[
\|1_{[0,r]}(A_n)\xi_n'-\xi_n\|^2= \rho_n(]r^{2n},\infty[)\rightarrow 0 \quad {\rm as} \quad
n\rightarrow\infty,
\]
and
\[
\|1_{]r,\infty[}(A_n)(T\cc)^n\eta_n'-(T\cc)^n\eta_n'\|^2=\sigma_n([0,r^{2n}])\rightarrow 0 \quad {\rm as} \quad
n\rightarrow\infty.
\]

Since $\|\xi_n'-\xi'\|\rightarrow 0$ and $\|\eta-\eta'\|<\frac{\eps}{2}$,
we get that
\[
\|1_{[0,r]}(A_n)\xi-\xi\|<\eps,
\]
eventually as $n\rightarrow \infty$, and similarly one argues that
\[
\|1_{]r,\infty[}(A_n)\eta-\eta\|<\eps,
\]
eventually as $n\rightarrow \infty$. Thus, with $P_n=1_{[0,r]}(A_n) $ and
$P=P_{E(T,r)}=1_{[0,r]}(A)$ we have shown that
\[
\lim_{n\rightarrow\infty}\|P_n\xi-\xi\|=0
\]
for all unit vectors $\xi\in P(\CH)$ and hence for all $\xi\in
P(\CH)$. Similarly,
\[
\lim_{n\rightarrow\infty}\|(\unit-P_n)\eta - \eta\|=0
\]
for all $\eta\in P(\CH)^\bot$. It then follows from Lemma~\ref{lemma a} that
\[
{\rm SO}-\lim_{n\rightarrow\infty}1_{[0,r]}(A_n)=1_{[0,r]}(A)
\]
for all $r>0$ with $\mu_T(\partial B(0,r))=0$. Then by Lemma~\ref{B} (with
$M=\|T\|$), $A_n\rightarrow A$ as $n\rightarrow \infty$ in strong
operator topology. This proves (a).

(b) According to (a) applied to $T\cc$, the limit
\[
B={\rm SO}-\lim_{n\rightarrow\infty}((T^n)\cc T^n)^{\frac{1}{2n}}
\]
exists and is uniquely determined by
\[
1_{[0,r]}(B)=P_{E(T\cc,r)}, \qquad (r>0).
\]
Moreover, for all but countably many $r>0$,
\[
1_{[r,\infty[}(B)= 1_{]r,\infty[}(B),
\]
and
\[
F(T,r)^\bot = E(T\cc,r).
\]
Hence, for these $r$'s,
\[
1_{[r,\infty[}(B)=1- 1_{[0,r]}(B)= 1-P_{E(T\cc,r)}=P_{F(T,r)}.
\]
Since $r\mapsto 1_{[r,\infty[}(B)$ and $r\mapsto P_{F(T,r)}$ are both
SO-continuous from the left, it follows that
\[
1_{[r,\infty[}(B)=P_{F(T,r)}
\]
for all $r>0$.

Conversely, if $B'\in\CM^+$ and $1_{[r,\infty[}(B')=P_{F(T,r)}$ for all
$r>0$, then the argument above may be reversed to show that
\[
1_{[0,r]}(B')=P_{E(T\cc,r)}, \qquad(r>0),
\]
and hence, by the uniqueness in (a), $B=B'$. $\endproof$ 

\vspace{.2cm}

\begin{example} According to Theorem~\ref{existence of A and B}, for every $T\in\CM$ the sequence $\big(((T\cc)^nT^n)^{\frac{1}{2n}}\big)_{n=1}^\infty$ converges in strong operator topology. The following example due to Voiculescu (personal communication 2004) shows that this is not the case for all $T\in B(\CH)$ when ${\rm dim}\CH=+\infty$. Indeed, let $T\in B(\ell^2(\N))$ be the weighted shift given by
\[
Te_n = c_n e_{n+1}, \qquad n\in\N,
\]
where $(e_n)_{n=1}^\infty$ is the standard basis for $\ell^2(\N)$, and $(c_n)_{n=1}^\infty$ is given by
\[
c_n = \left\{ \begin{array}{llll}1 &,& 2^k\leq n<2^{k+1}, & k\quad {\rm even},\\
2 &,& 2^k\leq n<2^{k+1}, & k\quad {\rm odd}.
\end{array}
\right.
\]
Then
\[
(T\cc)^nT^n e_1 = \Big(\prod_{i=1}^nc_i\Big)^2e_1,
\]
and thus
\[
((T\cc)^nT^n)^{\frac{1}{2n}} e_1 = \Big(\prod_{i=1}^nc_i\Big)^{\frac1n}e_1.
\]
Since
\[
\limsup_{n\rightarrow\infty} \Big(\prod_{i=1}^nc_i\Big)^{\frac1n} = 2^{\frac 23},
\]
and
\[
\liminf_{n\rightarrow\infty} \Big(\prod_{i=1}^nc_i\Big)^{\frac1n} = 2^{\frac 13},
\]
it follows that $\big(((T\cc)^nT^n)^{\frac{1}{2n}}\big)_{n=1}^\infty$ is not SO--convergent in $B(\ell^2(\N))$. 
\end{example}

%% file: decomposable.tex
\section{Local spectral theory and decomposability}

A bounded operator $T$ on a Hilbert space $\CH$ is said to be {\it
  decomposable} (cf. \cite[Definition~1.1.1]{LN}), if for every open cover $\C= U\cup V$ of
the complex plane, there are $T$-invariant closed subspaces $\CH'$ and
$\CH''$ of $\CH$, such that the spectra of the restrictions of $T$ satisfy
$\sigma(T|_{\CH'})\subseteq U$ and $\sigma(T|_{\CH''})\subseteq V$, and
such that $\CH = \CH'+\CH''$.

Given $T\in B(\CH)$, a {\it spectral capacity for $T$} is a mapping $E$ from the set of closed subsets of $\C$ into
the set of all closed, $T$-invariant subspaces of $\CH$, such that
\begin{itemize}
  \item[(i)] $E(\emptyset)=\{0\}$ and $E(\C)=\CH$,
  \item[(ii)] $E(\overline U_1)+ \cdots + E(\overline U_n)=\CH$ for every
  (finite) open cover $\{U_1, \ldots, U_n\}$ of $\C$,
  \item[(iii)] $E\Big(\bigcap_{n=1}^\infty F_n\Big)= \bigcap_{n=1}^\infty
  E(F_n)$ for every countable family $(F_n)_{n=1}^\infty$ of closed subsets
  of $\C$,
  \item[(iv)] $\sigma(T|_{E(F)})\subseteq F$ for every closed subset $F$ of
  $\C$ (with the convention that $\sigma(T|_{\{0\}}) = \emptyset$).
\end{itemize}

Finally, given $T\in B(\CH)$ and $\xi\in\CH$, the {\it local resolvent
  set}, $\rho_T(\xi)$, of $T$ at $\xi$ is the union of all open subsets $U$
  of $\C$, for which there exist holomorphic vector-valued functions
  $f_U:U\rightarrow \CH$, such that $(T-\lambda\unit)f_U(\lambda)=\xi$ for
  all $\lambda\in U$. Note that according to Neumann's lemma,
  \[
  \{z\in \C\,|\,|z|>\|T\|\} \subseteq \rho_T(\xi),
  \]
  and therefore, $\sigma_T(\xi):= \C\setminus \rho_T(\xi)$, the {\it local
  spectrum} of $T$ at $\xi$, is compact. For any subset $A$ of $\C$, the
  corresponding {\it local spectral subspace} of $T$ is
  \begin{equation}\label{localspectralsubspace}
    \CH_T(A)=\{\xi\in \CH\,|\,\sigma_T(\xi)\subseteq A\}.
  \end{equation}
  It is not hard to see that $\CH_T(A)$ is $T$-hyperinvariant.

  Now, the three definitions given above, i.e. that of decomposability,
  that of a spectral capacity and that of a local spectral subspace, are
  closely related, as the following theorem indicates:

  \begin{thm}\cite[Proposition~1.2.23]{LN} \label{connections} Let $T$ be a bounded
  operator on a Hilbert space $\CH$. Then the following are equivalent:
  \begin{itemize}
    \item[(i)] $T$ is decomposable,
    \item[(ii)] $T$ has a spectral capacity,
    \item[(iii)] for every closed subset $F$ of $\C$, $\CH_T(F)$ is closed
    and
    \[
    \sigma((\unit-p_T(F))T|_{\CH_T(F)^\bot})\subseteq
    \overline{\sigma(T)\setminus F},
    \]
    where $p_T(F)$ denotes the projection onto $\CH_T(F)$. 
  \end{itemize}
  Moreover, if $T$ is decomposable, then the map $F\mapsto \CH_T(F)$ is the
  unique spectral capacity for $T$.
\end{thm}

\vspace{.2cm}

Now, if the operator $T$ appearing in our generalized version of
Theorem~\ref{embedding} is decomposable, how is the local spectral
subspace $\CH_T(B)$ related to the $T$-invariant subspace $\CK_T(B)$ for
$B\in \B(\C)$?

\begin{prop}\label{dec1} If $T$ is a decomposable operator in the II$_1$-factor $\CM$,
  then for every $B\in \B(\C)$, $\CK_T(B)= \overline{\CH_T(B)}$.
\end{prop}

\proof Since $\sigma_T(\xi)$ is compact for every $\xi\in \CH$, we have that
\[
\CH_T(B)= \bigcup_{K\subseteq B, \; K\; {\rm compact}}\CH_T(K).
\]
Moreover, by definition
\[
\CK_T(B)= \overline{\bigcup_{K\subseteq B, \; K\; {\rm compact}}\CK_T(K)}.
\]
Hence, it suffices to prove that for every compact set $K\subseteq \C$,
$\CK_T(K)=\CH_T(K)$. Since $T$ is decomposable, Theorem~\ref{connections}
implies that for every closed subset $F$ of $\C$,
\[
\sigma(T|_{\CH_T(F)})\subseteq F,
\]
and
\[
\sigma((\unit-Q_T(F))T|_{\CH_T(F)^\bot})\subseteq
\overline{\sigma(T)\setminus F},
\]
where $Q_T(F)\in W\cc(T)$ denotes the projection onto $\CH_T(F)$. In
particular,
\[
\supp(\mu_{T|_{\CH_T(F)}})\subseteq F,
\]
and
\[
\supp(\mu_{T|_{(\unit-Q_T(F))T|_{\CH_T(F)^\bot}}})\subseteq \overline{\sigma(T)\setminus F},
\]
It follows that $Q_T(F)\leq P_T(F)$, where  $P_T(F)\in W\cc(T)$ denotes
the projection onto $\CK_T(F)$. Now,
\begin{eqnarray*}
\tau(P_T(F))&=&\mu_T(F)\\
&=& \tau(Q_T(F))\mu_{T|_{\CH_T(F)}}(F) +
\tau(Q_T(F)^\bot)\mu_{(\unit-Q_T(F))T|_{\CH_T(F)^\bot}}(F),
\end{eqnarray*}
and since $F\cap \overline{\sigma(T)\setminus F}\subseteq \partial F$, we
get that
\begin{eqnarray*}
\tau(P_T(F))
&\leq & \tau(Q_T(F)) +
\tau(Q_T(F)^\bot)\mu_{(\unit-Q_T(F))T|_{\CH_T(F)^\bot}}(\partial F).
\end{eqnarray*}
Hence, if $\mu_T(\partial F)=0$, then $\tau(P_T(F))\leq \tau(Q_T(F))$, and it
follows that $P_T(F)=Q_T(F)$.

For a general closed subset $F$ of $\C$, define
\[
F_t = \Big\{z\in\C\,\Big|\,{\rm dist}(z,F)\leq \frac1t\Big\}, \qquad (t>0).
\]
Then $F_t \searrow F$ as $t\nearrow \infty$. Take $0<t_1\leq t_2\leq
\ldots$, such that $t_n\rightarrow \infty$ and such that for all $n\in\N$,
$\mu_T(\partial F_{t_n})=0$. Then, since $F\mapsto Q_T(F)$ is a spectral capacity for $T$, we have:
\begin{eqnarray*}
  Q_T(F)&=&\bigcap_{n=1}^\infty Q_T(F_{t_n})\\
  &=& \bigcap_{n=1}^\infty P_T(F_{t_n})\\
  &=& P_T(F). \qquad \qquad \endproof
\end{eqnarray*}

\vspace{.2cm}

\begin{cor}\label{dec2} If $T\in\CM$ is decomposable, then $\supp(\mu_T)=\sigma(T)$.
\end{cor}

\proof We know from \cite{Br} that $\supp(\mu_T)\subseteq \sigma(T)$ always holds. Now, let $T\in\CM$ be decomposable, and let $F=\supp(\mu_T)$. Then $\CK_T(F)=\CH$. Hence, by Proposition~\ref{dec1} and Theorem~\ref{connections}, 
\[
\CH_T(F)=\overline{\CH_T(F)}=\CK_T(F)=\CH.
\]
Therefore, by condition (iv) in the definition of a spectral capacity, 
\[
\sigma(T)=\sigma(T|_{\CH_T(F)})\subseteq F =\supp(\mu_T). \qquad \endproof
\]

\vspace{.2cm}

\begin{cor} Every II$_1$--factor $\CM$ contains a non--decomposable operator.
\end{cor}

\proof According to \cite[Example~6.6]{DH}, $\CN=\bigoplus_{k=2}^\infty B(\C^k)$ contains an operator $T$ for which $\sigma(T)=\overline\D$ and $\mu_T=\delta_0$. Since $\CN$ embeds into the hyperfinite II$_1$--factor $\CR$, and since every II$_1$--factor contains $\CR$ as a von Neumann subalgebra, $\CM$ contains a copy of $T$. According to Corollary~\ref{dec2}, $T$ is not decomposable. $\endproof$

\vspace{.2cm}

\begin{remark} By \cite{DH3}, Voiculescu's circular operator is decomposable. More generally, every DT--operator is decomposable. 
\end{remark}

%% file: appendix.tex
\section{Appendix: Proof of Theorem~\ref{partition1}}

Consider a fixed map  $f:[a,b]\rightarrow \Lp$
which is H\"older continuous with exponent $\alpha >
\frac{1-p}{p}$. That is, $f$ satisfies
\eqref{eq3-15} for some positive constant $C$. For every closed subinterval $[c,d]$ of $[a,b]$ we let $\int_c^d f(x)\,\d x$ be given by Definition~\ref{define int}. 

\begin{lemma}\label{T_x-T_y} Let $a=x_0 < x_1 <\cdots <x_{m-1}<x_m =b$ and $a=y_0 < y_1
  <\cdots <y_{n-1}<y_n =b$ be partitions of the interval $[a,b]$, and
  define
  \[
  T_x = \sum_{i=1}^m f(x_i)(x_i-x_{i-1}),
  \]
  \[
   T_y = \sum_{j=1}^n f(y_j)(y_j-y_{j-1}).
  \]

  Then with $\delta_x= \max_{1\leq i\leq m}(x_i-x_{i-1})$ and $\delta_y=
  \max_{1\leq j\leq n}(y_j-y_{j-1})$ one has that
  \begin{equation}
    \|T_x-T_y\|_p^p \leq C^p(m+n)\max\{\delta_x, \delta_y\}^{p+\alpha p}.
  \end{equation}
\end{lemma}
  
\proof Let $a=z_0 < z_1 <\cdots <z_{r-1}<z_r =b$ be the paron of the
interval $[a,b]$ containing all of the points $ x_1 , \ldots , x_{m-1}$ and
$y_1, \ldots , y_{n-1}$. For $1\leq k\leq r$, let $i(k)\in \{1,\ldots,m\}$
and  $j(k)\in \{1,\ldots,n\}$  be those indices for which
\[
[z_{k-1},z_k]\subseteq [x_{i(k)-1}, x_{i(k)}]
\]
and
\[
[z_{k-1},z_k]\subseteq [y_{j(k)-1}, y_{j(k)}].
\]

Then
\begin{equation}\label{eq3-1}
  T_x = \sum_{k=1}^rf(x_{i(k)})(z_k-z_{k-1}),
\end{equation}
and
\begin{equation}\label{eq3-2}
  T_y = \sum_{k=1}^rf(y_{j(k)})(z_k-z_{k-1}).
\end{equation}

Since $z_k\in [x_{i(k)-1}, x_{i(k)}]\cap [y_{j(k)-1}, y_{j(k)}]$, both of
the points $ x_{i(k)}$ and $y_{j(k)}$ must belong to the interval $[z_k,
z_k+\max\{\delta_x, \delta_y\}]$, and it follows that
\begin{equation}\label{eq3-3}
\|f(x_{i(k)})-f(y_{j(k)})\|_p\leq C \max\{\delta_x, \delta_y\}^\alpha.
\end{equation}

Combining \eqref{eq3-1}, \eqref{eq3-2} and \eqref{eq3-3} we find that
\[
\|T_x-T_y\|_p^p \leq  \sum_{k=1}^r (z_k-z_{k-1})^p C^p  \max\{\delta_x,
\delta_y\}^{\alpha p} \leq (m+n)C^p \max\{\delta_x,
\delta_y\}^{\alpha p + p}. \endproof
\]

\vspace{.2cm}

\begin{lemma}\label{indskud} For every $c\in (a,b)$,
  \begin{equation}
    \int_a^b f(x)\,\d x =  \int_a^c f(x)\,\d x +  \int_c^b f(x)\,\d x .
  \end{equation}
\end{lemma}

\proof Let $(S_n)_{n=1}^\infty$, $(S_n^{(1)})_{n=1}^\infty$ and
$(S_n^{(2)})_{n=1}^\infty$ be the sequences defining $\int_a^b f(x)\,\d x$,
$ \int_a^c f(x)\,\d x$ and $\int_c^b f(x)\,\d x$, respectively (cf. Definition~\ref{define int}). That is,
with
\begin{eqnarray*}
x_i &=& a + i\frac{b-a}{2^n}, \qquad (i = 0,\ldots, 2^n),\\
y_j & = & a + j\frac{c-a}{2^n}, \qquad (j = 0,\ldots, 2^n),\\
y_j & = & c + (j-2^n)\frac{b-c}{2^n}, \qquad (j = 2^n+1,\ldots, 2\cdot2^n),
\end{eqnarray*}
one has:
\begin{eqnarray*}
S_n&=&\sum_{i=1}^{2^n}f(x_i)(x_i-x_{i-1}),\\
S_n^{(1)}+S_n^{(2)}&=& \sum_{j=1}^{2\cdot 2^n}f(y_j)(y_j-y_{j-1}).
\end{eqnarray*}

Then with $\delta_x = \max_{1\leq i\leq 2^n}(x_i-x_{i-1})=\frac{b-a}{2^n}$
and $\delta_y= \max_{1\leq j\leq 2\cdot 2^n}(y_j-y_{j-1})\leq
\frac{b-a}{2^n}$ we get from Lemma~\ref{T_x-T_y} that
\begin{eqnarray*}
\|S_n - (S_n^{(1)}+S_n^{(2)})\|_p^p &\leq& 3\cdot 2^n C^p
\Big(\frac{b-a}{2^n}\Big)^{p+\alpha p}\\
& = & 3C^p(b-a)^{p+\alpha p}2^{-n(p+\alpha p-1)},
\end{eqnarray*}

and since $p+\alpha p-1>0$, we conclude that $\|S_n - (S_n^{(1)}+S_n^{(2)})\|_p\rightarrow 0$
as $n\rightarrow\infty$. $\endproof$

\vspace{.2cm}

\begin{lemma}\label{minus}
  \begin{equation}
  \int_{-b}^{-a}f(-x)\,\d x =  \int_a^b f(x)\,\d x.
  \end{equation}
\end{lemma}

\proof $\int_{-b}^{-a}f(-x)\,\d x = \lim_{n\rightarrow \infty}S_n'$, where
\begin{eqnarray*}
S_n' &=& \frac{b-a}{2^n}\sum_{k=1}^{2^n} f\Big(b-k \frac{b-a}{2^n}\Big)\\
& = & \frac{b-a}{2^n}\sum_{k=1}^{2^n} f\Big(a+(2^n-k) \frac{b-a}{2^n}\Big)\\
& = & \frac{b-a}{2^n}\sum_{l=0}^{2^n-1} f\Big(a+l \frac{b-a}{2^n}\Big),
\qquad (n\in\N_0).
\end{eqnarray*}

Hence
\[
\|S_n-S_n'\|_p = \frac{b-a}{2^n}\|f(b)-f(a)\|_p\rightarrow 0 \quad {\rm as}
\; n\rightarrow \infty,
\]

and it follows that
\[
\int_{-b}^{-a}f(-x)\,\d x = \lim_{n\rightarrow \infty}S_n' =
\lim_{n\rightarrow \infty}S_n =  \int_a^b f(x)\,\d x. \endproof
\]

\vspace{.2cm}

\begin{lemma}\label{c} For every $c\in [a,b]$,
  \begin{equation}\label{eq3-5}
    \Big\|f(c)-\frac{1}{b-a}\int_a^b f(x)\,\d x\Big\|_p^p \leq
    \frac{C^p(b-a)^{p+\alpha p}}{\alpha p +p-1}.
  \end{equation}
\end{lemma}

\proof At first we consider the case $c=b$. Taking \eqref{eq3-4} into
account we obtain:
\begin{eqnarray*}
   \Big\|(b-a)f(b)-\int_a^b f(x)\,\d x \Big\|_p^p
  &=&\lim_{n\rightarrow\infty}\|S_0-S_n\|_ p^p\\
  & \leq & \sum_{n=1}^\infty \|S_n-S_{n-1}\|_p^p\\
  & \leq &  \frac{C^p(b-a)^{p+\alpha p}}{2}\sum_{n=1}^\infty
  2^{-n(p+p\alpha -1)},
\end{eqnarray*}

where
\begin{eqnarray*}
  \sum_{n=1}^\infty 2^{-n(p+p\alpha -1)} & = & \frac{1}{2^{p+p\alpha
  -1}-1}\\
& = & \frac{1}{\e^{\log2 (p+p\alpha-1)}-1}\\
& \leq & \frac{1}{\log2 (p+p\alpha-1)}\\
& \leq & \frac{2}{p+p\alpha-1}.
\end{eqnarray*}

Hence
\[
   \Big\|(b-a)f(b)-\int_a^b f(x)\,\d x \Big\|_p^p\leq \frac{C^p(b-a)^{p+\alpha
  p}}{p+p\alpha-1},
\]

and then by Lemma~\ref{minus},
\[
  \Big\|(b-a)f(a)-\int_a^b f(x)\,\d x \Big\|_p^p= \Big\|(b-a)f(-(-a))-\int_{-b}^{-a} f(-x)\,\d x \Big\|_p^p \leq \frac{C^p(b-a)^{p+\alpha
  p}}{p+p\alpha-1}.
\]

It follows now that for arbitrary $c\in (a,b)$,
\begin{equation}
   \Big\|(c-a)f(c)-\int_a^cf(x)\,\d x \Big\|_p^p \leq
   \frac{C^p(c-a)^{p+\alpha p}}{p+p\alpha-1}
\end{equation}

and
\begin{equation}
   \Big\|(b-c)f(c)-\int_c^bf(x)\,\d x \Big\|_p^p \leq
   \frac{C^p(b-c)^{p+\alpha p}}{p+p\alpha-1},
\end{equation}

whence
\begin{eqnarray*}
  \Big\|(b-a)f(c)-\int_a^bf(x)\,\d x \Big\|_p^p & =&
  \Big\|(b-c)f(c)+(c-a)f(c)-\int_a^cf(x)\,\d x - \int_c^bf(x)\,\d
  x\Big\|_p^p \\
  & \leq & \Big\|(b-c)f(c)- \int_a^cf(x)\,\d x\Big\|_p^p +
  \Big\|(c-a)f(c)-\int_a^cf(x)\,\d x \Big\|_p^p\\
  & \leq &\frac{C^p(c-a)^{p+\alpha p}}{p+p\alpha-1} +
  \frac{C^p(b-c)^{p+\alpha p}}{p+p\alpha-1}.
\end{eqnarray*}

Since $p+\alpha p>1$,
\[
(c-a)^{p+\alpha p} +(b-c)^{p+\alpha p}\leq
(b-a)^{p+\alpha p},
\]
and \eqref{eq3-5} follows. $\endproof$ 

\vspace{.2cm}

{\it Proof of Theorem~\ref{partition1}.} According to Lemma~\ref{indskud},
\begin{eqnarray*}
  \Big\|M-  \int_a^b f(x)\,\d x\Big\|_p^p&=& \Big\|\sum_{i=1}^n
  \Big[f(t_i)(x_i-x_{i-1}) - \int_{x_{i-1}}^{x_i}f(x)\,\d x\Big]\Big\|_p^p\\
  & \leq & \sum_{i=1}^n \Big\| f(t_i)(x_i-x_{i-1}) -
  \int_{x_{i-1}}^{x_i}f(x)\,\d x\Big\|_p^p,
\end{eqnarray*}

and by Lemma~\ref{c},
\begin{eqnarray*}
\Big\| f(t_i)(x_i-x_{i-1}) -
  \int_{x_{i-1}}^{x_i}f(x)\,\d x\Big\|_p^p &\leq&
  \frac{C^p(x_i-x_{i-1})^{p+\alpha p}}{p+\alpha p -1}\\
  &=&\frac{C^p(x_i-x_{i-1})(x_i-x_{i-1})^{p+\alpha p-1}}{p+\alpha p-1}\\
  &\leq &\frac{C^p(x_i-x_{i-1})\delta(I)^{p+\alpha p-1}}{p+\alpha p-1}.
\end{eqnarray*}

Consequently,
\[
\Big\|M-  \int_a^b f(x)\,\d x\Big\|_p^p \leq \frac{C^p}{p+\alpha p
  -1}\sum_{i=1}^n (x_i-x_{i-1})\delta(I)^{p+\alpha p-1} =
\frac{C^p(b-a)\delta(I)^{p+\alpha p-1}}{p+\alpha p-1}. \endproof
\]

%% file: bibliography.tex
{\small